\documentstyle[12pt]{article}
\begin{document}
\pagestyle{empty}
\def\lsim{\raise0.3ex\hbox{$<$\kern-0.75em\raise-1.1ex\hbox{$\sim$}}}
\def\gsim{\raise0.3ex\hbox{$>$\kern-0.75em\raise-1.1ex\hbox{$\sim$}}}
\def\noi{\noindent}
\def\nn{\nonumber}
\def\bea{\begin{eqnarray}}  \def\eea{\end{eqnarray}}
\def\beq{\begin{equation}}   \def\eeq{\end{equation}}
\def\sq{\hbox {\rlap{$\sqcap$}$\sqcup$}}
\def\beeq{\begin{eqnarray}} \def\eeeq{\end{eqnarray}}
\def\R{ {\rm R \kern -.31cm I \kern .15cm}}
\def\C{ {\rm C \kern -.15cm \vrule width.5pt \kern .12cm}}
\def\Z{ {\rm Z \kern -.27cm \angle \kern .02cm}}
\def\N{ {\rm N \kern -.26cm \vrule width.4pt \kern .10cm}}
\def\1{{\rm 1\mskip-4.5mu l} }
\def\lsim{\raise0.3ex\hbox{$<$\kern-0.75em\raise-1.1ex\hbox{$\sim$}}}
\def\gsim{\raise0.3ex\hbox{$>$\kern-0.75em\raise-1.1ex\hbox{$\sim$}}}
\vbox to 2 truecm {}
\centerline{\Large \bf Long Range Scattering and Modified}
\vskip 3 truemm  
\centerline{\Large \bf Wave Operators for some Hartree Type Equations*}

\vskip 1 truecm
\centerline{\bf J. Ginibre}
\centerline{Laboratoire de Physique Th\'eorique et Hautes Energies**} 
\centerline{Universit\'e de Paris XI, B\^atiment 210,
F-91405 Orsay Cedex, France}

\vskip 5 truemm
\centerline{\bf G. Velo}
\centerline{Dipartimento di Fisica}
\centerline{Universit\`a di Bologna and INFN, Sezione di Bologna, Italy}
\vskip 1 truecm 
\begin{abstract}
We study the theory of scattering for a class of Hartree type equations with long range
interactions in space dimension $n \geq 3$, including Hartree equations with potential $V(x) =
\lambda |x|^{- \gamma}$ with $\gamma < 1$. For $1/2 < \gamma < 1$ we prove the existence of
modified wave operators with no size restriction on the data and we determine the asymptotic
behaviour in time of solutions in the range of the wave operators.
   \end{abstract}
\vskip 3 truecm
\noi AMS Classification : Primary 35P25. Secondary 35B40, 35Q40, 81U99. \par
\noi Key words : Long range scattering, modified wave operators, Hartree equation.
\vskip 1 truecm

\noindent LPTHE Orsay 98-38 \par
\noindent June 1998 \par
\vfill
\noi * Work supported in part by NATO Collaborative Research Grant 972231. \par
\noi ** Laboratoire associ\'e
au Centre National de la Recherche Scientifique - URA D0063.

\newpage
\pagestyle{plain}
\baselineskip=24 pt

\section{Introduction}
\hspace*{\parindent}
In this paper, we study the theory of scattering and more precisely the existence of modified wave
operators for a class of long range Hartree type equations
$$i \partial_t u  + {1 \over 2} \Delta u = \widetilde{g} (|u|^2)u \eqno(1.1)$$
\noi where $u$ is a complex function defined in space time ${I\hskip-1truemm R}^{n+1}$, $\Delta$
is the Laplacian in ${I\hskip-1truemm R}^n$, and
$$\widetilde{g} (|u|^2) = \lambda t^{\mu - \gamma} \ \omega^{\mu - n} \ |u|^2 \eqno(1.2)$$
\noi with $\omega = (- \Delta)^{1/2}$, $\lambda \in {I\hskip-1truemm R}$, $0 < \gamma < 1$ and $0
< \mu < n$. The operator $\omega^{\mu - n}$ can also be represented by the convolution in $x$
$$\omega^{\mu - n} \ f = C_{n, \mu} \ |x|^{- \mu} * f \eqno(1.3)$$
\noi \cite {27r}, so that (1.2) is a Hartree type interaction with potential $V(x) = C |x|^{-
\mu}$. The more standard Hartree equation corresponds to the case $\gamma = \mu$. In that case,
the nonlinearity $\widetilde{g}(|u|^2)$ becomes
$$\widetilde{g}(|u|^2) = V * |u|^2 = \lambda |x|^{-\gamma} * |u|^2 \eqno(1.4)$$
\noi with a suitable redefinition of $\lambda$. \par

A large amount of work has been devoted to the theory of scattering for the Hartree equation
(1.1) with nonlinearity (1.4) as well as with similar nonlinearities with more general potentials
$V$ [8, 9, 11-17, 19-21, 25]. As in the case of the linear Schr\"odinger equation, one must
distinguish the short range case, corresponding to $\gamma > 1$, from the long range case
corresponding to $\gamma \leq 1$. Fairly satisfactory results exist in the short range case. In
particular it is known that the (ordinary) wave operators exist in suitable function spaces for
$\gamma > 1$ \cite{25r}. Furthermore for repulsive interactions, namely for $\lambda \geq 0$, it
is known that all solutions in suitable spaces admit asymptotic states in $L^2$ for $\gamma >
1$, and that asymptotic completeness holds in suitable spaces for $\gamma > 4/3$ \cite{20r}. In
the long range case $\gamma \leq 1$, the ordinary wave operators are known not to exist in any
reasonable sense \cite{20r}, and one expects that they sould be replaced by modified wave
operators including a suitable phase in their definition, as is the case for the linear
Schr\"odinger equation. A well developed theory of long range scattering exists for the latter.
See for instance \cite{3r} for a recent and comprehensive treatment and for an extensive
bibliography. In contrast with that situation, only preliminary results are available for the
Hartree equation (and even less results for the more difficult nonlinear Schr\"odinger (NLS)
equation), all of them restricted to the case of small solutions. On the one hand, the existence
of modified wave operators has been proved in the critical case $\gamma = 1$ for small solutions
\cite{9r}. On the other hand, it has been shown recently, first in the critical case $\gamma =
1$ \cite{12r,15r} and then in the whole range $\gamma < 1$ \cite{11r,13r,14r} that the global
solutions of the Hartree equation (1.1) (1.3) with small initial data exhibit an asymptotic
behaviour as $t \to \pm \infty$ of the expected scattering type characterized by scattering
states $u_{\pm}$ and including suitable phase factors that are typical of long range scattering.
In particular, in the framework of scattering theory, the results of \cite{11r,13r,14r} are just
the property of asymptotic completeness for small data. \par

In the simpler case of short range interactions, it is a fact of experience that for the Hartree
equation as well as for other equations such that the NLS equation or the nonlinear wave (NLW)
equation, asymptotic completeness for small data can be proved together with and by the same
method as the existence of the wave operators with no size restriction on the data. The results
of \cite{11r,13r,14r} therefore suggest that the same method that has been used there to prove
asymptotic completeness for small data can be used to prove the existence of suitably modified
wave operators, again with no size restriction on the data. It is the purpose of the present
paper to explore that possibility in the case of the Hartree type equation (1.1) (1.2). \par

The special choice of nonlinearity (1.2) has been dictated by the following reasons. In one
direction, at the present preliminary stage of long range nonlinear scattering theory, it seems
more appropriate to look for the basic facts on a rather specific example rather than on a general
class of nonlinearities of the type (1.4) with $V$ replaced by some general function satisfying
suitable smoothness and decay conditions. In the opposite direction, the strictly Hartree
interaction (1.4) has the drawback that one single parameter $\gamma$ serves two unrelated
purposes. On the one hand, it characterizes the long distance behaviour of the interaction,
leading in particular to the distinction between short range and long range cases. On the other
hand, it characterizes the local regularity of the interaction, and as a consequence the local
regularity that is required for the solutions. In order to avoid that confusion, we consider
the time dependent interaction (1.2) with two independent parameters $\gamma$ and $\mu$. With
that choice it turns out that $\gamma$ again characterizes the long distance or equivalently
long time behaviour, whereas $\mu$ characterizes the local regularity. Finally the nonlinearity
(1.2) is covariant under dilations, which is an important simplifying property. \par

The main result of the present paper is the existence of modified wave operators for the
equation (1.1) (1.2), together with a description of the asymptotic behaviour in time of
solutions in the ranges of those operators, with no size restriction on the data, in suitable
spaces and for suitable values of $\gamma$ and $\mu$. The method is an extension of the energy
method used in \cite{11r,13r,14r}, and uses in particular the equations introduced in \cite{13r}
to study the asymptotic behaviour of small solutions. The spaces of initial data, namely in the
present case of asymptotic states, are Sobolev spaces of finite order similar to those used in
\cite{14r}. The parameter $\mu$ characterizing the regularity of the interactions has to
satisfy the condition $\mu \leq n - 2$ and $\mu < 2$. The condition $\mu \leq n -2$ is the
really important one and is needed for the treatment of the problem in a neighborhood of
infinity in time. It restricts the theory to space dimension $n \geq 3$. That condition and
consequently the restriction $n \geq 3$ could most probably be relaxed at the cost of using
more complicated spaces such as those used in \cite{11r,13r}. The condition $\mu < 2$ is
imposed in order to make contact with the available treatments of the equation (1.1) at finite
times \cite{8r,18r,25r}. It cannot be avoided in the attractive case $\lambda \leq 0$, whereas
it can most probably be relaxed to $\mu < {\rm Min}(4,n)$ in the repulsive case $\lambda \geq
0$. The parameter $\gamma$ characterizing the long distance/long time behaviour of the
interaction will have to satisfy $1/2 < \gamma < 1$. The critical case $\gamma = 1$ can be
treated easily by the same methods, but is excluded here in order to simplify the
exposition, since it would require different formulas involving $\ell n t$ instead of $t^{1
- \gamma}$. The case $\gamma \leq 1/2$ can be treated by the methods of this paper, but it is more
complicated and requires a more careful analysis of the asymptotic behaviour of the solutions of
(1.1). It will be deferred to a subsequent paper. \par

The construction of the modified wave operators is too complicated to allow for a more precise
statement of the results at the present stage. That construction will be described in heuristic
terms in Section 2 below. It involves in particular the study of an auxiliary system of equations
involving a new function $w$ and a phase $\varphi$ instead of the original function $u$ and the
construction of local wave operators in a neighborhood of infinity for that system. After
collecting some notation and a number of preliminary estimates in Section 3, we shall study the
local Cauchy problem at finite times for the auxiliary system in Section 4. We shall then study the
Cauchy problem at infinity and the asymptotic behaviour of solutions for the auxiliary system in
Sections 5 and 6. In particular we shall essentially construct local wave operators at infinity
for that system. We shall then come back from the auxiliary system to the original equation
(1.1) for $u$ and construct the wave operators for the latter in Section 7, where the final
result will be stated in Proposition 7.5. A more detailed description of the technical sections
3 - 7 will be given at the end of Section 2. The reader who wants to get quickly to the heart of
the matter is invited to read Section 2, to skip most of Section 3 except for the notation at
the beginning and for the definition of admissibility (Definition 3.1) to look at Proposition
4.1 and skip its proof, and to proceed to Section 5 where the main construction starts.  \par

We conclude this section with some general notation which will be used freely throughout this
paper. We denote by $\parallel \cdot \parallel_r$ the norm in $L^r \equiv L^r ({I\hskip-1truemm
R}^n)$. For $1 \leq r \leq \infty$, we define $\delta (r) = n/2 - n/r$. For any interval $I$ and
any Banach space $X$, we denote by ${\cal C}(I, X)$ the space of strongly continuous functions
from $I$ to $X$ and by $L^{\infty}(I, X)$ (resp. $L_{loc}^{\infty}(I, X))$ the space of
measurable essentially bounded (resp. locally essentially bounded) functions from $I$ to $X$. For
real numbers $a$ and $b$, we use the notation $a \vee b = {\rm Max} (a, b)$, $a \wedge b = {\rm
Min}(a,b)$ and $[a] =$ integral part of $a$. Finally if $(p \cdot q)$ is the numbering of a
double inequality, we denote by $(p\cdot qL)$ and $(p\cdot qR)$ the left hand and right hand
inequality of $(p \cdot q)$ respectively. \par

Additional notation will be given at the beginning of Section 3.

\section{Heuristics.}
\hspace*{\parindent} 
In this section, we discuss in heuristic terms the construction of the modified wave operators for
the equation (1.1), as it will be performed below in this paper. The discussion applies to any
Schr\"odinger like equation of the type (1.1) where $\widetilde{g}(|u|^2)$ is real and depends on
$u$ only through $|u|$. In addition $\widetilde{g}$ may also at this stage depend on $x$ and $t$.
For instance $\widetilde{g}$ can be independent of $u$ and depend only on $(x, t)$, thereby
leading for (1.1) to a linear Schr\"odinger equation with time dependent potential~;
$\widetilde{g}$ can be given by (1.4) or (1.2), thereby leading to the Hartree or Hartree-type
equation considered here~; $\widetilde{g}$ can be a local function of $|u|^2$ and possibly $t$
such as
$$\widetilde{g} (|u|^2) = \lambda \ t^{\mu - \gamma} \ |u|^{2 \mu /n} \eqno(2.1)$$
\noi thereby leading to a NLS equation with power nonlinearity, and with the parameters $\gamma$
and $\mu$ playing the same role as in (1.2). Whatever the case, we assume that the Cauchy
problem for (1.1) is globally well posed at finite time, namely that for any $t_0 \in {I\hskip-1truemm
R}$ and any $u_0$ in a suitable space, (1.1) has a unique solution $u$ with $u(t_0) = u_0$,
defined for all $t$ in ${I\hskip-1truemm R}$ and depending continuously on $u_0$ in suitable
norms. \par

The problem addressed by scattering theory is first of all that of classifying the possible
asymptotic behaviour of the global solutions of (1.1) by relating them to a set of model functions
with suitably chosen and preferably simple asymptotic behaviour. The first question to be
considered is then the following one. For each function $v$ of the previous set, construct a
solution $u$ of (1.1) such that $u(t)$ behaves as $v(t)$ when $t \to + \infty$, for instance in
the sense that $u(t) - v(t)$ tends to zero when $t \to + \infty$ in suitable norms. A similar
question can be asked for $t \to - \infty$. From now on we restrict our attention to the case of
positive times. \par

A natural method to attack the previous question is the following one. Let $v$ be a fixed model
function. Take $t_0 > 0$, $t_0$ large. Using the wellposedness of (1.1) for finite initial time,
define $u_{t_0}$ as the solution of the Cauchy problem for (1.1) with initial data $u_{t_0}(t_0) =
v(t_0)$ at time $t_0$. For fixed $v$, take now the limit $t_0 \to \infty$. In favourable cases,
$u_{t_0}$ will tend to a limiting solution $u_{\infty}$ of (1.1) answering the previous question.
\par

Of special interest is the case where the set of model functions $v$ is the set of solutions of
an evolution problem which is globally well posed (and preferably simpler than (1.1)). In that
case, the set of model functions $v$ can be characterized by its initial data at some prescribed
time $T$, whereas the solutions $u = u_{\infty}$ as constructed above can also be characterized
by their values $u(T)$ at time $T$. The map $v \to u_{\infty}$ classifying (part of) the
solutions of (1.1) through their asymptotic behaviour is then equivalent to the simpler map
$\Omega_+$~: $v(T) \to u(T)$ relating the values of $v$ and $u$ at time $T$. That map is the wave
operator (for positive time). \par

At this stage there is some symmetry between the original evolution for $u$ and the model
evolution for $v$. In fact let $u$ be a solution of (1.1) constructed by the previous limiting
process. Then $v$ can be recovered from $u$ as follows. Take $t_0 > 0$, $t_0$ large. Define
$v_{t_0}$ as the solution of the Cauchy problem for the model evolution with initial data
$v_{t_0}(t_0) = u(t_0)$ at time $t_0$. For fixed $u$, take now the limit $t_0 \to \infty$. In
general $v_{t_0}$ will then tend to the original function $v$ from which $u$ was constructed, and
the limiting process will provide additional information on the asymptotic behaviour of $u$. \par

In the short range case, the previous scheme is implemented by taking for $v$ the solutions of
the equation 
$$i \partial_t v + (1/2) \Delta v = 0 \eqno(2.2)$$
\noi hereafter referred to as the free equation. The generic solution of that equation is 
$v(t) = U(t) u_+$ where $U(t)$ is the unitary group
$$U(t) = \exp (i(t/2)\Delta ) \quad . \eqno(2.3)$$
\noi It is then natural to take $T = 0$. The initial data $u_+$ for $v$ is called the asymptotic
state for the solution $u = u_{\infty}$ of (1.1) constructed by the method described above, and
that method yields the ordinary wave operator $\Omega_+$ : $u_+ \to u(0)$. Note also that the
asymptotic closeness of $u$ and $v$ as $t \to \infty$ can be better expressed in terms of the
function $\widetilde{u}(t) = U(-t) u(t)$ than in terms of $u$ itself. In fact that function is
expected to tend to $u_+$ in suitable norms, whereas $u(t)$ and $v(t)$ separately tend weakly to
zero in any reasonable sense. \par

In the long range case, it is known that the previous ordinary wave operators fail to exist,
which means that the previous set of $v$'s, namely the set of solutions of the free evolution,
is badly chosen. A better set of model functions $v$ is then obtained by modifying the previous
ones by a suitable phase. That modification can be done in several ways and uses some additional
structure of $U(t)$. In fact $U(t)$ can be written as 
$$U(t) = M(t) \ D(t) \ F \ M(t) \eqno(2.4)$$
\noi where $M(t)$ denotes the operator of multiplication by the function, also denoted $M(t)$,
$$M(t) = \exp (ix^2/2t) \quad , \eqno(2.5)$$
\noi $F$ is the Fourier transform and $D(t)$ is the dilation operator defined by 
$$D(t) \ f(x) = (it)^{-n/2} \ f(x/t) \quad . \eqno(2.6)$$
\noi Let now $\varphi_0 = \varphi_0 (x, t)$ be a real function of space and time, to be chosen
later, let $z_0(x, t) = \exp (- i \varphi_0 (x, t))$, let $\varphi_0 (t)$ and $z_0(t)$ be the
operators of multiplication by the function $\varphi_0 (x, t)$ and $z_0 (x, t)$ respectively
and let $\phi_0 (t)$ and $Z_0(t)$ be the operators
$$\phi_0 (t) = \varphi_0 (- i \nabla , t) = F^*\varphi_0(t) \ F \quad ,$$
$$Z_0(t) = z_0(- i \nabla , t) = F^* z_0(t) \ F \quad . \eqno(2.7)$$
\noi In what follows, we shall sometimes omit the time dependence of the various operators when no
confusion is likely to arise. \par

Instead of the free evolution $v(t) = U(t) u_+$ one can now consider the following three modified
free evolutions~: 
$$v_1(t) = U(t) \ Z_0(t) \ u_+ = U(t) \ F^* \ z_0(t) \ w_+ \quad ,
\eqno(2.8)$$
\noi where $w_+ = Fu_+$,
 $$\hskip 2 truecm\begin{array}{lll}
 v_2(t) &= U(t) \ M^*(t) \ Z_0(t) \ u_+ = M(t) \ D(t) \ F \ Z_0(t)
\ u_+ & \\
& & \\
&= M(t) \ D(t) \ z_0 (t) \ w_+ & \\
&& \\
&= M(t) \left ( D(t) \ z_0(t) \ D^*(t) \right ) D(t) \ w_+ \quad , &\hskip 3.5 truecm (2.9)
\end{array}$$

$$\hskip 1.5 truecm \begin{array}{lll} v_3(t) &= U(t) \ M^*(t) \ Z_0(t) \ M(t) u_+ = M(t) \ D(t)
\ F \ Z_0(t) \ M(t) \ u_+ &\\ && \\ &= M(t) \ D(t) \ z_0(t) \ U^*(1/t) \ w_+ & \\ & &\\ &=
M(t) \left ( D(t) \ z_0(t) \ D^*(t) \right ) D(t) \ U^*(1/t) \ w_+ &\hskip 2 truecm (2.10)
\end{array}$$ \noi where we have used the fact that
$$F\ M(t) \ F^* = U^*(1/t) \quad . \eqno(2.11)$$
\noi Note that $D(t) z_0(t) D(t)^*$ is the operator of multiplication by $z_0(x/t, t)$ so that the
modification due to $\varphi_0$ appears only as an overall phase factor in $v_2$ and $v_3$. \par

Most of the literature on long range scattering for the linear Schr\"odinger equation makes use of
$v_1$. The function $v_2$ has been introduced in \cite{28r} and further used in \cite{4r} in the
linear case. It has been introduced independently in \cite{26r} and used in \cite{9r} \cite{26r}
in the nonlinear case. The function $v_3$ is mentioned but not really used in \cite{9r}. \par

In the short range case with $v(t) = U(t) u_+$, it was hinted that the function $\widetilde{u}(t)
= U(-t) u(t)$ was a better object of study than $u(t)$. Similarly, in the long range case, it
will be useful to introduce new functions which will be better
suited than $u$ for comparison with $v_i(t)$. Furthermore, it will be useful to express those
functions as suitable combinations of a phase factor $z_i(t) = \exp [-i\varphi_i(t)]$ and of an
amplitude $\widetilde{u}_i(t)$ in such a way that the comparison of $u$ with $v_i$ can be reduced to
the facts that asymptotically $\varphi_i(t)$ behaves as $\varphi_0(t)$ and $\widetilde{u}_i(t)$ tends
to $u_+$ or equivalently $w_i(t)$ tends to $w_+$, where $w_i(t) = F \widetilde{u}_i(t)$ for $i =
1,2,3$. This is done through the following definitions, to be compared with (2.8) (2.9) (2.10)~:
$$u(t) = U(t) \ Z_1(t) \ \widetilde{u}_1(t) = U(t) \ F^* \ z_1(t) \ w_1(t) \quad , \eqno(2.12)$$
$$u(t) = U(t) \ M^*(t) \ Z_2(t) \ \widetilde{u}_2(t) = M(t) \ D(t) \ z_2(t) \ w_2(t) \quad ,
\eqno(2.13)$$  \vskip - 0.7 truecm 

$$\hskip 1 truecm\begin{array}{lll} u(t) &= U(t) \ M^*(t) \ Z_3(t) \ M(t) \ \widetilde{u}_3(t) =
 M(t) \ D(t) \ F \ Z_3(t) \ M(t) \ \widetilde{u}_3(t) & \\ & & \\ &= M(t) \ D(t) \ z_3(t) \
U^*(1/t) \ w_3(t) \quad . &\hskip 1.7 truecm (2.14)\end{array}$$

The study of the asymptotic behaviour of small solutions of the equation (1.1) has been performed
in \cite{12r} \cite{15r} by using $(w_1, \varphi_1)$ and in \cite{11r} \cite{13r} \cite{14r} by
using essentially $(w_2 , \varphi_2)$. \par

We now explain the construction of the wave operators performed in the present paper. For technical
reasons, that construction will use the variables $(w_3 , \varphi_3)$, but we first explain it on
the example of $(w_2, \varphi_2)$ because the algebra is slightly simpler. We shall then indicate
the necessary modifications needed to switch to $(w_3 , \varphi_3)$. \par

Instead of trying to construct directly the wave operators for $u$, we first try to construct wave
operators for $(w_2, \varphi_2)$ by using the general method given at the beginning of this section.
The equation (1.1) is equivalent to the following equation for $z_2 \ w_2$ 
$$ \left ( i \partial_t + (2t^2)^{-1} \Delta - D^* \widetilde{g} D \right ) (z_2 \ w_2) = 0
\eqno(2.15)$$
\noi as can be seen by an elementary computation. Note also that
$$|u(t)| = |D(t) \ w_2(t)| \eqno(2.16)$$
\noi by (2.13), so that
$$\widetilde{g} \equiv \widetilde{g} (|u|^2) = \widetilde{g} (|D w_2|^2) \eqno(2.17)$$
\noi and $\widetilde{g}$ depends only on $w_2$ (actually only on $|w_2|$), but not on $\varphi_2$.
Expanding the derivatives in (2.15), we obtain the equivalent form \cite{13r} 
$$\left \{ i \partial_t + (2t^2)^{-1} \Delta - i(2t^2)^{-1} \ \left ( 2 \nabla \varphi_2 \cdot
\nabla + (\Delta \varphi_2 ) \right ) \right \} w_2$$
$$+ \left ( \partial_t \varphi_2 - (2t^2)^{-1} |\nabla \varphi_2|^2 - D^* \widetilde{g} D \right )
w_2 = 0 \quad . \eqno(2.18)$$

We are now in the situation of a gauge theory. The equation of evolution (2.15) and therefore also
(2.18) is invariant under the transformation $(w_2 , \varphi_2) \to (w_2 \exp (i\sigma ),
\varphi_2 + \sigma )$, where $\sigma$ is an arbitrary function of space time, and the original
gauge invariant equation (2.15) is not sufficient to provide evolution equations for the two gauge
dependent quantities $(w_2 , \varphi_2)$. At this point, we arbitrarily add a gauge condition,
which will serve as a second evolution equation, and replace (2.18) by 
$$\hskip 1.5 truecm \left \{ \begin{array}{ll} \left \{ i \partial_t + (2t^2)^{-1} \Delta -
i(2t^2)^{-1} \left ( 2 \nabla \varphi_2 \cdot \nabla + (\Delta \varphi_2) \right ) \right \} w_2 = 0
&\hskip 3.5 truecm (2.19) \\ & \\ \partial_t \varphi_2 = (2t^2)^{-1} |\nabla \varphi_2|^2 + D^*
\widetilde{g} D &\hskip 3.5 truecm (2.20) \end{array} \right .$$
\noi where the second equation, namely the gauge condition, is one of the Hamiltan-Jacobi (HJ)
equations for the classical system associated with (1.1) \cite{3r,4r}. The situation here is
similar to that occurring for the Maxwell equations where for instance one can impose the Lorentz
gauge condition $\partial_{\mu} A^{\mu} = 0$ and use it as an evolution for $A_0$ in order to
reduce the gauge freedom in the study of the Cauchy problem. \par

We have now replaced the original evolution (1.1) by the system (2.19) (2.20) and we try to study
the asymptotic behaviour of its solutions and to construct wave operators for it by the same method
that we intended to use for (1.1). In contrast with (1.1) however, the Cauchy problem for (2.19)
(2.20) cannot be expected to be globally well posed. It turns out however, and that is sufficient
for our purposes, that this problem is locally well posed in a neighborhood of infinity in time.
Roughly speaking for given initial data of arbitrary size, the Cauchy problem is well posed for
initial time $t_0$ in some interval $[T, \infty )$ for some sufficiently large $T$ depending on
the size of the data. As a consequence, we shall be able to construct only local wave operators
for (2.19) (2.20) in a neighborhood of infinity. Wave operators for (1.1) will then be obtained
from those by switching back to $u$ and using the global wellposedness of (1.1) for finite times.
\par

In order to construct the local wave operators for (2.19) (2.20), we need to choose a set of
model functions playing the role of $v$, preferably defined through a model evolution. Keeping in
mind that $u$ is represented by (2.13) and should be asymptotic to $v_2$ defined by (2.9),
preferably with $w_2(t)$ tending to $w_+$ and $\varphi_0(t)$ asymptotic to $\varphi_2 (t)$, we
define the model evolution for a pair $(w_0, \varphi_0)$ corresponding to $(w_2, \varphi_2)$ by 
$$\left \{ \begin{array}{l} \partial_t \ w_0 = 0 \\ \\ \partial_t \varphi_0 = (2t^2)^{-1}
|\nabla \varphi_0|^2 + D^*\widetilde{g} (|Dw_0|^2)D  \quad . \end{array} \right . \eqno(2.21)$$ 
\noi The first equation is immediately solved by $w_0(t) = w_+$, thereby leading to the form
(2.9) of $v_2$ where now the phase $\varphi_0$ should be a solution of the equation
$$\partial_t \varphi_0 = (2t^2)^{-1} \ |\nabla \varphi_0|^2 + D^* \widetilde{g} (|Dw_+|^2) D
\quad . \eqno(2.22)$$
\noi As in the case of the system (2.19) (2.20), the Cauchy problem for (2.22) is well posed only
in a neighborhood of infinity in time, but in general not globally in time. Although (2.22) is not
the model evolution that we shall use later on, we use it to continue the heuristic discussion.
\par

The local wave operators at infinity for the system (2.19) (2.20) as compared with (2.22) are now
constructed by the general method described at the beginning of this section. Let $\Gamma_0 =
(w_+, \varphi_0)$ be a solution of (2.22), defined in some interval $[T, \infty )$ with $T$
sufficiently large, and depending on the initial data $\varphi_0(T)$. Let $t_0 > T$ and let
$\Gamma_{t_0} = (w_{2,t_0}, \varphi_{2,t_0})$ be the solution of (2.19) (2.20) with initial data
$(w_{2,t_0}(t_0) = w_+, \varphi_{2,t_0}(t_0) = \varphi_0 (t_0)$) at time $t_0$. Under suitable
assumptions, $\Gamma_{t_0}$ will be defined in the interval $[T, \infty )$ and will converge to
a well defined limit $\Gamma_{\infty} = \Gamma$ when $t_0 \to \infty$ for fixed $\Gamma_0$. This
could provide a basis for the definition of local wave operators at infinity for the system
(2.19) (2.20), although the fact that $T$ depends on the size of the data would cause some
difficulties, but we are actually interested in wave operators for $u$, and from the previous
construction we keep only the map $\Gamma_0 \to \Gamma$. Reconstructing $u$ from $\Gamma = (w_2,
\varphi_2)$ by the use of (2.13), we obtain a map $(w_+, \varphi_0(T)) \to u$ where $u$ is a
local solution of the equation (1.1) in a neighborhood of infinity, namely defined in $[T,
\infty )$, which behaves asymptotically as $v_2$ when $t \to \infty$ in a sense that is
expressed by the relation between $(w_+, \varphi_0) = \Gamma_0$ and $(w_2, \varphi_2) = \Gamma$
at infinity, as it follows from the previous construction. We can finally complete the
construction of $u$ by solving the Cauchy problem for (1.1) with initial data $u(T)$ at time $T$
obtained from the previous step down to time $t= 1$ and define accordingly a map $(w_+, \varphi_0(T))
\to u(1)$, which is a reasonable candidate for the wave operator for (1.1). \par

The map $(w_+, \varphi_0(T)) \to u$ however is not yet satisfactory for two reasons. Firstly $u$
depends on too many data. We want $u$ to depend only on $w_+$ and not in addition on an arbitrary
initial condition for $\varphi_0$. That defect is easily remedied, as in linear long range
scattering, by imposing arbitrarily some initial condition for $\varphi_0$. For instance, given
$w_+$, one could choose in some preassigned way some $T = T(w_+)$ sufficiently large for all
subsequent constructions to be possible, and choose for instance $\varphi_0 (T(w_+)) = 0$.
Secondly, because of gauge invariance, the map $(w_+, \varphi_0(T)) \to u$ has no chance of
being injective, since different $(w_+ , \varphi_0 (T))$ can very well produce different but
gauge equivalent $(w_2, \varphi_2)$, thereby leading to the same $u$. Fixing arbitrarily the
initial condition for $\varphi_0$ certainly will improve the injectivity, but at the risk of
restricting the set of $u$ obtained by the previous construction, namely the range of the wave
operator, and making it dependent on that initial condition. We now show that in principle,
this should not happen, and that fixing the initial condition for $\varphi_0$ exactly removes
the gauge freedom and ensures the injectivity of the map $\Gamma_0 \to u$ without restricting
its range. For that purpose we have to consider in more detail the gauge covariance of the
map $\Gamma_0 \to \Gamma$. Now the HJ gauge condition (2.20) does not entirely remove the
gauge freedom in the equation (2.18) but only reduces it from that associated with an
arbitrary function of space time to that associated with an arbitrary function of space
only, for instance with some initial condition for (2.20). In fact let $(w_2, \varphi_2)$
and $(w'_2 , \varphi '_2)$ be two gauge equivalent solutions of (2.19) (2.20), namely such
that $w_2 \exp (- i \varphi_2) = w'_2 \exp (- i \varphi '_2)$. Then the difference
$\varphi_- = \varphi'_2 - \varphi_2$ satisfies the equation
$$\partial_t \ \varphi_- = (2t^2)^{-1} \ \nabla \varphi_- \cdot \nabla \varphi_+
\eqno(2.23)$$       
\noi where $\varphi_+ = \varphi'_2 + \varphi_2$. Under suitable assumptions, it follows from (2.23)
that $\varphi_-(t)$ has a well defined limit $\sigma = \lim\limits_{t \to \infty} \varphi_- (t)$ as $t
\to \infty$, whereas both $\varphi_2$ and $\varphi'_2$ grow indefinitely as $t \to \infty$.
Conversely, given a solution $(w_2, \varphi_2)$ of (2.19) (2.20) and a suitable $\sigma$, the same
equation written in the form 
$$\partial_t \varphi_- = (2t^2)^{-1} \left ( 2 \nabla \varphi_2 \cdot \nabla \varphi_- + |\nabla
\varphi_-|^2 \right ) \eqno(2.24)$$
\noi with limiting condition $\lim\limits_{t \to \infty} \varphi_-(t) = \sigma$ can be used to determine
$\varphi_-$ in a neighborhood of infinity, and thereby determine $\varphi'_2 = \varphi_2 +
\varphi_-$ and $w'_2 = w_2 \exp (i \varphi_-)$ such that $(w'_2 , \varphi '_2)$ be a solution of
(2.19) (2.20) which is gauge equivalent to $(w_2, \varphi_2)$. This argument shows that the gauge
freedom left in (2.19) (2.20) is the invariance under gauge transformations $G_{\sigma}$ which
can be parametrized by the change $\sigma$ of $\varphi_2$ at infinity. \par

In a similar way, if $(w_+, \varphi_0)$ and $(w'_+, \varphi '_0)$ are two solutions of the model
evolution (2.22) with $|w_+| = |w'_+|$, then the difference $\varphi '_0 - \varphi_0 = \varphi_-$
again satisfies the equation (2.23) now with $\varphi_+ = \varphi '_0 + \varphi_0$, and has a limit
as $t \to \infty$. Conversely the same equation rewritten in analogy with (2.24) determines
$\varphi '_0$ from $\varphi_0$ and from the value of $\varphi_-$ at some large initial time,
possibly infinity, and possibly some preassigned time $T$.  \par

The map $\Gamma_0 \to \Gamma$ as we shall construct it later will be such that $\varphi_0 (t) -
\varphi_2(t) \to 0$ and $w_2 (t) \to w_+$ as $t \to \infty$. Consequently if we define the gauge
transformation $G_{\sigma}$ on the solutions of the model evolution (2.22) by $G_{\sigma} (w_+,
\varphi_0) = (w'_+, \varphi '_0)$ with $\lim\limits_{t \to \infty} \varphi '_0 - \varphi_0 (t) = \sigma$
and $w'_+ = w_+ e^{i \sigma}$, then the map $\Gamma_0 \to \Gamma$ is gauge covariant in the
sense that the image of $G_{\sigma} \Gamma_0$ is $G_{\sigma}\Gamma$ if $\Gamma$ is the
image of $\Gamma_0$. \par

>From the previous discussion it follows that imposing an initial condition on $\varphi_0$
exactly fixes the gauge, thereby ensuring the injectivity of the map $\Gamma_0 \to u$ without
restricting its range. Actually in practice that picture is clouded by the fact that all the
constructions involved produce a small loss of regularity which prevents a complete proof of
the previous statements. On the other hand the basic construction $\Gamma_0 \to \Gamma$ also
produces a similar loss of regularity, and it turns out that the former is hidden by the
latter, so that an entirely satisfactory discussion of gauge invariance can be given at the
level of regularity of the construction $\Gamma_0 \to \Gamma$.  \par

The previous heuristic discussion was based on the system (2.19) (2.20) for $(w_2, \varphi_2)$
and the model evolution (2.22) for $\varphi_0$. For technical reasons however, we shall use
different equations. The first two reasons are rooted in the basic construction of the method,
namely the construction of the solution $\Gamma_{t_0}$ of the full evolution coinciding at
$t_0$ with a given solution $\Gamma_0$ of the model evolution, and the third one is connected
with gauge invariance. \par
(1) We shall take $w_+ \in H^k$, where $H^k$ is the standard Sobolev space, and look for $w_2$
as a continuous function of time with values in $H^{k'}$ for some $k' \leq k$. In the
construction $\Gamma_0 \to \Gamma_{t_0}$, the term $\nabla \varphi_2 \cdot \nabla w_2$ in
(2.19) produces a loss of one derivative, which seems difficult to avoid. The term $\Delta
w_2$ on the other hand produces a loss of two derivatives but that loss is easily avoided by
switching from $(w_2, \varphi_2)$ to $(w_3 , \varphi_3)$ and we shall therefore use $(w_3,
\varphi_3)$ instead of $(w_2, \varphi_2)$. The equations for $(w_3, \varphi_3)$ could be
obtained easily from the equation for $u$, but it is simpler to deduce them from the system
(2.19) (2.20). We impose $\varphi_2 = \varphi_3 = \varphi$ and $w_3 (t) = U(1/t) w_2(t) =
w(t)$, which is consistent with (2.13) (2.14). The resulting system for $(w, \varphi )$ is
then     
$$\hskip 2 truecm \left \{ \begin{array}{ll} \partial_t w = (2t^2)^{-1} U(1/t) (2 \nabla
\varphi
 \cdot \nabla + (\Delta \varphi )) U^*(1/t) w &\hskip 4.5 truecm(2.25)\\ &\\
\partial_t \varphi = (2t^2)^{-1} \ |\nabla \varphi |^2 + D^* \widetilde{g} (|DU^*(1/t) w|^2 )
D \quad . &\hskip 4.5 truecm (2.26) \end{array}  \right .$$ 
\noi Correspondingly, the model evolution for $(w_0, \varphi_0)$ will replace (2.25) by $\partial_t
w = 0$, so as to produce model functions of the type (2.10), thereby leaving for $\varphi_0$ the
equation
$$\partial_t \varphi_0 = (2t^2)^{-1} \ |\nabla \varphi_0|^2 + D^* \widetilde{g} \left ( |DU^*(1/t)
w_+|^2 \right ) D \quad . \eqno(2.27)$$
\noi Since all the estimates on $w$ will be made in spaces $H^k$ where $U(1/t)$ is unitary, the
explicit occurrence of that operator in (2.25)-(2.27) will not make any difference in those
estimates. \par

(2) In this paper we restrict our attention to the case $\gamma > 1/2$. It is well known in linear
long range scattering theory that under that condition the correcting phase $\varphi_0$ need not be a
solution of the full HJ equation (2.27) and can be chosen simply according to the Dollard
prescription, namely as a solution of the simpler equation
$$\partial_t \varphi_0 = D^* \widetilde{g} \left ( |DU^*(1/t) w_+ |^2 \right ) D \quad .
\eqno(2.28)$$

We shall partly work with (2.28) instead of (2.27) in what follows. This produces a number of
simplifications in all the questions involving only $\Gamma_0$. In particular the Cauchy problem
for (2.28) is trivially solved globally by a simple integration, thereby allowing in particular
for imposing an initial condition for $\varphi_0$ at a fixed time independent of $w_+$. \par

Whereas the term $|\nabla \varphi_0|^2$ can be omitted from (2.27) for $\gamma > 1/2$, fiddling
with the term $|\nabla \varphi|^2$ in (2.26) may not be harmless. For instance shifting that term
from (2.26) to (2.25) would produce additional restrictions on $\gamma$ and require at least
$\gamma > 2/3$ in the construction $\Gamma_0 \to \Gamma$. \par

(3) The model evolution (2.27) and its simplified version (2.28) are best suited to study the
asymptotic behaviour of the system (2.25) (2.26), which we have chosen in order to minimize the
loss of regularity. On the other hand they are ill suited for a study of gauge invariance,
which plays an important role in the reconstruction of $u$. In particular a change $w_+ \to w_+
e^{i\sigma}$ produces a nontrivial change in the $\widetilde{g}$ term in (2.27), whereas no such
change occurs in (2.22). As a consequence, we cannot avoid using also (2.22), or rather its
simplified version
$$\partial_t \varphi_0 = D^* \widetilde{g} (|Dw_+|^2) D \eqno(2.29)$$
\noi obtained as previously by omitting the $|\nabla \varphi_0|^2$ term. We shall therefore use
both (2.29), which allows for a simple discussion of gauge invariance and for a cleaner
construction of $u$, and (2.28), which yields better asymptotic approximations. \par

The last technical modification is independent of the previous choices. \par

(4) The right-hand sides of (2.25) (2.26) (2.27) (2.28) (2.29) contain $\varphi$ and $\varphi_0$
only through their gradients, and the construction of $\Gamma$ can be discussed entirely in
terms of those variables. Only for the discussion of gauge invariance and for the reconstruction
of $u$ are $\varphi_0$ and $\varphi$ themselves needed. We therefore introduce the ${I\hskip-1truemm
R}^n$ valued functions $s= \nabla \varphi$ and $s_0 = \nabla \varphi_0$ and replace the basic
equations (2.26)-(2.29) by their gradients. Using the fact that
$$\partial_i |\nabla \varphi |^2 = \sum_j (\partial_i \partial_j \varphi ) (\partial_j \varphi
) = \sum_j (\partial_j \varphi ) (\partial_j \partial_i \varphi ) = (\nabla \varphi \cdot
\nabla ) \partial_i \varphi$$
\noi we obtain 
$$\hskip 3 truecm \left \{ \begin{array}{ll} \partial_t w = (2t^2)^{-1} U(1/t) (2s \cdot \nabla
+ (\nabla \cdot s)) U^* (1/t) w \quad , &\hskip 3 truecm(2.30)\\ & \\ \partial_t s = t^{-2} (s
\cdot \nabla ) s + \nabla D^* \widetilde{g} (|DU^*(1/t) w|^2)D \quad , &\hskip 3 truecm
(2.31) \end{array} \right . $$    
\noi and either
$$\partial_t s_0 = t^{-2} (s_0 \cdot \nabla ) s_0 + \nabla D^* \widetilde{g} \left ( |DU^*(1/t)
w_+|^2 \right ) D \quad , \eqno(2.32)$$
$$\partial_t s_0 = \nabla D^* \widetilde{g} \left ( |DU^*(1/t) w_+|^2 \right ) D \quad ,
\eqno(2.33)$$
\noi or
$$\partial_t s_0 = \nabla D^* \widetilde{g} (|Dw_+|^2)D \quad . \eqno(2.34)$$
The phase $\varphi$ itself will then be recovered from $s$ through the use of (2.26) as follows.
The equation (2.31) is an Euler like equation for $s$ and implies that the vorticity $\omega =
\nabla \times s$ remains zero for all time if it is zero for some initial time $t_0$, namely with an
initial condition $s(t_0 ) = \nabla \varphi (_0)$. In fact $\omega$ satisfies the linear equation
$$\partial_t \omega = t^{-2} \left ( (s \cdot \nabla ) \omega + A \omega + \omega A^T \right )
\eqno(2.35)$$
\noi where $A$ is the matrix with entries $A_{ij} = \partial_j s_i$, and that equation implies the
result just mentioned through the Gronwall inequality for sufficiently regular $s$. It follows then
from (2.26) and (2.31) that
$$\partial_t (s - \nabla \varphi ) = t^{-2} (s \times \omega ) = 0$$
\noi and therefore $s - \nabla \varphi$ vanishes for all $t$ if it vanishes for some $t_0$. \par

We are now in a position to describe in more detail the contents of the technical part of this
paper, Sections 3-7. In Section 3, we introduce some notation, define the relevant function
spaces needed to study the system (2.30) (2.31), and we derive a number of estimates which are
used throughout the paper. In Section 4, we study the Cauchy problem at finite times for the
system (2.30) (2.31), and we prove that that problem is locally well posed (Proposition 4.1). In
Section 5, we study the local Cauchy problem at infinity for the same system, and construct
local wave operators for it as compared with the model equation (2.33). We first solve the
Cauchy problem in a neighborhood of infinity for finite but large $t_0$ (Proposition 5.1) and
derive a uniqueness result for given asymptotic behaviour (Proposition 5.2). We then prove the
existence of asymptotic states for solutions $\Gamma = (w, s)$ thereby obtained, in the
following sense~: firstly $w(t)$ has a limit $w_+$ when $t \to \infty$ (Proposition 5.3).
Secondly the solution $\Gamma_{0,t_0}$ of (2.33) which coincides with $\Gamma$ at time $t_0$
satisfies estimates uniform in $t_0$ (Proposition 5.4) and has a limit $\Gamma_0$ when $t_0
\to \infty$ which is asymptotic to $\Gamma$ (Proposition 5.5). We then turn to the converse
construction, which is that of the local wave operators at infinity. For a fixed solution
$\Gamma_0$ of (2.33), we construct a solution $\Gamma_{t_0}$ of the system (2.30) (2.31)
which coincides with $\Gamma_0$ at $t_0$ and we estimate it uniformly in $t_0$ (Proposition
5.6). We then prove that when $t_0 \to \infty$, $\Gamma_{t_0}$ has a limit $\Gamma$ which is
asymptotic to $\Gamma_0$ in the same sense as in Proposition 5.5 (Proposition 5.7). We
conclude that section with some comments on the possible use of other equations, and in
particular on the modifications required to use the more complicated equation (2.32) instead
of (2.33). In Section 6, we perform exactly the same analysis of the system (2.30) (2.31) at
infinity, now however compared with the model equation (2.34), which yields less precise
asymptotics, but which is better suited for the study of gauge invariance and the
reconstruction of $u$. Propositions 6.1-6.4 are the exact analogues of Propositions 5.4-5.7
with (2.33) replaced by (2.34) and their proofs rely to a large extent on those of the
latter. \par

Finally in Section 7 we exploit the results of Sections 5 and 6, esp. Propositions 5.7 and 6.4, to
construct the wave operators for the equation (1.1) and to describe the asymptotic behaviour of
solutions in their range. We first supplement the constructions of Sections 5 and 6 with the
appropriate definitions in order to recover $\varphi$ and $\varphi_0$ from $s$ and $s_0$, both at
finite times and in their correspondence at infinity as it follows from Propositions 5.5, 5.7, 6.2
and 6.4. We then prove that the local wave operator at infinity for the system (2.30) (2.31) as
compared with (2.34) defined through Proposition 6.4 in Definition 7.1 is gauge covariant in the
sense of Definitions 7.2 and 7.3 in the best form that can be expected with the available
regularity (Propositions 7.2 and 7.3). With the help of some information on the Cauchy problem
for (1.1) at finite time (Proposition 7.1), we then define the wave operator $\Omega$ : $u_+ \to
u$ (Definition 7.4), we prove that it is injective and has the expected range (Proposition 7.4).
We then collect all the available information on $\Omega$ and on solutions of (1.1) in its range
in Proposition 7.5, which contains the main results of this paper. Finally, by using some
information on the global Cauchy problem at finite time for (1.1) (Proposition 7.6), we define
the usual wave operator $\Omega_1$ : $u_+ \to u(1)$ (Definition 7.5). \par

A question that we leave unsettled in this paper is that of the intertwining property of the wave
operator. That property can be stated in terms of $\Omega$ as the fact that for $t$ sufficiently
large and $\tau \geq 0$,
$$\left ( \Omega (U(\tau ) u_+) \right ) (t) = \left ( \Omega (u_+) \right ) (t + \tau ) \quad .
\eqno(2.36)$$
\noi That property is an asymptotic form of time translation invariance, and unfortunately that
invariance has been severely broken by the change of variables from $u$ to $(w, \varphi )$ in
(2.14). Therefore the method used in this paper is ill suited for a study of the intertwining
property and we leave that question open here. \par

We finally remark that the basic equations (2.19) (2.20) from \cite{13r}, of which we used the
modified from (2.25) (2.26), are very similar to the equations used in \cite{7r} \cite{10r} to
study the classical limit $\hbar \to 0$ of the nonlinear Schr\"odinger equation
$$i \hbar \partial_t u = - {1 \over 2} \hbar^2 \Delta u + \widetilde{g} (|u|^2) u \quad . $$
\noi This comes as no surprise, since the latter are also obtained by separating $u$ into an
amplitude and a phase, $u = w \exp (- i \varphi / \hbar )$. Accordingly, the same energy methods
can be applied to the small $\hbar$ problem \cite{10r} and to the large time problem.

\section{Notation and preliminary estimates.}
\hspace*{\parindent} In this section we introduce some additional notation and we collect a number
of estimates which will be used throughout this paper. We first define
$$g_0 (w_1, w_2) = \lambda \ {\rm Re} \ \omega^{\mu - n} \ w_1 \bar{w}_2 \quad , \eqno(3.1)$$
$$g (w_1, w_2) = g_0 \left ( U^*(1/t) w_1, U^* (1/t) w_2 \right )  \quad . \eqno(3.2)$$
\noi In particular
$$g_{(0)}(w_1, w_1) - g_{(0)}(w_2 , w_2) = g_{(0)}(w_-, w_+)$$
\noi where $w_{\pm} = w_1 \pm w_2$. By using the definition (1.2) and (2.6), we rewrite the
nonlinearity in (2.26) as
$$D^*\widetilde{g}\left ( |DU^*(1/t)w|^2 \right ) D = t^{- \gamma} \ g(w, w) \quad . \eqno(3.3)$$
\noi The nonlinearities in (2.31)-(2.34) can be rewritten in a similar way, so that the basic
equations (2.30)-(2.34) become respectively

$$\hskip 3 truecm \left \{ \begin{array}{ll} \partial_t w = (2t^2)^{-1} \ U(1/t) (2s \cdot \nabla +
(\nabla \cdot s)) \ U^*(1/t) w &\hskip 2 truecm (2.30) \equiv(3.4) \\ & \\
\partial_t s = t^{-2}(s\cdot \nabla )s + t^{- \gamma} \ \nabla g(w, w) 
&\hskip 2 truecm(2.31)\sim(3.5) \end{array} \right .$$
$$\hskip - 3 truecm \partial_t s_0 = t^{-2}(s_0 \cdot \nabla )s_0 + t^{- \gamma} \ \nabla g
(w_+, w_+) \eqno(2.32)\sim(3.6)$$
$$\hskip - 6 truecm \partial_t s_0 = t^{-\gamma} \ \nabla g(w_+, w_+) \eqno(2.33)\sim(3.7)$$
$$\hskip - 5.5 truecm \partial_t s_0 = t^{- \gamma} \ \nabla g_0 (w_+, w_+) \quad .
\eqno(2.34)\sim(3.8)$$ We next introduce the function spaces where we shall solve the basic
equations (3.4) (3.5). We denote multi-indices by greek letters $\alpha$, $\beta$, $\cdots$,
their lengths by $|\alpha |$, $|\beta |$, $\cdots$, and nonnegative integers by $j$, $k$,
$\ell$, $\cdots$ . For any function $u$ and any function space norm $\parallel \cdot \parallel$,
we define   $$\parallel \partial^j u \parallel \ = \sum_{\alpha : |\alpha | = j} \parallel
\partial^{\alpha} u \parallel \quad .$$
We shall use Sobolev spaces of integer order $H_r^k$ defined for $1 \leq r \leq \infty$ by 
$$H_r^k = \left \{ u: \parallel u ; H_r^k \parallel \ \equiv \sum_{0 \leq j \leq k} \parallel
\partial^j u \parallel_r \ < \infty \right \}$$
\noi and the associated homogeneous spaces $\dot{H}_r^k$ with norm 
$$\parallel u ; \dot{H}_r^k\parallel \ = \ \parallel \partial^k u \parallel_r \quad . $$
\noi The subscript $r$ will be omitted if $r = 2$. \par
We first recall the well known Sobolev-Gagliardo-Nirenberg inequalities \cite{6r} \cite{22r}
\cite{23r}. \\

\noi {\bf Lemma 3.1.} {\it Let $1 \leq p, q, r \leq \infty$. Let $j$ and $k$ be nonnegative integers
with $j < k$. If $p = \infty$, assume that $k - j > n/r$. Let $\sigma$ satisfy $j/k \leq \sigma
\leq 1$ and}
$$n/p - j = (1 - \sigma )n/q + \sigma (n/r - k) \quad .$$ 
\noi {\it Then the following inequality holds for any function $u \in L^q$ :}
$$\parallel \partial^j u \parallel_p \ \leq C \parallel u\parallel_q^{1 - \sigma } \ \parallel
\partial^k u \parallel_r^{\sigma}$$
\noi {\it except for $p < \infty$, $j = 0$ and $q = \infty$. In the latter case, for any function
$u \in L^{\infty}$ there exists a constant $c$ depending on $u$ such that} 
$$\parallel u - c \parallel_p \ \leq C \parallel u - c\parallel_{\infty}^{1 - \sigma} \ \parallel
\partial^k u \parallel_r^{\sigma} \quad .$$
\noi {\bf Remark 3.1.} The statement as just given may differ from the usual ones (see for instance
\cite{6r}) by unnecessarily excluding a few trivial cases. \\

We shall use extensively the following spaces. Let $\ell_0 = [n/2]$ and define $r_0$ by $\delta
(r_0) = \ell_0$ so that $r_0 = 2n$ for $n$ odd and $r_0 = \infty$ for $n$ even. Let $k$ and $\ell$ be nonnegative integers with $\ell \geq \ell_0-1$. Let $I \subset {I\hskip-1truemm
R}^+$ be an interval. We shall look for $w$ as a complex valued function in spaces
$L_{loc}^{\infty}(I, H^k)$ or ${\cal C}(I, H^k)$ and for $s$ as an ${I\hskip-1truemm R}^n$ or
$\C^n$ vector valued function in spaces $L_{loc}^{\infty}(I, X^{\ell})$ or ${\cal C}(I, X^{\ell})$
where
$$X^{\ell} = L^{r_0} \cap \dot{H}^{\ell_0} \cap \dot{H}^{\ell + 1} \quad . \eqno(3.9)$$

For $n$ odd, it follows from Lemma 3.1 that for $s \in X^{\ell}$, 
$$\parallel s \parallel_{r_0} \ \equiv \ \parallel s \parallel_{2n} \ \leq C \parallel s  ;
\dot{H}^{\ell_0} \parallel \ \equiv \ C \parallel s ; \dot{H}^{(n-1)/2} \parallel$$
\noi so that the $L^{r_0}$ norm will not need to be estimated separately. The space $L^{r_0}$ has
been included in the definition of $X^{\ell}$ only in order to make Lemma 3.1 applicable by
eliminating arbitrary polynomials of degree $\ell_0 - 1$ which are not seen by the other norms.
For $\ell \geq \ell_0$, the inclusion $X^{\ell} \subset L^{\infty}$ holds. In fact, by Lemma 3.1
again
$$\parallel s \parallel_{\infty} \ \leq C \left ( \parallel s \parallel_{r_0} \ \parallel s;
\dot{H}^{\ell_0 + 1}\parallel \right )^{1/2} \leq C \left ( \parallel s; \dot{H}^{\ell_0}
\parallel \ \parallel s; \dot{H}^{\ell_0 + 1} \parallel \right )^{1/2} \eqno(3.10)$$
\noi so that also the $L^{\infty}$ norm will not need to be estimated either. \par

For $n$ even, the $L^{r_0}$ norm, namely the $L^{\infty}$ norm, is not controlled by the
$\dot{H}^{\ell_0}$ norm, namely by the $\dot{H}^{n/2}$ norm, and will require separate estimates.
\par

The spaces $X^{\ell}$ obviously satisfy the embedding $X^{\ell '} \subset X^{\ell}$ for $\ell '
\geq \ell$. \par

Because of the presence of $L^{r_0}$ in the definition of $X^{\ell}$, one can replace
$\dot{H}^{\ell_0}$ in that definition by
$$K^{\ell_0} = \left \{ u:u \in \dot{H}^{\ell_0} \ \hbox{and} \ <x>^{-(n+1)/2} u \in L^2 \right \}$$
\noi which is a Hilbert space, and similarly one can replace $\dot{H}^{\ell + 1}$ by $K^{\ell +
1}$. As a consequence, $X^{\ell}$ is the intersection of the duals of (compatible) Banach spaces
and is therefore itself the dual of a Banach space. \par

We shall use systematically the short hand notation
$$|w|_k = \ \parallel w; H^k \parallel \qquad , \quad |s|_{\ell} = \ \parallel s;
X^{\ell}\parallel \eqno(3.11)$$
\noi and the meaning of the symbol $|a|_b$ will always be made unambiguous by the fact that the
pair $(a, b)$ contains either the pair $(w, k)$ or the pair $(s, \ell )$. \par

We shall need estimates of the solutions (in the sense of distributions) of transport diffusion
equations of the form
$$\partial_t u = \eta \Delta u + \nabla \cdot (uv) + h \eqno(3.12)$$
\noi where $u$ and $h$ are complex valued functions and $v$ a $\C^n$ vector valued function
defined in space time. \\ 

\noi {\bf Lemma 3.2.} {\it Let $2 \leq p \leq \infty$, let $I$ be an open interval. Let $u$, $h
\in L_{loc}^p(I, L^p)$ and $v \in L_{loc}^{\infty}(I, L^{\infty})$ with $\partial v \in
L_{loc}^{\infty}(I, L^{\infty})$ satisfy the equation (3.12) for some $\eta \geq 0$, for all $t
\in I$. Then for almost all $t_1$, $t_2 \in I$, with $t_1 \leq t_2$, the following estimate
holds~:} $$\parallel u(t_2)\parallel_p \ \leq \ \parallel u(t_1)\parallel_p \ + \int_{t_1}^{t_2}
dt' \Big \{ p^{-1} \parallel \nabla \cdot v(t')\parallel_{\infty} \ \parallel u(t') \parallel_p$$
$$+ C_0 \parallel \partial v(t')\parallel_{\infty} \ \parallel u(t')\parallel_p \ + \ \parallel
h(t') \parallel_p \Big \} \eqno(3.13)$$
\noi {\it for some absolute constant $C_0$. If $\eta = 0$, a similar estimate holds for $t_1 \geq
t_2$.} \\

\noi {\bf Proof.} The formal computation leading to (3.13) is given in Appendix A. The actual
proof is obtained by following the methods in \cite{5r}. \par \nobreak
\hfill $\sq$ \\

Note that the estimate (3.13) is independent of $\eta$. In subsequent applications of Lemma 3.2,
we shall for brevity state the results thereby obtained in the shorter differential form
corresponding to
$$\partial_t \parallel u \parallel_p \ \leq p^{-1} \parallel \nabla \cdot v(t) \parallel_{\infty}
\ \parallel u(t) \parallel_p \ + C_0 \parallel \partial v(t) \parallel_{\infty} \ \parallel
u(t)\parallel_p \ + \ \parallel h(t)\parallel_p \quad . \eqno(3.14)$$

We next give some preliminary estimates. They involve functions called $s$ and $w$ since that is
suggestive of subsequent applications, but at the present stage it is irrelevant whether those
functions are real or complex, scalar or vector valued. \\

\noi {\bf Lemma 3.3.} {\it Let $\alpha$ and $\beta$ be multi-indices with $\beta \leq \alpha$ and
let $\ell = |\alpha |$. Then the following estimate holds~:}
$$\parallel \partial^{\beta} s_1 \ \partial^{\alpha - \beta} s_2 \parallel_2 \ \leq \ C \parallel s_1
; L^{\infty} \cap \dot{H}^{\ell} \parallel \ \parallel s_2; L^{\infty} \cap \dot{H}^{\ell}
\parallel \quad . \eqno(3.15)$$  

\noi {\bf Proof.} We estimate by the H\"older inequality
$$\parallel \partial^{\beta} s_1 \ \partial^{\alpha - \beta} s_2 \parallel_2 \ \leq \ \parallel
\partial^{\beta} s_1 \parallel_{r_1} \ \parallel \partial^{\alpha - \beta} s_2 \parallel_{r_2}$$
\noi with $1/r_1 + 1/r_2 = 1/2$. \par

For $\beta = \alpha$, we take $r_1 = 2$, $r_2 = \infty$. For $\beta = 0$, we take $r_1 = \infty$
and $r_2 = 2$. For $0 < \beta < \alpha$, we apply Lemma 3.1 and obtain
$$\cdots \leq C \parallel s_1 \parallel_{\infty}^{1 - \sigma_1} \ \parallel s_1 ; \dot{H}^{\ell}
\parallel^{\sigma_1} \ \parallel s_2 \parallel_{\infty}^{1 - \sigma_2} \ \parallel s_2;
\dot{H}^{\ell}\parallel^{\sigma_2} \eqno(3.16)$$
\noi with
$$\sigma_1(n/2 - \ell ) = n/r_1 - |\beta| \quad ,$$
$$\sigma_2 (n/2 - \ell ) = n/r_2 - \ell + |\beta | \quad ,$$
$$|\beta |/\ell \leq \sigma_1 \leq 1 \quad , \quad 1 - |\beta |/\ell \leq \sigma_2 \leq 1 \quad .$$
\noi The equalities imply $\sigma_1 + \sigma_2 = 1$, while the inequalities imply $\sigma_1 +
\sigma_2 \geq 1$. As a consequence, all of them are satisfied for the unique choice
$$2/r_1 = \sigma_1 = |\beta|/\ell \quad , \quad 2/r_2 = \sigma_2 = 1 - |\beta |/\ell \quad .$$
\noi Then (3.15) follows from (3.16). \par \nobreak
\hfill $\sq$ \\

\noi {\bf Lemma 3.4.} {\it Let $\alpha$ and $\beta$ be multi-indices with $\beta \leq \alpha$,
 $|\alpha | \leq k$, $|\beta | \leq \ell$ and $\ell > n/2$. Then the following estimate
holds~:} $$\parallel \partial^{\beta} s \ \partial^{\alpha - \beta} w\parallel_2 \ \leq C
\parallel s; L^{\infty} \cap \dot{H}^{\ell} \parallel \ |w|_k \quad . \eqno(3.17)$$

\noi {\bf Proof.} It is sufficient to consider the case $\alpha = \beta$. The general case with
$\alpha \geq \beta$ will then follow therefrom by replacing $w$ by $\partial^{\alpha - \beta} w$.
We estimate by the H\"older inequality
$$\parallel \partial^{\beta} s \ w \parallel_2 \ \leq \ \parallel \partial^{\beta} s\parallel_{r_1}
\ \parallel w \parallel_{r_2}$$
\noi with $1/r_1 + 1/r_2 = 1/2$. \par

For $\beta = 0$, we take $r_1 = \infty$, $r_2 = 2$. For $|\beta | = \ell$, we take $r_1 = 2$, $r_2
= \infty$ and use the fact that in that case
$$\parallel w \parallel_{\infty} \ \leq \ C \parallel w; H^{|\beta |} \parallel \ = C \parallel
w; H^{\ell} \parallel$$
\noi since $\ell > n/2$. For $0 < |\beta | < \ell$, we apply Lemma 3.1 and obtain 
$$\cdots \leq C \parallel s \parallel_{\infty}^{1 - \sigma} \ \parallel
\partial^{\ell}s\parallel_2^{\sigma} \ \parallel w \parallel_{r_2} \eqno(3.18)$$
\noi with $\sigma (n/2 - \ell ) = n/r_1 - |\beta |$ and $| \beta |/\ell \leq \sigma \leq 1$. We
then choose $\sigma = |\beta |/\ell$ so that $2/r_2 = 1 - |\beta | / \ell$ and $\delta_2 = \delta
(r_2) = n |\beta |/2 \ell < |\beta |$ since $\ell > n/2$, so that by Lemma 3.1 again $\parallel w
\parallel_{r_2} \leq$\break \noindent $\parallel w; H^{|\beta |}\parallel$, which together with
(3.18), implies (3.17) in the special case $\alpha = \beta$. \par \nobreak
\hfill $\sq$ \\

\noi {\bf Lemma 3.5.} {\it Let $\varphi$ be a real function wtih $s = \nabla \varphi \in L^{\infty}
\cap \dot{H}^{\ell}$ for some $\ell > n/2$ and let $k \leq \ell + 1$. Then the following estimate
holds~:}
$$\left | e^{-i \varphi} w \right |_k \leq C \left ( 1 + \parallel \nabla \varphi ; L^{\infty} \cap
\dot{H}^{\ell} ) \parallel \right )^k \ |w|_k \quad . \eqno(3.19)$$
\noi {\it Let in addition $\varphi \in L^{\infty}$. Then the following estimate holds~:} 
$$\left | (e^{- i \varphi} - 1)w \right |_k \leq C \left ( \parallel \varphi \parallel_{\infty} \ +
\ \parallel \nabla \varphi ; L^{\infty} \cap \dot{H}^{\ell}\parallel \left ( 1 + \parallel \nabla
\varphi ; L^{\infty} \cap \dot{H}^{\ell} \parallel \right )^{k-1} \right ) |w|_k  \quad
. \eqno(3.20)$$

\noi {\bf Proof.} For any multi-indices $\alpha$, $\beta$ with $\beta \leq \alpha$ and $1 \leq
|\beta | \leq |\alpha | \leq k$, one has to estimate the $L^2$ norm of 
$$\partial^{\beta} \left ( e^{-i \varphi} \right ) \partial^{\alpha - \beta} w = e^{-i \varphi}
\sum_{\{\beta _i \}} (-i)^m \ C_{\{\beta_i \}} \left ( \prod_{1 \leq i \leq m} \partial^{\beta_i}
\varphi \right ) \partial^{\alpha - \beta} w \eqno(3.21)$$
\noi where the sum runs over all possible decompositions $\beta = \beta_1 + \cdots + \beta_m$ of
$\beta$ as the sum of $m \geq 1$ multi-indices. We have omitted the terms with $\beta = 0$ which
trivially satisfy the required estimate with obvious assumptions on $\varphi$. We estimate each
term in the RHS of (3.21) by Lemma 3.4 applied with $s = \nabla \varphi$, $\beta \leq \beta_1$,
$|\beta | = |\beta_1 | - 1$, $\alpha = \beta$ and $w$ replaced by the product of the last $m$
factors in (3.21). We obtain
$$\parallel \partial^{\beta} (e^{-i \varphi}) \partial^{\alpha - \beta} w \parallel_2 \ \leq C
\parallel \nabla \varphi ; L^{\infty} \cap \dot{H}^{\ell} \parallel 
\sum_{\{\beta '_i, \alpha '\}} \left | \left | \left ( \prod_{2 \leq i \leq m} \partial^{\beta
'_i} \varphi \right ) \partial^{\alpha '} w \right | \right |_2$$
\noi where $\beta '_i$ $(2 \leq i \leq m)$ and $\alpha '$ are multi-indices obtained by
distributing at most $|\beta_1 | - 1$ derivatives on $\beta_2, \cdots , \beta_m, \alpha - \beta$.
In particular, in the nontrivial case $m \geq 2$, one has $|\beta '_i| \leq |\beta_i| +
|\beta_1| - 1 \leq |\beta |-1 \leq |\alpha | - 1 \leq k - 1 \leq \ell$. One can then iterate the
process, thereby extracting the $m$-th power of $\parallel \nabla \varphi;L^{\infty} \cap
\dot{H}^{\ell}\parallel$ and obtaining for each $m \geq 1$ a contribution
$$\parallel \nabla \varphi ;L^{\infty} \cap \dot{H}^{\ell}\parallel^m \ |w|_{k-m} \quad .$$
\noi Taking the sum over $m$ and adding the contribution of the term $m = 0$ yields (3.19)
(3.20). \par \nobreak 
\hfill $\sq$ \\

The next estimate will be needed to estimate the nonlinearity $\widetilde{g}(|u|^2)$. \\

\noi {\bf Lemma 3.6.} {\it Let $n \geq 2$, $0 < \mu < n$, let $\ell$, $k_1$ and $k_2$ be
nonnegative integers with $k_1 \leq k_2$, and}
$$n/2 < \ell + \mu \leq (n/2 + k_1 + k_2) \wedge (n +k_1) \quad , \eqno(3.22)$$ 
\noi {\it and in addition $k_2 > n/2$ if $\ell + \mu = n + k_1$. Then the following estimate
holds~:} $$\parallel \partial^{\ell} \ \omega^{\mu - n} (w_1w_2)\parallel_2 \ \leq C |w_1|_{k_1}
\ |w_2|_{k_2} \quad . \eqno(3.23)$$

\noi {\bf Proof.} Let $m = [n - \mu ] \wedge \ell$ and define $r$ by 
$$\delta \equiv \delta (r) = n - \mu - m$$
\noi so that $0 \leq \delta (r) < n/2$ if $\ell + \mu \leq n$ by (3.22L), while $0 \leq \delta (r)
< 1$ if $\ell + \mu \geq n$. By the Hardy-Littlewood Sobolev (HLS) inequality (\cite{22r} p. 116)
if $m < n - \mu$ and by inspection if $m = n - \mu$, we estimate
$$\parallel \partial^{\ell} \ \omega^{\mu - n} (w_1w_2) \parallel_2 \ = \ \parallel \partial^m \
\omega^{\mu - n} \ \partial^{\ell - m} (w_1w_2)\parallel_2$$
$$\leq C \parallel \omega^{- \delta} \ \partial^{\ell - m} (w_1w_2)\parallel_2 \ \leq C \parallel
\partial^{\ell - m} (w_1 w_2) \parallel_{\bar{r}} \eqno(3.24)$$
\noi where $1/r + 1/\bar{r} = 1$ and $\bar{r} > 1$. By the Leibnitz formula and the H\"older
inequality, we continue (3.24) as
$$\cdots \leq C \sum_{0 \leq j \leq \ell - m} \parallel \partial^j w_1 \parallel_{r_1} \
\parallel \partial^{\ell - m - j} w_2 \parallel_{r_2} \eqno(3.25)$$
\noi where $\delta_i \equiv \delta (r_i)$, $i = 1, 2$, have to satisfy
$$\delta_1 + \delta_2 = n/2 - \delta = \mu + m - n/2 \quad . \eqno(3.26)$$
\noi One can then continue (3.25) as
$$\cdots \leq C |w_1|_{k_1} \ |w_2|_{k_2} \eqno(3.27)$$
\noi provided $\ell - m \leq k_1$ and provided for each $j$ one can choose $\delta_1$ and
$\delta_2$ satisfying (3.26) and
$$\left \{ \begin{array}{l} 0 \leq \delta_1 \leq k_1 - j \\ \\ 0 \leq \delta_2
\leq k_2 - (\ell - m - j) \end{array}\right . \eqno(3.28)$$
\noi with the RHS inequality being strict if the corresponding $\delta_i$ is equal to $n/2$. The
condition $\ell - m \leq k_1$ is equivalent to $\ell \leq k_1 + [n - \mu ]$ and therefore to $\ell
\leq k_1 + n - \mu$ since $\ell$ and $k_1$ are integers, and follows therefore from (3.22). The
compatibility of (3.26) (3.28) in $\delta_i$ reduces to 
$$\mu + m - n/2 \leq k_1 + k_2 - \ell + m \quad ,$$
\noi again a consequence of (3.22). Finally the only possible exceptional case corresponds to
$\delta = 0$, $j = \ell - m = k_1$, namely $\ell + \mu = n + k_1$, and requires $k_2 >
n/2$.\par \nobreak
\hfill $\sq$ \\

We shall use Lemma 3.6 through the following corollary. We recall that $\ell_0 = [n/2]$. \\

\noi {\bf Corollary 3.1.} {\it Let $0 < \mu < n$ and let $k$, $\ell$ be nonnegative integers
satisfying $\ell > n/2$ and}
$$\ell + 2 + \mu \leq (n/2 + 2k) \wedge (n + k) \eqno(3.29)$$
\noi {\it and in addition $k > n/2$ if $\ell + 2 + \mu = n + k$. Then the following estimates
hold~:}
$$\parallel \partial^{\ell_0 + 1} \ \omega^{\mu - n} \ w_1 w_2 \parallel_2 \ + \ \parallel
\partial^{\ell} \ \omega^{\mu - n} \ w_1 w_2 \parallel_2 \ \leq C |w_1|_{k-2} \ |w_2|_k
\quad , \eqno(3.30)$$
$$\parallel \partial^{\ell + 1} \ \omega^{\mu - n} \ w_1 w_2 \parallel_2 \ \leq C |w_1|_{k-1} \
|w_2|_k \quad , \eqno(3.31)$$
$$\parallel \partial^{\ell + 2} \ \omega^{\mu - n} \ w_1 w_2 \parallel_2 \ \leq C |w_1|_{k} \
|w_2|_k \quad , \eqno(3.32)$$
$$\parallel \partial^{\ell + 3} \ \omega^{\mu - n} \ |w|^2 \parallel_2 \ \leq C |w|_{k} \
|w|_{k+1} \quad , \eqno(3.33)$$
$$\parallel \partial \omega^{\mu - n} \ w_1 w_2 \parallel_{\infty} \ \leq C |w_1|_{k-1} \
|w_2|_k \quad , \eqno(3.34)$$
$$\parallel \omega^{\mu - n} \ w_1 w_2 \parallel_{\infty} \ \leq C |w_1|_{k-1} \
|w_2|_{k-1} \quad . \eqno(3.35)$$
\noi {\it If $n$ is even, assume in addition that the inequality (3.29) with $\ell$ replaced by
$n/2+1$ (which is the lowest allowed value) is strict, namely}
$$n/2 + 3 + \mu < (n/2 + 2k) \wedge (n + k) \quad .\eqno(3.36)$$
\noi {\it Then the following estimate holds :}
$$\parallel \partial \omega^{\mu - n} \ w_1 w_2 \parallel_{\infty} \ \leq C |w_1|_{k-2} \
|w_2|_k \quad . \eqno(3.37)$$

\noi {\bf Proof.} The estimates (3.30)-(3.32) are direct applications of Lemma 3.6 with $\ell$
replaced res\-pec\-ti\-ve\-ly by $\ell_0 + 1$ and $\ell$, by $\ell + 1$ and by $\ell + 2$, while
$(k_1, k_2)$ are replaced respectively by $(k - 2, k)$, by $(k - 1, k)$ and by $(k, k)$. The
condition (3.22L) follows from $n/2 < \ell_0 + 1 \leq \ell$, while the condition (3.22R) follows
from, actually reduces to (3.29) in all three cases. The estimate (3.33) follows from (3.32)
applied to $(\bar{w} , \partial w)$ and to $(w, \partial \bar{w})$. \par

In order to prove (3.34) (3.35) and (3.37), we note that by (1.3) $\omega^{\mu - n} w_1 w_2$ and
$\partial \omega^{\mu - n} w_1 w_2$ belong to $L^{r_+} + L^{r_-}$ for $w_1$, $w_2 \in H^1$ with
$n/r_{\pm} = \mu \pm \eta$. One can then estimate by Lemma 3.1
$$\parallel \omega^{\mu - n} \ w_1 w_2 \parallel_{\infty} \ \leq C \parallel \omega^{\mu - n}
\ w_1 w_2 \parallel_{n/\varepsilon}^{1/2} \ \parallel \partial \omega^{\mu - n} \ w_1
w_2\parallel_{n/(1 - \varepsilon )}^{1/2}$$
\noi which by the Hardy-Littlewood-Sobolev inequality implies 
$$\parallel\omega^{\mu - n} \ w_1 w_2 \parallel_{\infty} \ \leq C \parallel \omega^{\mu - n/2 -
\varepsilon} \ w_1 w_2 \parallel_{2}^{1/2} \ \parallel \omega^{\mu - n/2 + \varepsilon} \ w_1
w_2\parallel_{2}^{1/2} \eqno(3.38)$$
\noi for any $\varepsilon$, $0 < \varepsilon \leq 1$. Similarly one can write the estimate 
$$\parallel \partial \omega^{\mu - n} \ w_1 w_2 \parallel_{\infty} \ \leq C \parallel
\partial \omega^{\mu - n/2 - \varepsilon} \ w_1 w_2 \parallel_{2}^{1/2} \ \parallel
\partial \omega^{\mu - n/2 + \varepsilon} \ w_1 w_2\parallel_{2}^{1/2} \quad . \eqno(3.39)$$
\noi We rewrite each one of the four norms in the RHS of (3.38) and (3.39) as $\parallel
\partial^{m_1} \omega^{\mu - m_0 - n/2 \pm \varepsilon} w_1w_2\parallel_2$ with $m_1 = m_0$ in
(3.38) and $m_1 = 1 + m_0$ in (3.39), the values of $m_0$ to be chosen later. The proof is now
achieved by repeated application of Lemma 3.6 after identification of the relevant variables,
(3.35) being implied by (3.38) and (3.34) and (3.37) by (3.39). The quantity $\varepsilon$ will 
always be taken sufficiently small. \par

For $n$ even we choose $m_0 = n/2$. We apply (3.23) with $(\ell , \mu , k_1, k_2)$ replaced by
$(n/2, \mu \pm \varepsilon,$\break \noindent $k-1, k-1)$ to (3.38) and by $(n/2 + 1, \mu \pm
\varepsilon , k - 1, k)$ to (3.39), thereby obtaining (3.35) and (3.34) respectively. Similarly
we prove (3.37) by replacing $(\ell , \mu , k_1, k_2)$ by $(n/2 + 1 , \mu \pm
\varepsilon , k - 2 , k)$. The condition (3.22L) is obviously satisfied while (3.22R) follows
from (3.29) in the first two cases and from (3.36) in the last one. \par

For $n$ odd we choose $m_0 = (n+1)/2$ if $1/2 + \varepsilon < \mu < n$ and $m_0 = (n -
1)/2$ if $0 < \mu \leq 1/2 + \varepsilon$. In the first case we apply (3.23) with $(\ell
, \mu, k_1, k_2)$ replaced by $((n+1)/2, \mu - 1/2 \pm \varepsilon, k - 1, k-1)$ to (3.38) and by
$((n+3)/2, \mu - 1/2 \pm \varepsilon , k - 1, k)$ to (3.39), thereby obtaining (3.35) and (3.34)
respectively. In the second case we apply (3.23) with $(\ell , \mu , k_1, k_2)$ replaced by
$((n-1)/2, \mu + 1/2 \pm \varepsilon , k - 1 , k-1)$ to (3.38) and by $((n+1)/2, \mu + 1/2 \pm
\varepsilon , k - 1, k)$ to (3.39), thereby obtaining (3.35) and (3.34) respectively. In both
cases (3.22L) is obviously satisfied while (3.22R) follows from (3.29). \par \nobreak
\hfill $\sq$ \\

We now introduce the following definition \\

\noi {\bf Definition 3.1.} A pair of nonnegative integers $(k, \ell )$ will be called admissible
if it satisfies $k \leq \ell$, $\ell > n/2$ and (3.29) and in addition $k > n/2$ if $\ell + 2 +
\mu = n + k$, and (3.36) if $n$ is even. \\

Admissible pairs exist only if $\mu \leq n - 2$. For $\mu = n - 2$, admissible pairs are pairs
$(k, \ell )$ such that $k = \ell > n/2$. Admissible pairs always have $k > 1 + \mu /2$ and
therefore $k \geq 2$. If $(k, \ell )$ is an admissible pair, so is $(k + j, \ell + j)$ for any
positive integer $j$. For $n = 3$, $\mu = 1$, the pair $(k , \ell ) = (2, 2)$ is admissible. \par

We next derive a number of a priori estimates for solutions of the basic equations (3.4) (3.5), to
be understood in the distribution sense. In the rest of this section, we assume $n \geq 3$, $0 <
\mu \leq n - 2$, and we denote by $(k, \ell )$ an admissible pair. For solutions of (3.4) (3.5), we
shall derive a number of estimates similar to (3.13) in Lemma 3.2, and we shall state them in the
shorter differential form corresponding to (3.14). \\

\noi {\bf Lemma 3.7.} {\it Let $0 < \mu \leq n - 2$ and let $(k, \ell )$ be an admissible pair.
Let $I \subset {I\hskip-1truemm R}^+$ be an open interval, and let $(w, s)$ be a solution of
(3.4) (3.5) in $L_{loc}^{\infty}(I, H^k \oplus X^{\ell})$. Then $(w, s)$ satisfies the following
estimates~:}
$$\Big | \partial_t |w|_k \Big | \leq C \ t^{-2} \parallel \partial s ; L^{\infty} \cap
\dot{H}^{\ell} \parallel |w|_k \leq C \ t^{-2} |s|_{\ell} \ |w|_k \quad , \eqno(3.40)$$
$$\Big | \partial_t \parallel s; \dot{H}^{\ell + 1} \parallel \Big | \ \leq C \ t^{-2} \parallel
\partial s ; L^{\infty} \cap \dot{H}^{\ell} \parallel^2 \ + C \ t^{- \gamma} |w|_k^2 \quad ,
\eqno(3.41)$$ $$\Big | \partial_t \parallel s; \dot{H}^{\ell_0} \parallel \Big | \ \leq C \
t^{-2} \parallel s ; L^{\infty} \cap \dot{H}^{\ell_0} \parallel \ \parallel \partial s ;
L^{\infty} \cap \dot{H}^{\ell_0} \parallel \ + C \ t^{- \gamma} |w|_{k-2} \ |w|_k \quad ,
\eqno(3.42)$$ $$\Big | \partial_t \parallel s \parallel_{\infty} \Big | \ \leq C \ t^{-2}
\parallel s  \parallel_{\infty} \ \parallel \partial s \parallel_{\infty} \ + \ C \ t^{- \gamma}
|w|_{k-1} \ |w|_k \quad , \eqno(3.43)$$ $$\Big | \partial_t |s|_{\ell} \Big | \leq C \ t^{-2}
|s|_{\ell}^2 + C \ t^{- \gamma} |w|_k^2 \quad . \eqno(3.44)$$

\noi {\bf Proof.} (3.40L). For any multi-index $\alpha$ with $|\alpha | \leq k$, we obtain from
(3.4) 
$$\partial_t \ \partial^{\alpha} w = (2t^2)^{-1} U(1/t) \sum_{\beta \leq \alpha }
\left ( 2 \partial^{\beta} s \cdot \nabla + (\partial^{\beta} \nabla \cdot s) \right )
U^*(1/t) \partial^{\alpha - \beta } w \quad . \eqno(3.45)$$ 

We now use a minor extension of Lemma 3.2 with $p = 2$, taking advantage of the unitarity of
$U(1/t)$ in $L^2$, with $u = U^*(1/t) \partial^{\alpha} w$, $v = C t^{-2}s$, and $h$ the sum of all
terms not containing $\nabla \partial^{\alpha} w$, thereby obtaining
$$\Big | \partial_t \parallel \partial^{\alpha} w \parallel_2 \Big | \ \leq C t^{-2} \sum_{0 \leq
j \leq k} \parallel \partial^{j+1}s \ \partial ^{k -j}w \parallel_2 \eqno(3.46)$$
\noi from which (3.40L) follows by Lemma 3.4 applied with $s$ replaced by $\partial s$.\\

\noi (3.41). For any multi-index $\alpha$ with $|\alpha | = \ell + 1$, we obtain from (3.5) 
$$\partial_t \partial^{\alpha} s = t^{-2} \sum_{\beta \leq \alpha} (\partial^{\beta} s \cdot \nabla
) \partial^{\alpha - \beta} s + t^{- \gamma} \ \partial^{\alpha} \nabla g (w, w) \quad . 
\eqno(3.47)$$ 
\noi We now use Lemma 3.2 with $p = 2$, $u = \partial^{\alpha} s$, $v = t^{-2} s$, and
$h$ the sum of all terms with $|\beta | \geq 1$ and of the contribution of $g$. The terms quadratic
in $s$ all have at least one derivative on each $s$ and are estimated in $L^2$ by Lemma 3.3 with
$s_1$ and $s_2$ replaced by $\partial s$, while the term containing $g$ is estimated by (3.32) of
Corollary 3.1, thereby yielding (3.41). \\

\noi (3.42). For any multi-index $\alpha$ with $|\alpha | = \ell_0$, we obtain again (3.47) from
(3.5). We estimate the terms quadratic in $s$ directly by Lemma 3.3 with $s_1$ and $s_2$ replaced
by $s$ and $\nabla s$ respectively, and with $\ell$ replaced by $\ell_0$. We estimate the
contribution of $g$ by (3.30) of Corollary 3.1. This yields (3.42). \\

\noi (3.43). The estimate of the terms quadratic in $s$ is obvious and that of $g$ follows from
(3.34) of Corollary 3.1. \\

\noi (3.40R) and (3.44). The norms of $s$ appearing in the LHS of (3.41)-(3.43) are precisely the
norms defining $X^{\ell}$. The norms of $s$ appearing in the middle member of (3.40) and in the RHS
of (3.41)-(3.43) are again norms in the definition of $X^{\ell}$, with the exception of $\parallel
\partial s \parallel_{\infty}$. However by Lemma 3.1
$$\parallel \partial s \parallel_{\infty} \ \leq C \parallel s \parallel_{\infty}^{1 - \sigma} \
\parallel \partial^{\ell + 1} s \parallel_2^{\sigma} \ \leq C |s|K_{\ell} \eqno(3.48)$$
\noi with $1/(\ell + 1) < \sigma = 1/(\ell + 1 - n/2) < 1$ since $\ell > n/2$. This yields (3.40R)
and (3.44) from (3.40L) and (3.41)-(3.43) respectively. \par \nobreak
\hfill $\sq$ \\

The next result is a linear estimate of higher norms needed for regularity. \\

\noi {\bf Lemma 3.8.} {\it Let $0 < \mu \leq n - 2$ and let $(k, \ell )$ be an admissible pair.
Let $I \subset {I\hskip-1truemm R}^+$ be an open interval and let $(w, s)$ be a solution of (3.4)
(3.5) in $L_{loc}^{\infty}(I, H^{k+1} \oplus X^{\ell + 1})$. Then $(w, s)$ satisfies the
following estimates~:}
$$\Big | \partial_t \parallel w; \dot{H}^{k+1} \parallel \Big | \ \leq C \ t^{-2} \left \{
\parallel \partial s \parallel_{\infty} \ \parallel w; \dot{H}^{k+1} \parallel \ + \ \parallel
\partial^2s; L^{\infty} \cap \dot{H}^{\ell} \parallel |w|_k \right \} \quad , \eqno(3.49)$$
$$\Big | \partial_t |w|_{k+1} \Big | \leq C \ t^{-2} \left ( |s|_{\ell} \ |w|_{k+1} + |s|_{\ell +
1} \ |w|_k \right ) \quad , \eqno(3.50)$$
$$\Big | \partial_t \parallel s ; \dot{H}^{\ell + 2} \parallel \Big |  \ \leq C \ t^{-2}
\parallel \partial s; L^{\infty} \cap \dot{H}^{\ell} \parallel \ \parallel
\partial^2s;L^{\infty} \cap \dot{H}^{\ell} \parallel \ + \ C \ t^{-\gamma} |w|_k \ |w|_{k+1}
\quad , \eqno(3.51)$$  $$\Big | \partial_t |s|_{\ell + 1} \Big | \leq C \ t^{-2} |s|_{\ell} \
|s|_{\ell + 1} + C \ t^{- \gamma} |w|_k \ |w|_{k+1} \quad . \eqno(3.52)$$

\noi {\bf Proof.} The proof is very similar to that of Lemma 3.7 and will be sketched briefly.
\par
 
\noi (3.49). We take the $L^2$ norm of (3.45) now with $|\alpha | = k + 1$, we apply (the
same minor extension of) Lemma 3.2 with $p = 2$ and end up with (3.46) with $k$ replaced by $k +
1$. We separate the term $j = 0$ and estimate the remaining terms by Lemma 3.4 with $s$ replaced
by $\partial^2s$. This yields (3.49). \\

\noi (3.51). We take the $L^2$ norm of (3.47), now with $|\alpha | = \ell + 2$, we apply again
Lemma 3.2 followed by Lemma 3.3 with $s_1$ and $s_2$ replaced by $\partial s$ and $\partial^2s$
respectively, in order to estimate the terms quadratic in $s$, and we estimate the term containing
$g$ by (3.33) of Corollary 3.1. This yields (3.51). \\

\noi (3.50) and (3.52) follow immediately from (3.40) (3.49) and from (3.44) (3.51)
respectively, from (3.48) and from the fact that by Lemma 3.1 
$$\parallel \partial^2 s \parallel_{\infty} \ \leq C \parallel \partial s \parallel_{\infty}^{1 -
\sigma} \ \parallel \partial^{\ell + 2} s \parallel_2^{\sigma} \ \leq C |s|_{\ell + 1}
\eqno(3.53)$$
\noi with the same $\sigma$ as in (3.48). \par \nobreak
\hfill $\sq$ \\

We shall also need some estimates for the difference of two solutions of (3.4) (3.5). \\

\noi {\bf Lemma 3.9.} {\it Let $0 < \mu \leq n - 2$ and let $(k, \ell)$ be an admissible pair.
Let $I \subset {I\hskip-1truemm R}^+$ be an open interval and let $(w_1, s_1)$ and $(w_2 , s_2)$
be two solutions of (3.4) (3.5) in $L_{loc}^{\infty}(I, H^k \oplus X^{\ell})$. Let $w_{\pm} = w_1
\pm w_2$ and $s_{\pm} = s_1 \pm s_2$. Then the following estimates hold~:}
$$\begin{array}{ll} \Big | \partial_t |w_-|_{k-1} \Big | &\leq C \ t^{-2} \left \{ \parallel
\partial s_+ ; L^{\infty} \cap \dot{H}^{\ell} \parallel \ |w_-|_{k-1} \ + \ \parallel
s_-;L^{\infty} \cap \dot{H}^{\ell} \parallel |w_+|_k \right \} \\ & \\ &\leq C \ t^{-2} \left \{
|s_+|_{\ell} \ |w_-|_{k-1} + |s_-|_{\ell - 1} \ |w_+|_k \right \} \quad , \end{array}
\eqno(3.54)$$
$$\Big | \partial_t \parallel s_-; \dot{H}^{\ell} \parallel \Big | \ \leq C \ t^{-2} \parallel
s_- ; L^{\infty} \cap \dot{H}^{\ell} \parallel \ \parallel \partial s_+; L^{\infty} \cap
\dot{H}^{\ell} \parallel \ + \ C\ t^{- \gamma} |w_-|_{k-1} \ |w_+|_k \quad ,\eqno(3.55)$$ 
$$\Big | \partial_t \parallel s_- \parallel_{\infty} \Big | \ \leq C \ t^{-2} \parallel s_-
\parallel_{\infty} \ \parallel \partial s_+ \parallel_{\infty} + C\ t^{- \gamma}
|w_-|_{k-1} \ |w_+|_k \quad , \eqno(3.56)$$
$$\Big | \partial_t |s_-|_{\ell - 1} \Big | \leq C \ t^{-2} |s_-|_{\ell - 1} \ |s_+|_{\ell} + C
t^{- \gamma} |w_-|_{k-1} \ |w_+|_k \quad . \eqno(3.57)$$

\noi {\bf Proof.} The proof is again very similar to that of Lemma 3.7 and will be sketched
briefly. \par

\noi (3.54). Taking the difference of (3.4) for $w_1$ and $w_2$ and applying $\partial^{\alpha}$
with $|\alpha | \leq k - 1$ yields
$$\partial_t \ \partial^{\alpha} w_- = (4t^2)^{-1} U(1/t) \sum_{\beta \leq \alpha} \Big \{ \left
( 2 \partial^{\beta} s_+ \cdot \nabla + (\partial^{\beta} \nabla \cdot s_+) \right ) U^*(1/t)
\partial^{\alpha - \beta} w_-$$
$$+ \left ( 2 \partial^{\beta} s_- \cdot \nabla + (\partial^{\beta} \nabla \cdot s_-)\right )
U^*(1/t) \partial^{\alpha - \beta} w_+ \Big \} \quad . \eqno(3.58)$$
\noi We apply (the same minor extension of) Lemma 3.2 with $p = 2$, $u = \partial^{\alpha} w_-$ and
$v = (2t^2)^{-1} s_+$, followed by an application of Lemma 3.4 with $(w, s, k, \ell )$ replaced
by $(w_-, \partial s_+, k - 1, \ell )$ for the terms with $(w_-, s_+)$ and by $(w_+, s_-, k,
\ell )$ for the terms with $(w_+, s_-)$. This yields (3.54L), while (3.54R) follows from the
definition of $X^{\ell}$ and from (3.48). \\

\noi (3.55). Taking the difference of (3.5) for $s_1$ and $s_2$ and applying $\partial^{\alpha}$
with $|\alpha | = \ell$ yields
$$\partial_t \ \partial^{\alpha} s_- = (2t^2)^{-1} \sum_{\beta \leq \alpha} \left \{ 
(\partial^{\beta} s_+ \cdot \nabla ) \partial^{\alpha - \beta} s_- + (\partial^{\beta} s_- \cdot
\nabla ) \partial^{\alpha - \beta } s_+ \right \} + t^{-\gamma} \ \partial^{\alpha} \nabla g (w_-,
w_+) \quad . \eqno(3.59)$$
\noi We apply Lemma 3.2 with $p = 2$, $u = \partial^{\alpha} s_-$, $v = (2t^2)^{-1} s_+$, we 
estimate the terms quadratic in $s$ by Lemma 3.3 with $s_1$ and $s_2$ replaced by $s_-$ and
$\partial s_+$, and the contribution of $g$ by (3.31) of Corollary 3.1. This yields (3.55). \\

\noi (3.56). We proceed in the same way, applying Lemma 3.2 with $p = \infty$ to (3.59) with
$\alpha = 0$, $u = s_-$ and $v = (2t^2)^{-1} s_+$, we estimate the terms quadratic in $s$ in the
obvious way and the contribution of $g$ by (3.34) of Corollary 3.1 This yields (3.56). \\

Finally (3.57) follows from the definition of $X^{\ell}$, from (3.55) and its analogue with $\ell$
replaced by $\ell_0$, from (3.56) and from (3.48). \par \nobreak
\hfill $\sq$ \\ 

\noi {\bf Lemma 3.10.} {\it Let $0 < \mu \leq n - 2$ and let $(k, \ell )$ be an admissible pair.
Let $I \subset {I\hskip-1truemm R}^+$ be an open interval and let $(w_1, s_1)$ and $(w_2, s_2)$ be
two solutions of (3.4) (3.5) in $L_{loc}^{\infty}(I, H^{k+1} \oplus X^{\ell + 1})$ and
$L_{loc}^{\infty}(I, H^k \oplus X^{\ell})$ respectively. Let $w_{\pm} = w_1 \pm w_2$ and $s_{\pm}
= s_1 \pm s_2$. Then the following estimates hold~:}
$$\Big | \partial_t |w_-|_k \Big | \leq C t^{-2} \Big \{ \parallel \partial s_2 ; L^{\infty} \cap
\dot{H}^{\ell} \parallel \ |w_-|_k \ + \ \parallel \partial s_-; L^{\infty} \cap \dot{H}^{\ell}
\parallel |w_1|_k + \ \parallel s_- \parallel_{\infty} \ \parallel w_1 ; \dot{H}^{k+1} \parallel
\Big \}$$ $$\leq C t^{-2} \Big \{ |s_2|_{\ell} \ |w_-|_k + |s_-|_{\ell} \ |w_1|_k + \parallel
s_-\parallel_{\infty} \ |w_1|_{k+1} \big \} \quad , \eqno(3.60)$$ 
$$\Big | \partial_t \parallel s_-; \dot{H}^{\ell + 1} \parallel \Big | \ \leq C t^{-2} \Big \{
\parallel \partial s_- ; L^{\infty} \cap \dot{H}^{\ell} \parallel \sum_{i=1,2} \parallel
\partial s_i ; L^{\infty} \cap \dot{H}^{\ell} \parallel$$
$$+ \ \parallel s_- \parallel_{\infty} \ \parallel s_1 ; \dot{H}^{\ell + 2} \parallel \Big \} + C
t^{- \gamma} |w_-|_k \ |w_+|_k \quad , \eqno(3.61)$$
$$\Big | \partial_t |s_-|_{\ell} \Big | \leq C t^{-2} \Big \{ |s_-|_{\ell} \ (|s_1|_{\ell} +
|s_2|_{\ell}) + \ \parallel s_- \parallel_{\infty} \ |s_1|_{\ell + 1}\Big \} + C t^{-\gamma}
|w_-|_k \ |w_+|_k \quad . \eqno(3.62)$$

\noi {\bf Proof.} The proof is again very similar to that of Lemma 3.7 and will be sketched
briefly. \par

\noi (3.60). We rewrite the difference of (3.4) for $w_1$ and $w_2$ as follows
$$\partial_t w_- = (2t^2)^{-1} U(1/t) \Big \{ (2 s_2 \cdot \nabla + (\nabla \cdot s_2)) U^*(1/t)
w_- + (2s_- \cdot \nabla + (\nabla \cdot s_-)) U^* (1/t) w_1 \big \} \quad . \eqno(3.63)$$
We apply $\partial^{\alpha}$ to (3.63) with $|\alpha | \leq k$, use (the same minor extension
of) Lemma 3.2 with $p = 2$, $u = \partial^{\alpha} w_-$ and $v = (2t^2)^{-1} s_2$, estimate all
the resulting terms by Lemma 3.4 with $(w, s)$ replaced by $(w_-, \partial s_2)$ or by $(w_1 ,
\partial s_-)$, except for the terms $s_- \cdot \nabla \partial^{\alpha} U^*(1/t) w_1$ which we
estimate directly in an obvious way. This yields (3.60L), while (3.60R) follows from the definition
of $X^{\ell}$ and from (3.48). \\

\noi (3.61). We rewrite the difference of (3.5) for $s_1$ and $s_2$ similarly as
$$\partial_t s_- = t^{-2} \left ( (s_2 \cdot \nabla ) s_- + (s_- \cdot \nabla ) s_1 \right ) +
t^{- \gamma} g (w_-, w_+) \quad . \eqno(3.64)$$
\noi We apply $\partial^{\alpha}$ to (3.64) with $| \alpha | = \ell + 1$, we use Lemma 3.2 with $p =
2$, $u = \partial^{\alpha} s_-$ and $v = t^{-2} s_2$, we estimate all the resulting terms quadratic
in $s$ by Lemma 3.3 with $(s_1, s_2)$ replaced by $(\partial s_2, \partial s_-)$ or by $(\partial
s_-, \partial s_1)$, except for the terms $(s_- \cdot \nabla ) \partial^{\alpha} s_1$ which we
estimate directly in an obvious way. The contribution of $g$ is estimated by (3.32) of Corollary
3.1. This yields (3.61). \\

Finally (3.62) follows from the definition of $X^{\ell}$, from (3.57), (3.61) and from (3.48). \par
\nobreak
\hfill $\sq$ \\

\noi {\bf Lemma 3.11.} {\it Let $0 < \mu \leq n - 2$ and let $(k, \ell )$ be an admissible pair.
Let $I \subset {I\hskip-1truemm R}^+$ be an open interval and let $(w_1, s_1)$ and $(w_2, s_2)$ be
two solutions of (3.4) (3.5) in $L_{loc}^{\infty}(I, H^{k} \oplus X^{\ell})$. Let $w_{\pm} = w_1 \pm
w_2$ and $s_{\pm} = s_1 \pm s_2$. Then the following estimates hold~:}
$$\Big | \partial_t |w_-|_{k-2} \Big | \leq C t^{-2} \Big \{ |s_+|_{\ell} \ |w_-|_{k-2} +
|s_-|_{\ell - 2} \  |w_+|_k \Big \} \quad , \eqno(3.65)$$
$$\Big | \partial_t |s_-|_{\ell -2} \Big | \leq C t^{-2} |s_-|_{\ell - 2} \ |s_+|_{\ell} + C t^{-
\gamma} |w_-|_{k - 2} \  |w_+|_k  \quad . \eqno(3.66)$$

\noi {\bf Proof.} The proof is similar to that of Lemma 3.9. \par

\noi (3.65). We apply (the same minor extension of) Lemma 3.2 with $p = 2$, $u = \partial^{\alpha}
w_-$ and $v = (2t^2)^{-1} s_+$ to (3.58) with now $|\alpha | \leq k - 2$. Applying Lemma 3.4 to the
terms containing $(w_-, s_+)$ and omitting the irrelevant operator $U^*(1/t)$ for brevity in the
terms containing $(w_+ , s_-)$, we obtain
$$\Big | \partial_t |w|_{k-2} \Big | \leq C t^{-2} \Big \{ \parallel \partial s_+; L^{\infty}
\cap \dot{H}^{\ell} \parallel \ |w_-|_{k-2} \ + \ \parallel s_- \ w_+ ; H^{k-1} \parallel \Big
\}  \quad . \eqno(3.67)$$ \noi We now distinguish two cases. \par

If $\ell > n/2 + 1$, we estimate the last norm in (3.67) by Lemma 3.4 with $(w, s, k, \ell )$
replaced by $(w_+, s_-, k-1, \ell - 1)$ so that
$$\begin{array}{ll} \parallel s_- \ w_+ ; H^{k-1} \parallel  &\leq C \parallel s_- ;
L^{\infty} \cap \dot{H}^{\ell - 1} \parallel \ |w_+|_{k-1} \\ & \\ &\leq C |s_-|_{\ell -
2} \ |w_+|_k \quad .\end{array} \eqno(3.68)$$

If $n/2 < \ell \leq n/2 + 1$, namely for the lowest admissible value $\ell = \ell_0 + 1$, we
estimate directly
$$\parallel s_- \ w_+ ; H^{k-1} \parallel \ \leq C \sum_{j+j' \leq k-1} \parallel \partial^j s_-
\ \partial^{j'} w_+ \parallel_2$$
$$\leq C \sum_{j+j' \leq k-1} \parallel \partial^j s_- \parallel_{r_1} \ \parallel \partial^{j'}
w_+ \parallel_{r_2} \eqno(3.69)$$
\noi with
$$0 \leq j \leq k - 1 \leq \ell - 1 \leq n/2 \quad ,$$

$$\left \{ \begin{array}{ll} 0 \leq \delta (r_1) = \ell - 1 - j \leq n/2 -
j \quad , \\ \\ 0 \leq \delta (r_2) = n/2 - \ell + 1 + j \quad . \end{array} \right . \eqno(3.70)$$
\noi By Lemma 3.1, we then estimate
$$\parallel \partial^j s_-\parallel_{r_1} \ \leq C \parallel \partial^{\ell - 1} s_- \parallel_2
\ = C \parallel s_-;\dot{H}^{\ell_0} \parallel \eqno(3.71)$$
\noi with the only exception of the case of even $n$, $\ell = n/2 + 1$ and $j = 0$ where $r_1 =
\infty$ and that norm reduces to $\parallel s_- \parallel_{\infty}$, so that in all cases
$$\parallel \partial^j s_- \parallel_{r_1} \ \leq C |s_-|_{\ell - 2} \quad . \eqno(3.72)$$
\noi On the other hand
$$\delta (r_2) + j' \leq n/2 - \ell + 1 + k - 1 < k$$
\noi and therefore
$$\parallel \partial^{j'} w_+ \parallel_{r_2} \ \leq C |w_+|_k \quad . \eqno(3.73)$$
\noi Substituting (3.72) and (3.73) into (3.69), substituting either the result thereof or (3.68)
into (3.67) and using (3.48) yields (3.65). \\

\noi (3.66). We apply Lemma 3.2 with $p = 2$, $u = \partial^{\alpha} s_-$, $v = (2t^2)^{-1} s_+$
to (3.59) with $|\alpha | = \ell - 1$ and obtain
$$\Big | \partial_t \parallel \partial^{\ell - 1} s_- \parallel_2 \Big | \ \leq C t^{-2}
\parallel \partial^{\ell - 1} ((\partial s_+) s_-)\parallel_2 + C t^{-\gamma} \parallel
\partial^{\ell} g (w_- , w_+)\parallel_2 \quad . \eqno(3.74)$$
\noi We distinguish again two cases. \par

If $\ell > n/2 + 1$, we apply Lemma 3.3 with $(s_1, s_2, \ell )$ replaced by $(\partial s_+, s_-,
\ell - 1)$ to estimate
$$\begin{array}{ll} \parallel \partial^{\ell - 1} ((\partial s_+) s_-) \parallel_2 &\leq
C \parallel \partial s_+; L^{\infty} \cap \dot{H}^{\ell - 1} \parallel \ \parallel s_-;
L^{\infty} \cap \dot{H}^{\ell - 1} \parallel \\ & \\ &\leq C |s_+|_{\ell - 1} \
|s_-|_{\ell - 2} \quad . \end{array} \eqno(3.75)$$

If $n/2 < \ell \leq n/2 + 1$, we estimate directly
$$\parallel \partial^{\ell - 1} ((\partial s_+)s_-) \parallel_2 \ \leq \ C \sum_{j \leq \ell - 1}
\parallel \partial^j s_- \parallel_{r_1} \ \parallel \partial^{\ell - j} s_+ \parallel_{r_2}
\eqno(3.76)$$
\noi with $r_1$ and $r_2$ again given by (3.70), so that (3.72) holds as before, while $\delta
(r_2) + \ell - j = n/2 + 1$ so that
$$\parallel \partial^{\ell - j} s_+ \parallel_{r_2} \ \leq C |s_+|_{\ell } \quad . \eqno(3.77)$$
\noi Substituting (3.72) (3.77) into (3.76), and either the result thereof or (3.75) into (3.74)
and estimating the contribution of $g$ by (3.30) of Corollary 3.1 yields
$$\Big | \partial_t \parallel \partial^{\ell - 1} s_- \parallel_2 \Big | \ \leq C t^{-2}
|s_-|_{\ell - 2} \ |s_+|_{\ell} + C t^{- \gamma} |w_-|_{k-2} \  |w_+|_k \quad . \eqno(3.78)$$

In the case of even $n$, we estimate in addition 
$$\Big | \partial_t \parallel s_- \parallel_{\infty} \Big | \ \leq C t^{-2} \parallel s_-
\parallel_{\infty} \ \parallel \partial s_+ \parallel_{\infty} + Ct^{- \gamma} |w_-|_{k-2} \
|w_+|_k \eqno(3.79)$$ \noi in the same way as in the proof of (3.56), but estimating now the
contribution of $g$ by (3.37) instead of (3.34) of Corollary 3.1. \par

Collecting (3.78), its analogue with $\ell$ replaced by $\ell_0 + 1$, and in addition (3.79) for
even $n$ yields (3.66). \par \nobreak
\hfill $\sq$ \\

\noi {\bf Remark 3.2.} In Lemmas 3.7-3.11, not all the properties in the definition of
admissibility for $(k , \ell )$ are used in every single estimate. The condition $\ell > n/2$ is
used in many places. However the condition $k \leq \ell$ is used only in the estimates of $w$ or
$w_-$ from (3.4), but not in the estimates of $s$ or $s_-$ from (3.5). Conversely the condition
(3.29) is used only in the estimates of $s$ or $s_-$ from (3.5), but not in the estimates of $w$
or $w_-$ from (3.4). Furthermore that condition is used only through the estimates
(3.30)-(3.35) of Corollary 3.1, so that (3.29) could be replaced by (3.30)-(3.35) in the
definition of admissibility, thereby opening the possibility of treating more general
nonlinearities $\widetilde{g}$ than simply (1.2). Finally the condition (3.36) for $n$ even is
used only through (3.37) in Lemma 3.11.

\section{The Cauchy problem at finite times.}
\hspace*{\parindent} In this section, we study the Cauchy problem for the basic system (3.4) (3.5)
at finite times. We use the basic spaces $H^k$ and $X^{\ell}$ defined at the beginning of Section
3, as well as the notation (3.11) for the norm in those spaces. Admissible pairs $(k, \ell)$ are
defined in Definition 3.1. We prove that the Cauchy problem for the system (3.4) (3.5) is locally
well posed in ${I\hskip-1truemm R}^+\setminus \{0\}$ for positive initial time $t_0$ and initial
data in $H^k \oplus X^{\ell}$ for admissible $(k, \ell)$. A similar result holds in ${I\hskip-1truemm
R}^-\setminus \{ 0\}$. We make no effort to study the situation as $t \to 0$, since that is of no
interest for later purposes, and since the system (3.4) (3.5) is singular at $t = 0$ anyway
because of the choice (1.2) of $g_0$ and of the change of variables (2.14) from $u$ to $(w, s)$.
The main result can be stated as follows. \\

\noi {\bf Proposition 4.1.} {\it Let $\gamma > 0$, $n \geq 3$ and $0 < \mu \leq n - 2$. Let $(k, \ell
)$ be an admissible pair. Let $t_0 > 0$. Then for any $(w_0, s_0) \in H^k \oplus X^{\ell}$, there
exist $T_{\pm}$ with $0 \leq T_- < t_0 < T_+ \leq \infty$ such that~:} \par

{\it (1) The system (3.4) (3.5) has a unique solution $(w , s) \in {\cal C}(I, H^k \oplus X^{\ell})$
with $(w, s) (t_0) = (w_0, s_0)$, where $I = (T_- , T_+)$. If $T_- > 0$ (resp. $T_+ < \infty$),
then $|w(t)|_k + |s(t)|_{\ell} \to \infty$ when $t$ decreases to $T_-$ (resp. increases to
$T_+$).} \par

{\it (2) If $(w_0, s_0) \in H^{k'} \oplus X^{\ell '}$ for some admissible pair $(k', \ell ')$ with
$k' \geq k$ and $\ell ' \geq \ell$, then $(w, s) \in {\cal C}(I, H^{k'} \oplus X^{\ell '})$.} \par

{\it (3) For any compact subinterval $J \subset \subset I$, the map $(w_0, s_0) \to (w, s)$ is
continuous from $H^{k-1} \oplus X^{\ell - 1}$ to $L^{\infty}(J, H^{k-1} \oplus X^{\ell - 1})$
uniformly on the bounded sets of $H^k \oplus X^{\ell}$, and is pointwise continuous from $H^k
\oplus X^{\ell}$ to $L^{\infty}(J,H^k \oplus X^{\ell})$.} \\

\noi {\bf Proof.} Most of the proof proceeds by standard arguments, and we shall mainly
concentrate on those which are not. We concentrate on the case of increasing time, namely $t \geq
t_0$. The case of decreasing time $t \leq t_0$ can be treated in the same way, possibly after
changing $t$ to $1/t$ and $s$ to $- s$, thereby transforming (3.4) (3.5) into the system
$$\partial_t w = (1/2) U(t)(2s \cdot \nabla + (\nabla \cdot s)) U(-t) w \eqno(4.1)$$
$$\partial_ts = (s \cdot \nabla ) s + t^{\gamma - 2} \ \nabla g(w, w) \eqno(4.2)$$
\noi and considering that system for increasing time. \par

The (negative) powers of $t$ in the coefficients of (3.4) (3.5) are bounded on compact subintervals
of $[t_0, \infty )$, actually bounded on $[t_0, \infty )$, and play no role in the present problem.
We omit them for brevity. \par

The proof proceeds in several steps. \\

\noi {\bf Step 1.} We introduce a parabolic regularization and consider the system
$$\partial_t w = \eta \Delta w + U(1/t) (s \cdot \nabla + (1/2)(\nabla \cdot s)) U^{*}(1/t) w \equiv
\eta \Delta w + F(w, s) \eqno(4.3)$$
$$\partial_t s = \eta \Delta s + (s \cdot \nabla ) s + \nabla g (w, w) \equiv \eta \Delta s +
G(w, s) \eqno(4.4)$$
\noi with $\eta > 0$. The Cauchy problem for the system (4.3) (4.4) can be recast in
the integral form 
$${w \choose s} (t) = V_{\eta}(t - t_0) {w_0 \choose s_0} + \int_{t_0}^t dt' \ V_{\eta}(t - t')
{F(w,s) \choose G(w, s)} (t') \eqno(4.5)$$
\noi where $V_{\eta}(t) = \exp (\eta t \Delta )$. The operator $V_{\eta}(t)$ is a contraction in
$H^k \oplus X^{\ell}$, while the operator $\nabla V_{\eta}(t)$ satisfies the bound 
$$\parallel \nabla V_{\eta}(t) ; {\cal L}(H^k \oplus X^{\ell}) \parallel \ \leq C (\eta t)^{-1/2}
\quad . \eqno(4.6)$$
\noi From these facts and from estimates on $F$, $G$ similar to, but simpler than, those in Lemma
3.7, it follows by  a standard contraction argument that there exists $T > 0$ depending only on
$\eta$ and on $|w_0|_k + |s_0|_{\ell}$ such that the system (4.5) has a unique solution
$(w_{\eta},s_{\eta}) \in {\cal C}([t_0, t_0 + T], H^k \oplus X^{\ell})$. \\

\noi {\bf Step 2.} Estimates uniform in $\eta$. We estimate $(w_{\eta}, s_{\eta})$ by
Lemma 3.7, taking into account the fact that by Lemma 3.2, the term $\eta \Delta w$ in (4.3)
(4.4) does not contribute to the estimates. Let
$$y(t) = |w_{\eta}(t)|_k \quad , \quad z(t) = |s_{\eta}(t)|_{\ell} \quad .
\eqno(4.7)$$
\noi We obtain from Lemma 3.7 (with the powers of $t$ omitted)
$$\left \{ \begin{array}{l} \partial_t y \leq Cyz \\ \\ \partial_t z \leq C(y^2 + z^2)
\end{array} \right .$$
\noi and by integration, with $y_0 = y(t_0)$, $z_0 = z(t_0)$,
$$y(t) + z(t) \leq (y_0 + z_0) \left ( 1 - C(t - t_0) (y_0 + z_0) \right )^{-1} \leq 2 (y_0 +
z_0)$$
\noi for $2C(t - t_0)(y_0 + z_0) \leq 1$, so that for some $T$ depending only on $|w_0|_k +
|s_0|_{\ell}$, $(w_{\eta}(t), s_{\eta}(t))$ is estimated a priori in ${\cal C}([t_0, t_0 + T],
H^k \oplus X^{\ell})$ uniformly in $\eta$. By a standard globalisation argument, the solution
constructed in Step 1 can be extended to that new interval, now independent of $\eta$. \\

\noi {\bf Step 3.} Limit $\eta \to 0$. Let $I_0 = [t_0, t_0 + T]$ be the interval,
independent of $\eta$, obtained in Step 2. We now prove that $(w_{\eta}, s_{\eta})$ converges in
norm in $L^{\infty}(I_0, H^{k-1} \oplus X^{\ell - 1})$. We know already that $(w_{\eta},
s_{\eta})$ is estimated uniformly in $\eta$ according to
$$\parallel w_{\eta} ; L^{\infty}(I_0, H^k)\parallel \ \leq a \quad , \eqno(4.8)$$
$$\parallel s_{\eta} : L^{\infty} (I_0, X^{\ell}) \parallel \ \leq b \quad .
\eqno(4.9)$$
\noi Let now $\eta_1, \eta_2 > 0$ and let $w_i = w_{\eta_i}$, $i = 1,2$. We estimate the
difference $(w_1 - w_2, s_1 - s_2)$ by Lemma 3.9, except for the contribution of the terms coming
from $(\eta \Delta w, \eta \Delta s)$ in (4.3) (4.4), which are estimated directly
as follows. For $|\alpha | \leq k - 1$, we estimate
\begin{eqnarray*}
\partial_t \parallel \partial^{\alpha}(w_1 - w_2)\parallel_2^2 &\leq & 2 {\rm Re} < \partial^{\alpha}
(w_1 - w_2), \eta_1 \partial^{\alpha} \Delta w_1 - \eta_2 \partial^{\alpha} \Delta w_2 > \\
&+&  \hbox{other terms} \\
&\leq & 4a^2 (\eta_1 + \eta_2) + \hbox{other terms}
\end{eqnarray*}
\noi where the other terms are those coming from Lemma 3.9. The contribution of the terms $\eta
\Delta s$ from (4.4) to the estimate of $s_1 - s_2$ are treated in the same way. Defining now
$$y(t) = |(w_1 - w_2)(t)|_{k-1} \quad , \quad z(t) = |(s_1 - s_2)(t)|_{\ell - 1} \quad ,
\eqno(4.10)$$
\noi and combining the previous estimates with Lemma 3.9, we obtain
$$\left \{ \begin{array}{l} \partial_t y^2 \leq C (a^2(\eta_1 + \eta_2 ) + by^2 + ayz) \\ \\
\partial_t z^2 \leq C (b^2 (\eta_1 + \eta_2) + bz^2 + ayz) \end{array} \right . \eqno(4.11)$$
\noi for $t \in I_0$, which together with $y(t_0) = 0$, $z(t_0) = 0$, implies that $y$ and $z$ tend
to zero uniformly in $I_0$ when $\eta_1, \eta_2 \to 0$ by Gronwall's Lemma. As a consequence
$(w_{\eta}, s_{\eta})$ converges to a limit $(w, s)$ in norm in $({\cal C} \cap L^{\infty}) (I_0,
H^{k-1} \oplus X^{\ell - 1})$. Clearly, $(w, s)$ is a solution of the system (3.4)(3.5) with the
appropriate initial data. Furthermore, since $(w_{\eta}, s_{\eta})$ is uniformly bounded in
$L^{\infty}(I_0, H^k \oplus X^{\ell})$, it follows from a standard compactness argument that $(w,
s)$ belongs to that space with the same bound and that $(w_{\eta}, s_{\eta})$ converges to $(w,
s)$ in that space in the weak-$*$ sense.  \\

\noi {\bf Step 4.} Uniqueness. That step is independent of the previous ones and could
equally well have been made at the very beginning. Actually uniqueness in $L^{\infty}(\cdot , H^k
\oplus X^{\ell})$ follows immediately from Lemma 3.9 and from Gronwall's Lemma.\\

\noi {\bf Step 5.} Regularity. It follows immediately from Lemmas 3.7 and 3.8 and from
Gronwall's Lemma that a solution $(w, s) \in L^{\infty}(I, H^k \oplus X^{\ell})$ with initial data
$(w, s) (t_0) \in H^{k+1} \oplus X^{\ell + 1}$ belongs to $L^{\infty}(I, H^{k+1} \oplus X^{\ell +
1})$ for the same interval $I$. A similar regularity with general $(k', \ell ')$ follows by
iteration. \\

Using the previous five steps and standard arguments, one can then prove most of Proposition 4.1,
with the only restriction that Parts (1) (3) hold only with ${\cal C}(\cdot , H^k \oplus X^{\ell})$
replaced by ${\cal C}(\cdot, H^{k-1} \oplus X^{\ell - 1}) \cap L_{loc}^{\infty}(\cdot , H^k \oplus
X^{\ell})$, with continuity in Part (3) being in norm in the former space and in the weak-$*$
sense in the latter, while Part (2) holds only with a similar restriction. \par

We now turn to the proof of the missing continuities, which is more delicate. We follow a method
used in \cite{1r}. We first derive an additional estimate for the difference of two solutions of
(3.4) (3.5) in $L^{\infty} (I, H^k \oplus X^{\ell})$ for some interval $I$. For brevity we
introduce the short hand notation $y = (w, s)$, $y_0 = (w_0 , s_0)$ and for two solutions $y_i =
(w_i , s_i)$, $y_{0i} = (w_{0i}, s_{0i})$, $i = 1, 2$, $y_- = y_1 - y_2$ and $y_{0-} = y_{01} -
y_{02}$. Furthermore for $(k, \ell )$ an admissible pair and $\theta$ an integer, $-2 \leq \theta
\leq 1$, we denote $j + \theta = (k + \theta , \ell + \theta )$ and $Y^{j + \theta} = H^{k+ \theta}
\oplus X^{\ell + \theta}$. Let now $y_i \in L^{\infty} (I, Y^j)$, $i = 1,2$, be two solutions of
(3.4) (3.5) with initial data $y_{0i}$ at time $t_0 \in I$ for some compact interval $I$, satisfying
the estimate
$$\parallel y_i ; L^{\infty}(I, Y^j)\parallel \ \leq a < \infty \quad , \quad i = 1,2 \quad .
\eqno(4.12)$$ 
\noi Assume furthermore that $y_{01} \in Y^{j+1}$, so that by Step 5 $y_1 \in L^{\infty} (I,
Y^{j+1})$. We now estimate $y_-$ in $Y^j$ by Lemma 3.10. From (3.60) (3.62) we obtain
$$\partial_t |y_-|_j \leq C \left ( a |y_-|_j + \parallel s_- \parallel_{\infty} \ |y_1|_{j+1} \right
) \quad . \eqno(4.13)$$
\noi By Lemma 3.8, esp. (3.50) (3.52), we estimate
$$\partial_t |y_1|_{j+1} \leq C a |y_1|_{j+1}$$
\noi and therefore by Gronwall's Lemma, for all $t \in I$,
$$|y_1|_{j+1} \leq C(a, |I|)|y_{01}|_{j+1} \quad .
\eqno(4.14)$$
\noi We next estimate for $\ell > n/2$, possibly by using (3.10),
$$\left \{ \begin{array}{ll} \parallel s_-\parallel_{\infty} \ \leq \ C|s_-|_{\ell - 2} &\quad
\hbox{for} \ \ell \geq n/2 + 1 \quad , \\ & \\
\parallel s_- \parallel_{\infty} \ \leq \ C|s_-|^{1/2}_{\ell - 2} \ |s_-|^{1/2}_{\ell - 1}
&\quad \hbox{for $n$ odd,} \ \ell = (n+1)/2 \quad . \end{array} \right . \eqno(4.15)$$
\noi We now estimate
$$\partial_t |y_-|_{j - \theta} \leq C a |y_-|_{j - \theta}$$
\noi with $\theta = 1$, by Lemma 3.9, esp. (3.54) (3.57), and with $\theta = 2$ by Lemma 3.11,
esp. (3.65) (3.66), so that by Gronwall's Lemma again, for $\theta = 1,2$ and for all $t \in I$
$$|y_-|_{j - \theta } \leq C (a, |I|) |y_{0-}|_{j - \theta} \eqno(4.16)$$
\noi and therefore by (4.15), for all $t \in I$,
$$\parallel s_- \parallel_{\infty} \ \leq C(a, |I|)\parallel y_{0-}\parallel_b \eqno(4.17)$$
\noi where
$$\left \{ \begin{array}{ll} \parallel y_{0-}\parallel_b \ = \ |y_{0-}|_{j-2} &\quad \hbox{for}
\ \ell \geq n/2 + 1 \quad , \\ \\
\parallel y_{0-} \parallel_b \ = \ |y_{0-}|_{j-2}^{1/2}\  |y_{0-}|_{j-1}^{1/2} &\quad \hbox{for
$n$ odd,} \ \ell = (n+1)/2 \quad . \end{array} \right . \eqno(4.18)$$
\noi Substituting (4.14) and (4.17) into (4.13) and applying Gronwall's Lemma again, we obtain for
all $t \in I$
$$|y_-|_j \leq C(a, |I|) \left ( |y_{0-}|_j + \parallel y_{0-}\parallel_b \ |y_{01}|_{j+1}
\right ) \quad .\eqno(4.19)$$

We now come back to the proof of the missing continuities, which will make an essential use of the
estimate (4.19). \\

\noi {\bf Step 6.} Continuity of the solutions in $H^k \oplus X^{\ell}$. Let $I \subset
{I\hskip-1truemm R}^+ \setminus \{ 0 \}$ be a compact interval and let $(w, s) \in {\cal C}(I,
H^{k-1} \oplus X^{\ell - 1}) \cap L^{\infty}(I, H^k \oplus X^{\ell})$ be solution of the system
(3.4) (3.5) with initial data $(w_0, s_0) \in H^k \oplus X^{\ell}$ at some time $t_0 \in I$. We
shall prove that $(w, s) \in {\cal C}(I, H^k \oplus X^{\ell})$. We use the short hand notation $y$,
$y_0$, $Y^j$, etc. introduced above. We introduce a regularisation defined as follows. We choose
a function $\psi_1 \in {\cal S}({I\hskip-1truemm R}^n)$ such that $\int dx \psi_1 (x) = 1$ and
such that $|\xi |^{-m} (\widehat{\psi}_1(\xi ) - 1 )|_{\xi = 0} = 0$ for $m = 1,2$, we define
$\psi_{\varepsilon}(x) = \varepsilon^{-n} \psi_1(x/\varepsilon )$ so that
$\widehat{\psi}_{\varepsilon} (\xi ) = \widehat{\psi}_1(\varepsilon \xi )$ and we define the
regularisation by $f \to f_{\varepsilon} = \psi_{\varepsilon} * f$ for all $f \in {\cal S}'$.
Clearly the regularisation is a bounded operator with norm at most $\parallel \psi_1
\parallel_1$ and tends strongly to the unit operator when $\varepsilon \to 0$ in $Y^j$ for all
relevant $j$. \par

We now regularize the initial data $y_0$ to $y_{0 \varepsilon}$. By the previous steps,
$y_{0\varepsilon}$ generates a solution $y_{\varepsilon}$ of (3.4) (3.5). For $\varepsilon$
sufficiently small, and possibly after a small restriction of $I$, that solution can be assumed
to be in $L^{\infty}(I, Y^j)$ for the same interval $I$ as $y$ and to be bounded there uniformly
in $\varepsilon$, namely 
$$\parallel y; L^{\infty} (I, Y^j) \parallel \vee \parallel y_{\varepsilon} ; L^{\infty} (I,
Y^j)\parallel \ \leq a < \infty \quad . \eqno(4.20)$$
\noi Furthermore, since $y_{0 \varepsilon} \in Y^{j+1}$, by regularity (Step 5),
$y_{\varepsilon} \in {\cal C}(I, Y^j) \cap L^{\infty}(I, Y^{j+1})$. In order to prove that $y
\in {\cal C}(I, Y^j)$ it is therefore sufficient to prove that $y_{\varepsilon}$ converges to
$y$ in norm in $L^{\infty}(I, Y^j)$. For that purpose we apply the estimate (4.19) with
$y_1 = y_{\varepsilon}$, $y_2 = y$, $y_- = y_{\varepsilon} - y$. Now for all $f$     
$$\parallel \partial f_{\varepsilon}\parallel_2 \ = \ \parallel \partial (\psi_{\varepsilon} *
f)\parallel_2 \ \leq \ \parallel \partial \psi_{\varepsilon} \parallel_1 \ \parallel f\parallel_2 \
\leq \varepsilon^{-1} \parallel \partial \psi_1 \parallel_1 \ \parallel f\parallel_2
\quad , \eqno(4.21)$$ $$\begin{array}{ll}
\parallel f_{\varepsilon} - f \parallel_2  &=  \parallel \psi_{\varepsilon} * f - f \parallel_2 \
= \ \parallel \left ( \widehat{\psi}_{\varepsilon} (\xi ) - 1 \right ) \widehat{f}(\xi ) \parallel_2
\\  & \\ &\leq  \varepsilon^{\theta} \parallel |\xi |^{- \theta} \left ( \widehat{\psi}_1 (\xi )
- 1 \right ) \parallel_{\infty} \ \parallel f; \dot{H}^{\theta} \parallel
\end{array}\eqno(4.22)$$  \noi for $\theta = 1, 2$. Furthermore for $n$ even and $s \in
L^{\infty} \cap \dot{H}^{n/2} \cap \dot{H}^{n/2 + \theta}$ so that $|\xi |^{n/2} \widehat{s}(\xi
) \in L^2$, and for $\theta = 1,2$, $$\parallel \psi_{\varepsilon} * s - s\parallel_{\infty} \
\leq \ \parallel (\widehat{\psi}_{\varepsilon} - 1) \widehat{s} \parallel_1 \ \leq
\varepsilon^{\theta} \parallel |\xi |^{-n/2 - \theta} \left ( \widehat{\psi}_1 (\xi ) - 1 \right
) \parallel_2 \ \parallel s; \dot{H}^{n/2+\theta} \parallel \quad . \eqno(4.23)$$

It follows from (4.21)-(4.23) that for $\theta = 1,2$,
$$\left \{ \begin{array}{l} |y_{0 \varepsilon}|_{j+1} \ \leq \ C \ \varepsilon^{-1} |y_0|_j \quad ,
\\ \\
|y_{0-}|_{j - \theta} \ \leq \ C \ \varepsilon^{\theta} |y_0|_j \quad . \end{array} \right .
\eqno(4.24)$$
\noi Substituting the estimate (4.24) into (4.18) (4.19) and using the fact that
$|y_{0-}|_j \to 0$ when $\varepsilon \to 0$ shows that $|y_-|_j$ tends to zero when $\varepsilon
\to 0$ uniformly for $t \in I$, which completes the proof. \\

\noi {\bf Step 7.} Continuity with respect to initial data in ${\cal C}(\cdot , H^k
\oplus X^{\ell})$. From Steps 1-5 and general arguments, it follows only that the solution $(w,
s)$ so far constructed is norm continuous in ${\cal C}(\cdot , H^{k-1} \oplus X^{\ell - 1})$ and
weak-$*$ continuous in $L^{\infty}(\cdot , H^k \oplus X^{\ell})$ as a function of the initial
data $(w_0, s_0) \in H^k \oplus X^{\ell}$. We now prove strong continuity in ${\cal C}(\cdot ,
H^k \oplus X^{\ell})$. We use again the short hand notation $y$, $y_0$, $Y^j$, etc. introduced
above. Let $I \subset {I\hskip-1truemm R}^+ \setminus \{ 0 \}$ be a compact interval and let $y
\in ({\cal C} \cap L^{\infty})(I, Y^j)$ be a fixed solution of the system (3.4) (3.5) with
initial data $y_0 \in Y^j$ at some time $t_0 \in I$. Let $y'_0$ be an initial data in a small
neighborhood of $y_0$ in $Y^j$ and let $y'$ be the solution of (3.4) (3.5) thereby generated. We
also consider the regularized initial data $y_{0 \varepsilon}$ and $y'_{0 \varepsilon}$, and the
solutions $y_{\varepsilon}$ and $y'_{\varepsilon}$ thereby generated. By taking $y'_0$
sufficiently close to $y_0$ and $\varepsilon$ sufficiently small and possibly after a small
restriction of $I$, we can assume that $y'$ and $y_{\varepsilon}$, $y'_{\varepsilon}$ are in
$L^{\infty}(I, Y^j)$ for the same interval $I$ as $y$, and are bounded there uniformly in $y'_0$
and in $\varepsilon$ so that both (4.20) and its analogue with $y$ replaced by $y'$ hold. For
all $t \in I$, we estimate $$|y - y'|_j \ \leq \ |y_{\varepsilon} - y|_j + |y_{\varepsilon} -
y'_{\varepsilon}|_j + |y'_{\varepsilon} - y'|_j  \eqno(4.25)$$ 
\noi and we estimate the three norms in the RHS by (4.19) with $(y_1, y_2)$ replaced by
$(y_{\varepsilon}, y),(y_{\varepsilon}, y_{\varepsilon '})$ and $(y_{\varepsilon '}, y')$
respectively. Using in addition the first inequality in (4.24), we obtain
$$|y - y'|_j \ \leq \ C (a, |I|) \Big ( |y_{0 \varepsilon} - y_0|_j + |y_{0 \varepsilon} - y'_{0
\varepsilon }|_j + |y'_{0 \varepsilon} - y'_0|_j$$
$$+ a \varepsilon^{-1} \left ( \parallel y_{0 \varepsilon} - y_0 \parallel_b \ + \ \parallel
y_{0 \varepsilon} - y'_{0 \varepsilon} \parallel_b \ + \ \parallel y'_{0 \varepsilon} - y'_0
\parallel_b \right ) \Big ) \quad .\eqno(4.26)$$
\noi We now estimate
\begin{eqnarray*}
|y_{0 \varepsilon} - y'_{0 \varepsilon}|_j &\leq & C |y_0 - y'_0|_j \\
|y'_{0 \varepsilon} - y'_0|_j &\leq & |y'_{0 \varepsilon} - y_{0 \varepsilon}|_j + |y_{0
\varepsilon} - y_0|_j + |y_0 - y'_0|_j \\
&\leq & |y_{0 \varepsilon} - y_0|_j + C |y'_0 - y_0|_j  
\end{eqnarray*}
\noi and similarly 
\begin{eqnarray*}
\parallel y_{0 \varepsilon} - y'_{0 \varepsilon}\parallel_b &\leq & C \parallel y_0 - y'_0
\parallel_b \ \leq C |y_0 - y'_0|_j \\ 
\parallel y'_{0 \varepsilon} - y'_0\parallel_b &\leq & \parallel y_{0 \varepsilon} -
y_{0}\parallel_b \ + \ C \parallel y'_{0} - y_0\parallel_b \\
&\leq & \parallel y_{0 \varepsilon} - y_0\parallel_b \ + \ C |y'_0 - y_0|_j \quad . \end{eqnarray*}
\noi Substituting those estimates into (4.26) yields
$$|y - y'|_j \leq C(a, |I|) \left ( |y_{0 \varepsilon} - y_0|_j + a \varepsilon^{_-1}
\parallel y_{0 \varepsilon} - y_0 \parallel_b + (1 + a \varepsilon^{-1}) |y_0 - y'_0|_j \right )
\eqno(4.27)$$  
\noi which can be made to tend to zero uniformly for $t \in I$ by letting $y'_0$ tend to $y_0$ in
$Y^j$ and letting $\varepsilon$ tend to zero, in that order. \par \nobreak
\hfill $\sq$ \\

\noi {\bf Remark 4.1.} Whereas the map $y_0 \to y$ is uniformly continuous from $Y^j$ to $({\cal C}
\cap L^{\infty})(\cdot , Y^{j-1})$ on the bounded sets of $Y^j$, it is only pointwise continuous
from $Y^j$ to $({\cal C} \cap L^{\infty})(\cdot , Y^j)$. In fact Step 7 is performed for fixed
$y_0$, and does not yield an estimate of $\parallel y - y';L^{\infty}(\cdot, Y^j)\parallel$ in
terms of $|y_0 - y'_0|_j$. This is a standard situation in that kind of problems.

 \section{The auxiliary system at infinite time. Existence and asymptotics I.}
\hspace*{\parindent} 
In this section we study the existence of solutions in a neighborhood of infinity in time for the
auxiliary system
$$\hskip 1.5 truecm \left \{ \begin{array}{ll} \partial_tw = (2t^2)^{-1} \ U(1/t) (2s \cdot \nabla +
(\nabla \cdot s)) U^*(1/t) w &\hskip 2 truecm (2.30)\equiv (3.4)\equiv(5.1) \\ & \\
\partial_t s = t^{-2}(s \cdot \nabla ) s + t^{-\gamma} \ \nabla g(w, w) &\hskip 2 truecm
(2.31)\sim (3.5) \equiv (5.2) \end{array} \right .$$
\noi where $g$ is defined by (3.1) (3.2), and we study the asymptotic behaviour in time of those
solutions by essentially constructing local wave operators at infinity for the system (5.1) (5.2)
as compared with the auxiliary free equation
$$\partial_t s_0 = t^{- \gamma} \ \nabla g(w_+, w_+) \quad .
\eqno(2.33) \sim (3.7) \equiv(5.3)$$ \noi The general solution of (5.3) can be written as
$$s_0(t) = s_0(1) + \int_1^t dt' \ t'^{-\gamma} \ \nabla g (w_+ , w_+) \quad . \eqno(5.4)$$

\noi Since from (5.2) and (5.4) the functions $s(t)$ and $s_0(t)$ are expected to increase as
$t^{1- \gamma}$, we define the functions
$$\widetilde{s}(t) = t^{\gamma - 1} \ s(t) \qquad , \quad \widetilde{s}_0 (t) = t^{\gamma - 1} \
s_0(t) \eqno(5.5)$$
\noi which are expected to be bounded in time. We shall use those functions throughout this
section. \noi It follows from Corollary 3.1 that for admissible $(k, l)$ and for $(w_+, s_0 (1)) \in H^k
\oplus X^{\ell}$, $\widetilde{s}_0 (t) \in ({\cal C} \cap L^{\infty})([1 , \infty ), X^{\ell})$ and
that $s_0(t)$ satisfies the estimate
$$\parallel \widetilde{s}_0 ; L^{\infty}([1, \infty ); X^{\ell})\parallel \ \leq |s_0(1)|_{\ell} +
C(1 - \gamma)^{-1} |w_+|_k^2 \quad .$$

We shall use the basic spaces $H^k$ and $X^{\ell}$ defined at the beginning of Section 3, as well
as the notation (3.11) for the norms in those spaces. We recall that admissible pairs $(k, \ell)$
are defined in Definition 3.1. In all this section, we assume that $n \geq 3$ and $0 < \mu \leq n
- 2$. The letter $C$ in subsequent estimates will denote various constants depending on $n$,
$\mu$ and possibly on an admissible pair $(k, \ell)$. On the other hand we shall keep the
dependence of the estimates on $\gamma$ sufficiently explicit for the constants $C$ to be
uniform in $\gamma$ in the range of $\gamma$ where the estimates are stated. For instance a
factor $\gamma^{-1}$ will be kept explicitly in estimates valid for all $\gamma > 0$, but will
be included in $C$ for estimates valid for $\gamma > 1/2$. \par

In the first three propositions, we study the existence of solutions of the system (5.1) (5.2)
defined in a neighborhood of infinity and some of their asymptotic pro\-per\-ties. All those results
hold for all $\gamma > 0$, but we restrict our attention to $0 < \gamma < 1$ in order to simplify
the exposition, as explained in the introduction. Most of those results are consequences of
Proposition 4.1 and of a priori estimates where we now carefully keep track of the time dependence.
We begin with the existence of solutions defined in a neighborhood of infinity in time. \\

\noi {\bf Proposition 5.1.} {\it Let $0 < \gamma < 1$ and let $(k, \ell )$ be an admissible pair. Let
$(w_0, \widetilde{s}_0) \in H^k \oplus X^{\ell}$ and define $y_0 = |w_0|_k$ and $\widetilde{z}_0 =
|\widetilde{s}_0|_{\ell}$. Then there exists $T_0 < \infty$, depending on $y_0$, $\widetilde{z}_0$,
such that for all $t_0 \geq T_0$, there exists $T \leq t_0$, depending on $y_0$, $\widetilde{z}_0$
and $t_0$, such that the system (5.1) (5.2) with initial data $w(t_0) = w_0$, $s(t_0) =
\widetilde{s}_0 t_0^{1 - \gamma}$, has a unique solution $(w, s)$ such that $(w, \widetilde{s})
\in ({\cal C} \cap L^{\infty})([T, \infty ), H^k \oplus X^{\ell})$. One can take} 
$$\gamma T_0^{\gamma} = C \left ( \widetilde{z}_0 + (1 - \gamma )^{-1} y_0^2 \right )
\quad , \eqno(5.6)$$ $$T =  \gamma \ T_0^{\gamma} \ t_0^{1 - \gamma} \quad , \eqno(5.7)$$
\noi {\it and the solution $(w, s)$ is estimated by}
$$|w|_k \leq 2 y_0 \eqno(5.8)$$
$$|s|_{\ell} \leq \left ( 2 \widetilde{z}_0 + C(1 - \gamma )^{-1} y_0^2 \right ) \left ( t_0 \vee
t \right )^{1 - \gamma} \eqno(5.9)$$
\noi {\it for all $t \geq T$.} \\

\noi {\bf Proof.} The result follows from Proposition 4.1 and standard globalisation arguments,
provided we can derive (5.8) (5.9) as a priori estimates under the assumptions made on $t_0$ and
$t$. \par

Let $(w, s)$ be the maximal solution of (5.1) (5.2) with the appropriate initial condition at
$t_0$ obtained by Proposition 4.1 and define $y = |w|_k$ and $z = |s|_{\ell}$. By Lemma 3.7, $y$
and $z$ satisfy
$$\left \{  \begin{array}{l} |\partial_t y | \leq C \ t^{-2} \ yz \\ \\ |\partial_t z| \leq C \
t^{-2}\  z^2 + C t^{- \gamma} \ y^2 \end{array}\right .\eqno(5.10)$$
\noi and we estimate $y$ and $z$ from those inequalities, taking $C = 1$ for the rest of the
proof. We distinguish two cases. \\

\noi {\bf Case ${\bf t} \geq {\bf t_0}$.} Let $\bar{t} > t_0$ and define $Y =
Y(\bar{t}) = \parallel y; L^{\infty}([t_0, \bar{t}\ ])\parallel$ and $Z = Z(\bar{t}) =$\break
\noindent $\parallel t^{\gamma - 1} z;L^{\infty}([t_0, \bar{t}\ ])\parallel$. Then for all $t \in
[t_0, \bar{t}\ ]$
$$\left \{ \begin{array}{l} \partial_t y \leq t^{-1 - \gamma} \ YZ \\ \\ \partial_t z \leq t^{-2
\gamma} \ Z^2 + t^{-\gamma} \ Y^2 \end{array} \right . \eqno(5.11)$$
\noi and therefore by integration with the appropriate initial condition at $t_0$
$$\left \{ \begin{array}{l}  Y \leq y_0 + \gamma^{-1} \ t_0^{- \gamma} \ YZ \\ \\ Z \leq
\widetilde{z}_0 + t_0^{- \gamma} \ Z^2 + (1 - \gamma )^{-1} \ Y^2\end{array}\right .\eqno(5.12)$$
\noi where we have used the fact that the function $f(\gamma ) = \int_{t_0}^t dt' \ t'^{-2
\gamma}$ is logarithmically convex in $\gamma$ and therefore satisfies
$$f(\gamma ) \leq f(0)^{1 - \gamma} \ f(1)^{\gamma} = (t- t_0) (t_0t)^{-\gamma} \leq t_0^{-
\gamma} \ t^{1 - \gamma} \quad . \eqno(5.13)$$
\noi Now (5.12) defines a closed subset ${\cal R}$ of ${I\hskip-1truemm R}^+ \times {I\hskip-1truemm
R}^+$ in the $(Y, Z)$ variables, containing the point $(y_0, \widetilde{z}_0)$, and $(Y, Z)$ is a
continuous function of $\bar{t}$ starting from that point for $\bar{t} = t_0$. If we can find an
open region ${\cal R}_1$ of ${I\hskip-1truemm R}^+ \times {I\hskip-1truemm R}^+$ containing $(y_0 ,
\widetilde{z}_0)$ and such that $\overline{{\cal R} \cap {\cal R}_1} \subset {\cal R}_1$, then $(Y,
Z)$ will remain in ${\cal R} \cap {\cal R}_1$ for all time, because ${\cal R} \cap {\cal R}_1 =
\overline{{\cal R} \cap {\cal R}_1}$ is both open and closed in ${\cal R}$. We take for ${\cal R}_1$
the strip ${\cal R}_1 = \{ (Y, Z) : Z <(1/2) \gamma t_0^{\gamma}\}$, so that in $\overline{{\cal R}
\cap {\cal R}_1}$
$$\left \{ \begin{array}{l} Y \leq 2 y_0 \\ \\ Z \leq 2 \widetilde{z}_0 + 2(1 - \gamma)^{-1} \
Y^2 \leq 2 \widetilde{z}_0 + 8(1 - \gamma )^{-1} \ y_0^2 \end{array} \right . \eqno(5.14)$$ 
\noi and the condition $\overline{{\cal R} \cap {\cal R}_1} \subset {\cal R}_1$ is ensured by
$$\gamma t_0^{\gamma} > 4 \widetilde{z}_0 + 16(1 - \gamma )^{-1} \ y_0^2 \quad . \eqno(5.15)$$

\noi {\bf Case ${\bf t} \leq {\bf t_0}$.} Let $\bar{t} < t_0$ and define $Y =
Y(\bar{t}) = \parallel y;L^{\infty}([\bar{t}, t_0])\parallel$ and $Z = Z(\bar{t}) =$\break \noindent
$\parallel z ; L^{\infty} ([\bar{t} , t_0])\parallel$. Then for all $t \in [\bar{t} , t_0]$
$$ \left \{ \begin{array}{l} - \partial_t y \leq t^{-2} \ YZ \\ \\ - \partial_t z \leq
t^{-2} \ Z^2 + t^{-\gamma} \ Y^2 \end{array}\right . \eqno(5.16)$$  
\noi and therefore by integration with the appropriate initial condition at $t_0$
$$\left \{ \begin{array}{l} Y \leq y_0 + t^{-1} \ YZ \\ \\ Z \leq \left (
\widetilde{z}_0 + (1 - \gamma )^{-1} \ Y^2 \right ) t_0^{1 - \gamma} + t^{-1} Z^2 \end{array}\right .
\eqno(5.17)$$
\noi so that for $Z < t/2$
$$\left \{ \begin{array}{l} Y < 2y_0 \\ \\ Z < 2 \left ( \widetilde{z}_0 + (1 - \gamma
)^{-1} \ Y^2 \right ) t_0^{1 - \gamma} < \left ( 2 \widetilde{z}_0 + 8 (1 - \gamma )^{-1}
\ y_0^2 \right ) t_0^{1 - \gamma} \quad , \end{array} \right . \eqno(5.18)$$
\noi which ensures the condition $Z < t/2$ provided
$$t > \left ( 4 \widetilde{z}_0 + 16 (1 - \gamma)^{-1} y_0^2 \right ) t_0^{1 - \gamma}
\quad . \eqno(5.19)$$

Putting back constants at appropriate places, we obtain (5.8) (5.9) from (5.14) (5.18), while
(5.15) (5.19) allow for the choice (5.6) (5.7). \par \nobreak
\hfill $\sq$ \\

\noi {\bf Remark 5.1.} The weakness of Proposition 5.1 is immediately apparent on (5.7) (5.9). In
fact, for the construction of wave operators, we want to solve the Cauchy problem for (5.1) (5.2)
for an initial time $t_0 \to \infty$. However when $t_0 \to \infty$, the interval $[T, \infty )$
where the solution is defined disappears, while the estimate for $|s|_{\ell}$ blows up. This
defect will be remedied in Proposition 5.6 below, but only for $\gamma > 1/2$. For the next
results, we shall need the following estimate. \\

\noi {\bf Lemma 5.1.} {\it Let $0 < \gamma < 1$, let $a > 0$, $b > 0$, $t_0 > 0$ and let $y$, $z$
be nonnegative continuous functions satisfying $y(t_0) = y_0$, $z(t_0) = z_0$, and}

$$\left \{ \begin{array}{l} |\partial_t y | \leq t^{-1 - \gamma} \ by + t^{-2} \
az \\ \\ |\partial_t z| \leq t^{-1 - \gamma} \ bz + t^{-\gamma}\ a y \quad . \end{array}\right .
\eqno(5.20)$$

\noi {\it Define $\bar{y}$, $\bar{z}$ by}

$$(y, z) = (\bar{y} , \bar{z}) \exp \left [ b \gamma^{-1} |t^{- \gamma} - t_0^{- \gamma} |\right ]
\quad . \eqno(5.21)$$

\noi {\it Then, for $\gamma (t_0^{\gamma} \wedge t^{\gamma}) \geq 2a^2$, the following estimates
hold~:}  

$$\left \{ \begin{array}{l} \bar{y} \leq 2 \left ( y_0 + a \ z_0 \ t_0^{-1} \right )
\\ \\ \bar{z} \leq z_0 + 2 (1 - \gamma )^{-1} a \left ( y_0 + a\ z_0 \ t_0^{-1} \right )
t^{1 - \gamma} \end{array}\right . \eqno(5.22)$$

\noi {\it for $t \geq t_0$, and}

$$\left \{ \begin{array}{l} \bar{y} \leq y_0 + 2a \left ( z_0 + (1 - \gamma)^{-1}
\ a\ y_0 \ t_0^{1 - \gamma} \right ) t^{-1} \\ \\ \bar{z} \leq 2 \left ( z_0 + (1 -
\gamma)^{-1} \ a\ y_0 \ t_0^{1 - \gamma} \right ) \end{array}\right .\eqno(5.23)$$

\noi {\it for $t \leq t_0$.} \\

\noi {\bf Proof.} The inequalities (5.20) have been written in differential form, but should be
understood in integrated form. Changing the variables from $(y, z)$ to $(\bar{y}, \bar{z})$ yields
$$\left \{ \begin{array}{l} \pm \partial_t \bar{y} \leq t^{-2} a\bar{z} \\ \\ \pm \partial_t \bar{z}
\leq t^{- \gamma} a \bar{y} \end{array}\right . \eqno(5.24)$$  
\noi for $t \displaystyle{\mathrel{\mathop >_{<}}} t_0$, and in integrated form
$$\left \{ \begin{array}{l} \bar{y} \leq y_0 + a \left | \int_{t_0}^t dt' \ t'^{-2}\ 
\bar{z} (t') \right | \\ \\ \bar{z} \leq z_0 + a \left | \int_{t_0}^t dt' \
t'^{-\gamma} \ \bar{y}(t') \right | \quad . \end{array}\right . \eqno(5.25)$$
\noi Substituting the second component of (5.25) into the first one yields an inequality for
$\bar{y}$ alone~:
$$\bar{y} \leq y_0 + a\ z_0 |t^{-1} - t_0^{-1}| + a^2 \left | \int_{t_0}^t dt' \ t'^{-\gamma}
\left ( t'^{-1} - t^{-1} \right ) \bar{y}(t') \right | \quad . \eqno(5.26)$$
\noi We shall estimate $\bar{y}$ by (5.26) and estimate $\bar{z}$ by substituting the result into
the second component of (5.25). (One could equally well obtain an equation for $\bar{z}$ alone
similar to (5.26) and estimate $\bar{z}$ directly therefrom). We distinguish the cases $t
\displaystyle{\mathrel{\mathop >_{<}}} t_0$. \\

\noi {\bf Case ${\bf t} \geq {\bf t_0}$.} Let $\bar{t} > t_0$. We define $Y =
\parallel \bar{y} ; L^{\infty} ([t_0 , \bar{t}])\parallel$. It then follows from (5.26) that
$$Y \leq y_0 + a\ z_0 \ t_0^{-1} + a^2 \gamma^{-1} \ t_0^{- \gamma} \ Y$$
\noi and therefore
$$Y \leq 2 \left ( y_0 + a\ z_0 \ t_0^{-1} \right )$$
\noi for $\gamma t_0^{\gamma} \geq 2a^2$, while by (5.25)
$$\bar{z} \leq z_0 + a(1 - \gamma )^{-1} \ t^{1 - \gamma} \ Y \quad ,$$
\noi which immediately yields (5.22). \\ 

\noi {\bf Case ${\bf t} \leq {\bf t_0}$.} Let $\bar{t} < t_0$. We define $Y =
\parallel t(\bar{y} - y_0);L^{\infty}([\bar{t} , t_0]) \parallel$. It then follows from (5.26) that
$$\bar{y} \leq y_0 + a\ z_0 \ t^{-1} + a^2 y_0 \ t^{-1} (1 - \gamma )^{-1} \ t_0^{1 - \gamma} +
a^2 Y t^{-1} \ \gamma^{-1} \ t^{-\gamma}$$
\noi so that
$$Y \leq a \ z_0 + a^2 y_0 (1 - \gamma )^{-1} \ t_0^{1 - \gamma} + a^2 \gamma^{-1} \ t^{-\gamma}
\ Y$$ \noi and therefore
$$Y \leq 2a \left ( z_0 + a\ y_0 (1 - \gamma )^{-1} \ t_0^{1 - \gamma} \right )$$
\noi for $\gamma t^{\gamma} \geq 2a^2$. Substituting that result into the second component of (5.25)
yields 
\begin{eqnarray*}
\bar{z} &\leq & z_0 + a\ y_0 (1 - \gamma)^{-1} \ t_0^{1 - \gamma} + a \gamma^{-1} \ t^{- \gamma}
\ Y \\ &\leq & 2 \left ( z_0 + a\ y_0 (1 - \gamma)^{-1} \ t_0^{1 - \gamma} \right )
\end{eqnarray*}
\noi by using again the condition $\gamma t^{\gamma} \geq 2a^2$. The previous estimates
immediately yield (5.23). \\

\noi {\bf Remark 5.2.} Since (5.20) is a linear system, it is clear that one could obtain
estimates of $y$ and $z$ valid for all $t_0 > 0$, $t > 0$, but in view of Proposition 5.1, we are
not interested in estimates for small $t_0$, $t$. The apparent lower restrictions on $t_0$, $t$
come from the fact that the ansatz $Y$ for $\bar{y}$ is inadequate for small $t$. \\

As a first consequence of Lemma 5.1, we obtain a uniqueness result at infinity for the system
(5.1) (5.2). \\

\noi {\bf Proposition 5.2.} {\it Let $0 < \gamma < 1$ and let $(k, \ell )$ be an admissible pair. Let
$(w_i, s_i)$, $i = 1,2$ be two solutions of the system (5.1) (5.2) such that $(w_i,
\widetilde{s}_i) \in L^{\infty} ([T, \infty ), H^k \oplus X^{\ell})$ for some $T > 0$ and such
that $|w_1 - w_2|_{k-1} t^{1 - \gamma}$ and $|s_1 - s_2|_{\ell - 1}$ tend to zero when $t \to
\infty$. Then $(w_1, s_1) = (w_2, s_2)$.} \\

\noi {\bf Proof.} Define
$$a = \ \mathrel{\mathop {\rm Max}_i} \parallel w_i ; L^{\infty}([T, \infty ), H^k) \parallel
\qquad , \quad b = \ \mathrel{\mathop {\rm Max}_i} \parallel \widetilde{s}_i : L^{\infty}([T,
\infty ), X^{\ell}) \parallel \quad , \eqno(5.27)$$ 
$$y = |w_1 - w_2 |_{k-1} \quad , \quad z = |s_1 - s_2|_{\ell - 1} \quad . \eqno(5.28)$$
\noi Then by Lemma 3.9, $y$ and $z$ satisfy the inequalities
$$\left \{ \begin{array}{l} |\partial_t y | \leq C \ t^{-1 - \gamma} \ by + C \
t^{-2} \ az \\ \\ |\partial_t z | \leq C \ t^{-1 - \gamma} \ bz + C \ t^{-\gamma} \
ay\end{array}\right . \eqno(5.29)$$
\noi and therefore (up to constants), the estimates (5.21) (5.23) for all $T \leq t \leq t_0$,
$t$ and $t_0$ sufficiently large. Taking the limit $t_0 \to \infty$ in (5.23) for fixed $t$ shows
that $y(t) = 0$, $z(t) = 0$. \par \nobreak
\hfill $\sq$ \\

We now begin the study of the asymptotic behaviour of solutions of the system (5.1) (5.2) by
showing that for the solutions obtained in Proposition 5.1 (actually for slightly more general
solutions, if any) $w(t)$ tends to a limit when $t \to \infty$. \\

\noi {\bf Proposition 5.3.} {\it Let $0 < \gamma < 1$ and let $(k, \ell )$ satisfy $k \leq \ell + 1$
and $\ell > n/2$. Let $(w, s)$ be a solution of the system (5.1) (5.2) such that $(w,
\widetilde{s}) \in ({\cal C} \cap L^{\infty})([T, \infty ), H^k \oplus X^{\ell})$ for some $T >
0$. Let} 
$$a = \parallel w; L^{\infty}([T, \infty ), H^k)\parallel \quad , \quad b = \parallel
\widetilde{s} ; L^{\infty}([T, \infty ), X^{\ell}) \parallel \quad . \eqno(5.30)$$

{\it Then there exists $w_{+} \in H^k$ such that $w(t)$ tends to $w_{+}$ strongly in
$H^{k-1}$ and weakly in $H^k$ when $t \to \infty$. Furthermore the following estimates hold~:} 
$$ |w_{+}|_k \leq a \eqno(5.31)$$
$$|w(t_0) - w(t)|_{k-1}   \leq C \ ab\ \gamma^{-1}(t_0 \wedge t )^{- \gamma} \eqno(5.32)$$  
$$|w(t) - w_{+}|_{k-1}   \leq C \ ab\ \gamma^{-1}\ t^{- \gamma} \eqno(5.33)$$  
\noi {\it for $t_0$, $t$ large enough, namely $\gamma (t_0 \wedge t)^{\gamma} \geq Cb$ or $\gamma
t^{\gamma} \geq Cb$.} \\

\noi {\bf Proof.} Let $T \leq t_0 \leq t$ and $w_0 = w(t_0)$. By exactly the same method as in
Lemma 3.9 (see esp. (3.54)) we obtain from (5.1)
\begin{eqnarray*}
\partial_t |w - w_0|_{k-1} &\leq & C \ t^{-2} |s|_{\ell} \left ( |w - w_0|_{k-1} + |w_0|_k \right )
\\ &\leq & C \ t^{-1 - \gamma} \ b\left ( |w - w_0|_{k-1} + a \right ) \end{eqnarray*}
\noi and by integration between $t_0$ and $t$
\begin{eqnarray*}
|w - w_0|_{k-1} &\leq& a \left ( \exp \left ( C \ b \ \gamma^{-1} \ t_0^{- \gamma} \right ) - 1
\right ) \\ &\leq & C(e - 1) ab \ \gamma^{-1} \ t_0^{- \gamma} \end{eqnarray*}
\noi for $\gamma t_0^{\gamma} \geq Cb$. This proves (5.32), from which it follows that $w(t)$ has a
limit $w_{+} \in H^{k-1}$ such that (5.33) holds. Since in addition $w(t)$ is uniformly
bounded in $H^k$, it follows by a standard compactness argument that $w_{+} \in H^k$, that
$w_{+}$ satisfies (5.31) and that $w(t)$ tends weakly to $w_{+}$ in $H^k$. \par \nobreak
\hfill $\sq$ \\

\noi {\bf Remark 5.3.} Note that the condition on $(k, \ell )$ in Proposition 5.3 is weaker than
admissibility, since in particular we do not use (5.2) and therefore do not require estimates of
$g$. \\

We now turn to the proof of existence of asymptotic states for the solutions of (5.1) (5.2)
constructed in Proposition 5.1. As explained in Section 2, we want to perform the following
construction on such a solution. We take $t_0 \geq T$ and we consider the solution $(w_0, s_{0,
t_0})$ of the free equation (5.3) which coincides with $(w, s)$ at time $t_0$. That solution is
defined by $w_0 = w(t_0)$ and
$$s_{0, t_0}(t) = s(t_0) + \int_{t_0}^t dt' \ t'^{- \gamma} \ \nabla g (w_0, w_0) \eqno(5.34)$$
\noi or equivalently, by (5.2)
$$s_{0,t_0}(t) = s(t) - \int_{t_0}^t dt' \left \{ t'^{- 2}(s \cdot \nabla ) s+ t'^{- \gamma}
\left ( \nabla g(w, w) - \nabla g (w_0, w_0) \right ) \right \} \quad . \eqno(5.35)$$

We want to prove that $(w_0, s_{0, t_0})$ has a limit when $t_0 \to \infty$. Proposition 5.3
already provides us with the limit $w_+ = \lim w(t_0)$, and it remains only to be proved
that also $s_{0,t_0}$ has a limit. Taking formally the limit $t_0 \to \infty$ in (5.35) leads us to
expect that that limit should be defined by
$$s_0(t) = s(t) + \int_t^{\infty} dt' \left \{ t'^{-2} (s \cdot \nabla ) s + t'^{- \gamma} \left (
\nabla g (w, w) - \nabla g (w_+ , w_+) \right ) \right \} \quad . \eqno(5.36)$$

>From Corollary 3.1 and from (5.34), it follows that for $(w, \widetilde{s}) \in ({\cal C} \cap
L^{\infty})([T, \infty ), H^k \oplus X^{\ell})$, also $(w_0, \widetilde{s}_{0,t_0})$ belongs to
the same space, actually with constant $w_0$, and by (5.30), satisfies the estimate
$$|s_{0,t_0}(t)|_{\ell} \leq b \ t_0^{1 - \gamma} + C \ a^2(1 - \gamma )^{-1} (t_0 \vee t)^{1 -
\gamma} \quad . \eqno(5.37)$$

That estimate however is not bounded uniformly in $t_0$, and is therefore useless to take the
limit $t_0 \to \infty$. In order to proceed further, we shall need stronger assumptions, and in
particular $\gamma > 1/2$. Under that condition and by using (5.35), we can indeed obtain
estimates on $s_{0,t_0}$ that are uniform in $t_0$. \\

\noi {\bf Proposition 5.4.} {\it Let $1/2 < \gamma < 1$ and let $(k, \ell )$ be an admissible pair.
Let $(w, s)$ be a solution of the system (5.1) (5.2) such that $(w, \widetilde{s}) \in ({\cal C} \cap
L^{\infty})([T, \infty ), H^k \oplus X^{\ell})$ for some $T \geq 1$ and define $a$ and $b$ by
(5.30). Let $t_0 \geq T$, let $w_0 = w(t_0)$ and define $s_{0, t_0}(t)$ by (5.34). Then
$s_{0,t_0}(t)$ satisfies the estimates}
$$|s_{0, t_0}(t) - s(t)|_{\ell - 1} \leq C \left ( b^2 + (1 - \gamma )^{-1} a^2 b \right ) t_0^{-
\gamma} \ t^{1 - \gamma} \eqno(5.38)$$
\noi {\it for $t \geq t_0$ and $t_0^{\gamma} \geq Cb$,}
$$|s_{0, t_0}(t) - s(t)|_{\ell - 1} \leq C (2 \gamma - 1)^{-1} \left ( b^2 + a^2b \right ) t^{1 - 2
\gamma} \eqno(5.39)$$
\noi {\it for $T \leq t \leq t_0$ and $t^{\gamma} \geq Cb$.} \par

{\it Assume in addition that $(1 - \gamma ) t_0^{\gamma} \geq Ca^2$ and $(2 \gamma - 1) t^{\gamma}
\geq C(a^2 + b)$. Then $s_{0, t_0}$ satisfies the estimate}
$$|s_{0, t_0}(t)|_{\ell - 1} \leq C \ b \ t^{1 - \gamma} \quad . \eqno(5.40)$$
\noi {\bf Proof.} By the same method as in Lemma 3.9 (see esp. (3.57)), we estimate the integrand
in (5.35) by
$$|\{ \cdot \} |_{\ell - 1} \leq C \  t^{-2} |s|_{\ell - 1} \  |s|_{\ell} + C \ t^{- \gamma} \  |w -
w_0|_{k-1} \  |w + w_0|_k \eqno(5.41)$$
\noi which by (5.30) and Proposition 5.3, esp. (5.32), can be continued as 
$$\cdots \leq C \ b^2 \ t^{-2 \gamma} + C \ a^2b \ t^{- \gamma} (t \wedge t_0)^{-\gamma}
\eqno(5.42)$$
\noi for $(t \wedge t_0)^{\gamma} \geq Cb$. Integrating (5.42) between $t_0$ and $t$ and using
(5.13) if $t \geq t_0$ yields (5.38) for $t \geq t_0$ and (5.39) for $t \leq t_0$.  \par

Finally (5.40) follows from (5.30), from (5.38) (5.39) and from the additional conditions on $t$,
$t_0$. \par \nobreak
\hfill $\sq$ \\

We next prove that $s_0(t)$ is actually well defined by (5.36) and is the limit of $s_{0,t_0}(t)$
when $t_0 \to \infty$. \\

\noi {\bf Proposition 5.5.} {\it Let $1/2 < \gamma < 1$ and let $(k, \ell )$ be an admissible pair.
Let $(w, s)$ be a solution of the system (5.1) (5.2) such that $(w, \widetilde{s}) \in ({\cal C}
\cap L^{\infty})([T, \infty ), H^k \oplus X^{\ell})$ for some $T \geq 1$ and define $a$ and $b$
by (5.30). Let $w_+$ be the limit of $w(t)$ when $t \to \infty$ obtained in
Proposition 5.3. Then} \par

{\it (1) The integral in (5.36) is absolutely convergent in $X^{\ell - 1}$ and defines a solution
$s_0$ of the equation (5.3) such that $\widetilde{s}_0 \in ({\cal C} \cap L^{\infty})([1, \infty ),
X^{\ell - 1})$ and such that the following estimate holds for $t \geq T$, $t^{\gamma} \geq
Cb$~:} $$|s_0(t) - s(t)|_{\ell - 1} \leq C (2 \gamma - 1)^{-1} \left ( b^2 + a^2b \right ) t^{1
- 2 \gamma} \quad . \eqno(5.43)$$
\noi {\it If in addition $(2 \gamma - 1) t^{\gamma} \geq C (b + a^2)$, the following estimate
holds~:} 
$$|s_0(t)|_{\ell - 1} \leq C \ b \ t^{1- \gamma} \quad . \eqno(5.44)$$

{\it (2) The function $s_{0,t_0}$ defined by (5.34) converges to $s_0$ in norm in $X^{\ell - 1}$ when
$t_0 \to \infty$ for $t \geq T$, $t^{\gamma} \geq Cb$, uniformly in compact intervals, and the
following estimate holds for $t \wedge t_0 \geq T$, $(t \wedge t_0)^{\gamma} \geq Cb$~:}
$$|s_{0,t_0}(t) - s_0(t)|_{\ell - 1} \leq C \left ( ( 1 - \gamma)^{-1} + (2 \gamma - 1)^{-1} \right
) \left ( b^2 + a^2 b \right ) t_0^{- \gamma} (t \vee t_0)^{1 - \gamma} \quad . \eqno(5.45)$$

\noi {\bf Proof.} {\bf Part (1).} By the same method as in the proof of Proposition 5.4, we
estimate the integrand in (5.36) by (5.41) with $w_0$ replaced by $w_+$ and therefore by (5.30)
(5.33)
$$|\{ \cdot \} |_{\ell - 1} \leq C \left ( b^2 + a^2 b \right ) t^{-2 \gamma} \eqno(5.46)$$
\noi which proves the convergence of the integral in $X^{\ell - 1}$ and yields the estimate (5.43).
Finally (5.44) follows from (5.30) (5.43) and from the additional condition on $t$. \\

\noi {\bf Part (2).} For $t \geq t_0$, we estimate
$$|s_{0,t_0}(t) - s_0(t)|_{\ell - 1} \leq |s_{0,t_0}(t) - s(t)|_{\ell - 1} + |s_0(t) - s(t)|_{\ell
- 1} \quad . \eqno(5.47)$$
\noi We estimate the first norm in the RHS by (5.38) and the second norm by (5.43). This yields
(5.45) for $t \geq t_0$. \par

For $t \leq t_0$, we start from 
$$s_{0,t_0}(t) - s_0(t) = s(t_0) - s_0(t_0) + \int_t^{t_0} dt' \ t'^{- \gamma} \left ( \nabla
g(w_+, w_+) - \nabla g (w (t_0), w(t_0)) \right ) \quad . \eqno(5.48)$$
\noi We estimate the first difference in the RHS by (5.43) with $t = t_0$ and the integrand by
Corollary 3.1 and (5.33) as
$$|\nabla g (w_+ , w_+) - \nabla g (w (t_0), w(t_0))|_{\ell - 1} \leq C |w(t_0) - w_+|_{k-1} \ 
 |w(t_0) + w_+|_k \leq C \ a^2 b \ t_0^{- \gamma} \eqno(5.49)$$
\noi which after integration yields (5.45) for $t \leq t_0$. \par

The convergence of $s_{0,t_0}$ to $s_0$ in the norms indicated follows from (5.45). \par
\nobreak \hfill $\sq$ \\

\noi {\bf Remark 5.4.} Note that in Proposition 5.5 we are loosing one degree of regularity in
$s$, namely a solution $(w, s)$ in $H^k \oplus X^{\ell}$ has an asymptotic free solution $(w_+ ,
s_0)$ in $H^k \oplus X^{\ell - 1}$ only. We have no uniform estimate in $t_0$ of $s_{0,t_0}$ in
$X^{\ell}$, the best we have being (5.37) and we are therefore unable to assert that the limiting
$s_0$ remains in $X^{\ell}$. \par

We now turn to the main and more difficult question of existence of the wave operator. The
construction is now the converse of that performed in Propositions 5.3-5.5. We consider a fixed
solution $(w_+ , s_0)$ of (5.3), we construct a solution $(w_{t_0}, s_{t_0})$ of (5.1) (5.2)
which coincides with $(w_+, s_0)$ at time $t_0$, and we take the limit of that solution when
$t_0 \to \infty$. Here again we shall encounter a loss of derivative, now more severe than in
the previous case. In particular we shall need to start with a free solution in $H^{k+1} \oplus
X^{\ell +1}$ for admissible $(k, \ell )$ and end up with a solution $(w, s)$ which is only in
$H^k \oplus X^{\ell}$. The crucial step is the construction of $(w_{t_0}, s_{t_0})$ and will be
performed by an extension of the energy method used in Proposition 5.1. \\

\noi {\bf Proposition 5.6.} {\it Let $1/2 < \gamma < 1$ and let $(k, \ell )$ be an admissible pair.
Let $(w_+ , s_0(1)) \in H^{k+1} \oplus X^{\ell + 1}$, define $s_0$ by (5.4) and let}
$$a = |w_+|_{k+1} \quad , \quad b = \parallel \widetilde{s}_0 ; L^{\infty} ([1, \infty
), X^{\ell + 1} ) \parallel \quad .\eqno(5.50)$$

\noi {\it Then, there exist $T_0$ and $T$, $1 \leq T_0$, $T < \infty$, depending only on $(\gamma ,
a, b)$, such that for all $t_0 \geq T_0 \vee T$, the system (5.1) (5.2) with initial data $w(t_0) =
w_+$, $s(t_0) = s_0(t_0)$, has a unique solution $(w_{t_0}, s_{t_0})$ such that $(w_{t_0},
\widetilde{s}_{t_0}) \in ({\cal C} \cap L^{\infty} ) ([T, \infty ), H^k \oplus X^{\ell})$. One can
take}
$$T_0^{\gamma} = C \left ( b + (1 - \gamma )^{-1} a^2 \right ) \eqno(5.51)$$
$$T^{\gamma} = C (2 \gamma - 1 )^{-1} \ (b + a^2) \eqno(5.52)$$
\noi {\it and the solution satisfies the estimates}
$$\left \{ \begin{array}{l} |w_{t_0}(t) - w_+|_k \leq C \ ab \ t_0^{- \gamma} \\ \\
|s_{t_0}(t) - s_0(t) |_{\ell} \leq C \left ( b^2 + (1 - \gamma )^{-1} a^2 b \right )
t_0^{- \gamma} \ t^{1 - \gamma} \end{array} \right . \eqno(5.53)$$
\noi {\it for $t \geq t_0$,}
$$\left \{ \begin{array}{l} |w_{t_0}(t) - w_+|_k \leq C \ ab \ t^{- \gamma} \\ \\
|s_{t_0}(t) - s_0(t) |_{\ell} \leq C (2 \gamma - 1 )^{-1} ( b^2 + a^2 b )
t^{1 - 2 \gamma}  \end{array} \right . \eqno(5.54)$$
\noi {\it for $T \leq t \leq t_0$, and}
$$|w_{t_0}(t)|_k \leq C \ a \quad , \quad |s_{t_0}(t)|_{\ell} \leq C \ b \ t^{1 - \gamma} 
\eqno(5.55)$$
\noi {\it for all $t \geq T$.} \\

\noi {\bf Proof.} The result follows from Proposition 4.1 and standard globalisation arguments
provided we can derive (5.53)-(5.54) as a priori estimates under the assumptions of the
proposition. Let $(w_{t_0}, s_{t_0})$ be the maximal solution of (5.1) (5.2) with the appropriate
initial condition at $t_0$ and define $y = |w_{t_0} - w_+|_k$ and $z = |s_{t_0} - s_0|_{\ell}$. We
rewrite (5.1) and the difference between (5.2) (5.3) as 
$$\hskip 1 truecm \left \{ \begin{array}{ll} \partial_t (w_{t_0} - w_+) = (2t^2)^{-1} U(1/t) \left
( 2 s_{t_0} \cdot \nabla + (\nabla \cdot s_{t_0}) \right ) U^*(1/t) w_{t_0} &\hskip 2 truecm
(5.56)\\ \\ \partial_t (s_{t_0} - s_0) = t^{-2} (s_{t_0} \cdot \nabla ) s_{t_0} + t^{-\gamma}
\left ( \nabla g (w_{t_0}, w_{t_0}) - \nabla g (w_+ , w_+) \right ) \quad .
&\hskip 2 truecm (5.57)\end{array}\right .$$ \noi From (5.56) (5.57), by exactly the same method as
in Lemma 3.8, we obtain $$\left \{ \begin{array}{l} \Big | \partial_t |w_{t_0} - w_+|_k \Big |
\leq C t^{-2} |s_{t_0}|_{\ell} \left ( |w_{t_0} - w_+|_k + |w_+|_{k+1} \right ) \\ \\
\Big | \partial_t |s_{t_0} - s_0|_{\ell} \Big | \leq C t^{-2} |s_{t_0}|_{\ell} \left ( |s_{t_0} -
s_0|_{\ell} + |s_0|_{\ell + 1} \right ) + C t^{- \gamma} |w_{t_0} - w_+|_k \ |w_{t_0} + w_+|_k
\end{array}\right . \eqno(5.58)$$
\noi and therefore by (5.50)
$$\left \{ \begin{array}{l} |\partial_t y | \leq C \ t^{-2} (z + b\ t^{1 - \gamma})
(y + a) \\ \\ |\partial_t z| \leq C \ t^{-2} (z + b \ t^{1 - \gamma} )^2 + Ct^{-
\gamma} \ y(y + a) \quad . \end{array}\right . \eqno(5.59)$$
\noi We estimate $y$ and $z$ from (5.59), taking $C = 1$ for the rest of the proof. We distinguish
again two cases. \\

\noi {\bf Case ${\bf t} \geq {\bf t_0}$.} Let $\bar{t} > t_0$ and define $Y =
Y(\bar{t}) = \parallel y; L^{\infty}([t_0, \bar{t}\ ])\parallel$ and $Z = Z(\bar{t}) =$\break
\noindent $\parallel t^{\gamma - 1} z ; L^{\infty} ([t_0, \bar{t}\ ])\parallel$. Then for all $t
\in [t_0, \bar{t}\ ]$  $$\left \{ \begin{array}{l} \partial_t y \leq t^{-1 - \gamma} (Z + b) (Y
+ a) \\ \\ \partial_t z \leq t^{-2 \gamma} (Z + b)^2 + t^{-\gamma} Y(Y + a) \end{array}\right .
\eqno(5.60)$$ \noi and therefore by integration with the appropriate initial condition at $t_0$
$$\left \{ \begin{array}{l} Y \leq 2 t_0^{- \gamma} (Z + b) (Y + a) \\ \\ Z \leq
t_0^{- \gamma} \left ( Z (Z + 2b ) + b^2 \right ) + (1 - \gamma )^{-1} Y (Y + a) \end{array}\right .
\eqno(5.61)$$
\noi where we have used the condition $\gamma > 1/2$ and (5.13) again. \par

As in the proof of Proposition 5.1, we impose an additional condition
$$2t_0^{- \gamma} (Z + 2b) < \lambda \leq 1/2 \eqno(5.62)$$
\noi and we obtain from (5.61)
$$\left \{ \begin{array}{l} Y \leq (1 - \lambda )^{-1} \lambda a \\ \\ Z \leq (1 - \lambda
)^{-1} b^2 t_0^{- \gamma} + (1 - \lambda)^{-3} (1 - \gamma)^{-1} \lambda a^2 \quad . \end{array}
\right . \eqno(5.63)$$
\noi We now choose $\lambda = 8b t_0^{- \gamma}$ so that (5.62L) reduces to $Z < 2b$, and we obtain
from (5.63)
$$\left \{ \begin{array}{l} Y \leq 16 ab\ t_0^{- \gamma} \\ \\ Z \leq \left ( 2b^2 + 64 (1 -
\gamma )^{-1} a^2b \right ) t_0^{- \gamma} \quad . \end{array} \right . \eqno(5.64)$$
 \noi The condition $Z < 2b$ is ensured by (5.51) for $t_0 \geq T_0$, and the estimate (5.53)
follows from (5.64). \\

\noi {\bf Case ${\bf t} \leq {\bf t_0}$.} Let $\bar{t} \leq t_0$ and define $Y =
Y(\bar{t}) = \parallel t^{\gamma}y;L^{\infty}([\bar{t}, t_0])\parallel$ and $Z = Z(\bar{t})
=$\break \noindent $\parallel t^{2 \gamma - 1}z;L^{\infty} ([\bar{t}, t_0])\parallel$. Then for all
$t \in [\bar{t}, t_0]$
$$\left \{ \begin{array}{l} - \partial_t y \leq t^{-1 - \gamma} (Z t^{- \gamma} +
b) (Yt^{-\gamma} + a) \\ \\ - \partial_t z \leq t^{-2 \gamma} (Z t^{- \gamma} + b)^2
+ t^{-2 \gamma} Y(Yt^{- \gamma} + a)\end{array}\right . \eqno(5.65)$$
\noi and therefore by integration with the appropriate initial condition at $t_0$
$$\left \{ \begin{array}{l} y \leq 2abt^{- \gamma} + (aZ + bY)t^{-2 \gamma} + YZ
t^{-3 \gamma} \\ \\ z \leq (2 \gamma - 1)^{-1} (b^2 + aY) t^{1 - 2 \gamma} + 2(Y^2 +
2bZ)t^{1 - 3 \gamma} + Z^2 t^{1 - 4 \gamma} \end{array}\right . \eqno(5.66)$$
\noi where we have used the condition $\gamma > 1/2$, so that
$$\left \{ \begin{array}{l} Y \leq 2ab + aZt^{-\gamma} + (b + Zt^{- \gamma} )
Yt^{- \gamma} \\ \\ Z \leq (2 \gamma - 1)^{-1} (b^2 + a Y) + 2Y^2t^{- \gamma} + (4b +
Zt^{- \gamma}) Zt^{- \gamma} \quad .\end{array}\right . \eqno(5.67)$$

We impose the additional conditions $Zt^{-\gamma} < 2b$ and $t^{\gamma} \geq 12b$ and we obtain
from (5.67)
$$\left \{ \begin{array}{l} Y \leq 8ab \\ \\ Z \leq 2 \left ( (2 \gamma - 1)^{-1}
(b^2 + aY) + (6b)^{-1} Y^2 \right )\end{array}\right . \eqno(5.68)$$
\noi which yields (5.54), while the conditions $Zt^{-\gamma} < 2b$ and $t^{\gamma} \geq 12b$ are
ensured by (5.52). \par

Finally (5.55) follows from (5.50)-(5.54). \par \nobreak
\hfill $\sq$ \\

\noi {\bf Remark 5.5.} The major improvement of Proposition 5.6 over Proposition 5.1, which has
been achieved by estimating the difference of $(w, s)$ and $(w_+, s_0)$ instead of estimating
$(w, s)$ alone, is that now $(w, s)$ is defined in a time interval $[T, \infty )$ which is
independent of $t_0$, and is estimated in that interval uniformly with respect to $t_0$. The
condition $\gamma > 1/2$ plays a crucial role in obtaining that improvement (cf Remark 5.1
above). \\

\noi {\bf Remark 5.6.} One may wonder whether the loss of regularity from $(k + 1, \ell + 1)$ to
$(k, \ell )$ when constructing $(w_{t_0}, s_{t_0})$ from $(w_+, s_0)$ is unavoidable. Actually
$(w_{t_0}(t_0), s_{t_0}(t_0)) \in H^{k+1} \oplus X^{\ell + 1}$ and by regularity (Proposition 4.1
part (2)), $(w_{t_0}, s_{t_0}) \in {\cal C} ([T, \infty ), H^{k+1} \oplus X^{\ell + 1})$. One can
then estimate $y = |w_{t_0}|_{k+1}$ and $z = |s_{t_0}|_{\ell + 1}$. By Lemma 3.8 and by (5.1)
(5.2) (5.55), $y$ and $z$ satisfy the system (5.29) and therefore are estimated up to constants
by (5.22) (5.23). However under natural assumptions $y_0 = O(1)$ and $z_0 = O(t_0^{1 - \gamma})$,
the estimate (5.23), which is the one really relevant for large $t_0$, is not uniform in $t_0$,
so that the estimate is lost when one takes the limit $t_0 \to \infty$, which is what we shall do
next. As a consequence the regularity at the level of $(k+1, \ell +1)$ is also lost in that
limit, and we have therefore made no effort to keep track of it at the stage of Proposition 5.6.
\\

We can now take the limit $t_0 \to \infty$ of the solution $(w_{t_0}, s_{t_0})$ constructed in
Proposition 5.6 for fixed $(w_+, s_0)$. \\

\noi {\bf Proposition 5.7.}  {\it Let $1/2 < \gamma < 1$ and let $(k, \ell )$ be an admissible
pair.} \par

{\it (1) Let $(w_+, s_0(1)) \in H^{k+1} \oplus X^{\ell + 1}$ and define $s_0$, $a$, $b$, by (5.4)
(5.50). Then there exists $T$, $1 \leq T < \infty$, depending only on $(\gamma , a, b)$ and a unique
solution $(w, s)$ of the system (5.1) (5.2) such that $(w, \widetilde{s}) \in ({\cal C} \cap
L^{\infty}) ([T, \infty ), H^k \oplus X^{\ell})$, satisfying (5.36) and such that the following
estimates hold for all $t \geq T$~:}
$$\left \{ \begin{array}{l} |w(t) - w_+|_k \ \leq \ C \ ab \ t^{- \gamma} \\ \\ |s(t) -
s_0(t)|_{\ell} \ \leq \ C(2 \gamma - 1)^{-1} (b^2 + a^2b) t^{1 -2\gamma} \quad , \end{array}\right .
\eqno(5.69)$$   
$$|w(t)|_k \ \leq \ C \ a \qquad , \quad |s(t)|_{\ell} \ \leq \ C \ b \ t^{1 - \gamma} \quad .
\eqno(5.70)$$ \noi {\it One can take}
$$T^{\gamma} = C (2 \gamma - 1)^{-1} (b+a^2) \quad . \eqno(5.71)$$
{\it (2) Let $(w_{t_0}, s_{t_0})$ be the solution of the system (5.1) (5.2) constructed in
Proposition 5.6 for $t_0 \geq T_0 \vee T$ and in particular such that $(w_{t_0},
\widetilde{s}_{t_0}) \in ({\cal C} \cap L^{\infty}) ([T, \infty ), H^k \oplus X^{\ell})$. Then
$(w_{t_0}, s_{t_0})$ converges to $(w, s)$ in norm in $L^{\infty} (J, H^{k-1} \oplus X^{\ell -
1})$ and in the weak-$*$ sense in $L^{\infty}(J, H^k \oplus X^{\ell})$ for any compact interval
$J \subset [T, \infty )$, and in the weak-$*$ sense in $H^k \oplus X^{\ell}$ pointwise in
$t$.}\par

{\it (3) The map $(w_+ , s_0(1)) \mapsto (w, s)$ defined in Part (1) is continuous on the
bounded sets of $H^{k+1} \oplus X^{\ell + 1}$ from the norm topology of $(w_+ , s_0(1))$ in
$H^{k-1} \oplus X^{\ell - 1}$ to the norm topology of $(w, s)$ in $L^{\infty}(J, H^{k-1} \oplus
X^{\ell - 1})$ and to the weak-$*$ topology in $L^{\infty}(J, H^k \oplus X^{\ell})$ for any
compact interval $J \subset [T, \infty )$, and to the weak-$*$ topology in $H^k \oplus X^{\ell}$
pointwise in $t$.} \\

\noi {\bf Proof.} {\bf Parts (1) (2)} will follow from the convergence of $(w_{t_0}, s_{t_0})$
when $t_0 \to \infty$ in the topologies stated in Part (2). We recall that $(w_{t_0}, s_{t_0})$
satisfies the estimates (5.53) (5.54) which we rewrite more briefly as
$$\left \{ \begin{array}{l} [w_{t_0} - w_+|_k \ \leq \ M_1 (t \wedge t_0)^{-\gamma} \\ \\ |s_{t_0}
- s_0|_{\ell} \ \leq \ M_2 (t \wedge t_0)^{-\gamma} \ t^{1 - \gamma} \end{array} \right .
\eqno(5.72)$$
\noi where $M_1$, $M_2$ depend on $(\gamma , a, b)$ and can be read from (5.53) (5.54), and
satisfies the estimate (5.55). \par

Let now $T_0 \vee T \leq t_0 \leq t_1$. From (5.72) it follows that for all $t \geq t_0$
$$\left \{ \begin{array}{l} |w_{t_0} - w_{t_1}|_k \ \leq \ 2 M_1 \ t_0^{- \gamma} \\ \\ |s_{t_0}
- s_{t_1}|_{\ell} \ \leq \ 2 M_2 \ t_0^{- \gamma} \ t^{1 - \gamma} \quad . \end{array} \right .
\eqno(5.73)$$
\noi We now estimate $(w_{t_0} - w_{t_1}, s_{t_0} - s_{t_1})$ in $H^{k-1} \oplus X^{\ell - 1}$ for
$t \leq t_0$. Let
$$y = |w_{t_0} - w_{t_1}|_{k-1} \qquad , \quad z = |s_{t_0} - s_{t_1}|_{\ell - 1} \quad .
\eqno(5.74)$$
\noi From Lemma 3.9 and from (5.55), it follows that $y$ and $z$ satisfy the system (5.29).
Integrating that system for $t \leq t_0$ with initial data at $t_0$ estimated by (5.73), we obtain
from Lemma 5.1, esp. (5.23)
$$\left \{ \begin{array}{l} y(t) \leq M \left ( t_0^{- \gamma} + t_0^{1 - 2 \gamma} \ t^{-1}
\right ) \\ \\ z(t) \leq M \ t_0^{1 - 2 \gamma} \end{array} \right . \eqno(5.75)$$
\noi for some $M$ depending on $(\gamma ,a , b)$ and for $T \leq t \leq t_0$. From (5.75) it
follows that there exists $(w, s) \in {\cal C}([T, \infty ), H^{k-1} \oplus X^{\ell - 1})$ such
that $(w_{t_0}, s_{t_0})$ converges to $(w, s)$ in $L^{\infty}(J, H^{k-1} \oplus X^{\ell - 1})$ for
all compact intervals $J \subset [T, \infty )$. From that convergence, from (5.54) (5.55) and from
standard compactness arguments, it follows that $(w, \widetilde{s}) \in (L^{\infty} \cap {\cal
C}_{w^*})([T, \infty ), H^k \oplus X^{\ell})$, that $(w, s)$ satisfies the estimates (5.69) (5.70)
for all $t \geq T$, and that $(w_{t_0}, s_{t_0})$ converges to $(w, s)$ in the other topologies
considered in Part (2). Furthermore, by the local result of Proposition 4.1, part (1), $(w, s)
\in {\cal C}([T, \infty ), H^k \oplus X^{\ell})$. Obviously $(w, s)$ is a solution of (5.1)
(5.2). We now prove that $(w, s)$ satisfies (5.36). Let $s'_0(t)$ be the RHS of (5.36). By
Proposition 5.5, part (1), $s'_0(t)$ is well defined and satisfies the analogue of (5.43).
Furthermore $s'_0(t)$ satisfies (5.3) so that $s_0(t) - s'_0(t)$ is constant in time. By (5.43)
for $s'_0(t)$ and (5.69), that constant is zero, namely $s'_0(t) = s_0(t)$. This proves (5.36).
\par

>From the uniqueness result of Proposition 5.2, it follows that $(w, s)$ is unique under the
condition (5.69). This completes the proof of Parts (1) and (2). \\

\noi {\bf Part (3).} Let $(w_+, s_0(1))$ and $(w'_+, s'_0(1))$ belong to $H^{k+1} \oplus
X^{\ell + 1}$, define $s_0$ and $s'_0$ by (5.4) and its analogue, assume that 
$$|w_+|_{k+1} \vee |w'_+|_{k+1} \leq a$$
$$\parallel t^{\gamma - 1} s_0 ; L^{\infty} ([1, \infty ), X^{\ell + 1})\parallel \vee \parallel
t^{\gamma - 1} s'_0; L^{\infty}([1, \infty ), X^{\ell + 1})\parallel \ \leq b \quad ,$$
\noi let $(w, s)$ and $(w', s')$ be the associated solutions of the system (5.1) (5.2)
constructed in Part (1), satisfying (5.69) and its analogue with the same $(a, b)$ for $t \geq
T$ with the same $T$ defined by (5.71). We take $(w'_+, s'_0(1))$ close to $(w_+, s_0(1))$ in
$H^{k-1} \oplus X^{\ell - 1}$ in the sense that for some small $\varepsilon$, $\varepsilon_0 > 0$
$$|w_+ - w'_+|_{k-1} \leq \varepsilon a \eqno(5.76)$$
$$|s_0(1) - s'_0(1)|_{\ell - 1} \leq \varepsilon_0 b \eqno(5.77)$$
\noi and therefore by (5.4) and estimates from Corollary 3.1
$$|s_0 - s'_0|_{\ell - 1} \leq \varepsilon_0 b + C(1 - \gamma)^{-1} \ a^2 \ \varepsilon \ t^{1 -
\gamma} \quad . \eqno(5.78)$$
\noi From (5.69) and its analogue for $(w', s')$ and from (5.76) (5.78), it follows that for all
$t \geq T$ $$\left \{ \begin{array}{l} |w - w'|_{k-1} \leq \varepsilon a + C \ ab \ t^{-\gamma}
\\ \\ |s - s'|_{\ell - 1} \leq \varepsilon_0b + C(1 - \gamma )^{-1} \ a^2 \ \varepsilon \ t^{1 -
\gamma} + C(2 \gamma - 1)^{-1} (b^2 + a^2b)t^{1 - 2 \gamma} \quad . \end{array} \right .
\eqno(5.79)$$
\noi We now define $t_0$ by $t_0^{- \gamma} = \varepsilon b^{-1}$ so that for $t \geq t_0$, 
$$|w - w'|_{k-1} \leq C \ \varepsilon \ a \eqno(5.80)$$
$$|s - s'|_{\ell - 1} \leq \varepsilon_0 b + C \left [ (1 - \gamma )^{-1} a^2 + (2 \gamma -
1)^{-1} (b + a^2) \right ] \varepsilon t^{1 - \gamma}$$
\noi and in particular
$$|s(t_0) - s'(t_0)|_{\ell - 1} \leq \varepsilon_0 b + M \ \varepsilon^{(2\gamma - 1)/\gamma}
\eqno(5.81)$$
\noi for some $M$ depending on $(\gamma , a, b)$. We now estimate $(w - w' , s - s')$ in $H^{k-1}
\oplus X^{\ell - 1}$ for $t \leq t_0$. Let 
$$y = |w - w'|_{k-1} \qquad , \quad z = |s - s'|_{\ell - 1} \quad .$$
\noi From Lemma 3.9, it follows that $y$ and $z$ satisfy the system (5.29). Integrating that
system for $t \leq t_0$ with initial data at $t_0$ estimated by (5.80) (5.81), we obtain from
Lemma 5.1, esp. (5.23),
$$\left \{ \begin{array}{l} y(t) \leq C \ \varepsilon \ a + \left ( C \ \varepsilon_0 \ b + M \
\varepsilon^{(2 \gamma - 1)/\gamma} \right ) a t^{-1} \\ \\ z(t) \leq C \ \varepsilon_0 \ b + M \
\varepsilon^{(2 \gamma - 1)/\gamma} \end{array} \right . \eqno(5.82)$$
\noi for some $M$ depending on ($\gamma , a, b)$ and for all $t$, $T \leq t \leq t_0$. This
implies the continuity of $(w, s)$ as a function of $(w_+, s_0(1))$ in the norm topology of
$L^{\infty}(J, H^{k-1} \oplus X^{\ell - 1})$ for all compact intervals $J \subset [T, \infty )$.
The other continuities follow therefrom and from the boundedness of $(w, \widetilde{s})$ in
$L^{\infty}([T, \infty ), H^k \oplus X^{\ell})$ by standard compactness arguments. \par \nobreak
\hfill $\sq$ \\

\noi {\bf Remark 5.7.} By analogy with Proposition 4.1, part (3), one expects the map $(w_+,
s_0(1)) \to (w, s)$ to be also continuous from the norm topology in $H^k \oplus X^{\ell}$ to
the norm topology in $L^{\infty}(J, H^k \oplus X^{\ell})$ for compact $J \subset [T, \infty )$.
A proof of that fact would require a combination of Steps 6 and 7 in the proof of Proposition
4.1 with the proof just given of Proposition 5.7, part (3) with however estimates at the level
$(k, \ell )$ instead of $(k - 1, \ell - 1)$. However, the coupling between $t_0$ and
$\varepsilon$ in the latter has the effect that, when the time dependence is taken into account
in the former, the resulting dependence of the estimates on $\varepsilon$ or on $t_0$ is not
sufficiently good to establish the result without additional assumptions on $\gamma$, namely
without assuming $\gamma$ sufficiently close to 1. Since the argument is rather complicated for
a result of restricted validity, we refrain from pushing it any further. \\

In Propositions 5.5 and 5.7, we have defined two maps $(w, s) \to (w_+, s_0)$ and $(w_+, s_0)
\to (w, s)$ between solutions $(w, s)$ of the system (5.1) (5.2) and solutions of the auxiliary
free equation (5.3) defined in a neighborhood of infinity in time. Both of these maps suffer
from the loss of one derivative. Nevertheless, they are inverse of each other (and in particular
injective) whenever they can be applied successively. We state that fact in the following
proposition. \\

\noi {\bf Proposition 5.8.} {\it Let $1/2 < \gamma < 1$ and let $(k, \ell)$ be an admissible pair.} 
\par

{\it (1) Let $(w_+, s_0)$ be a solution of (5.3) such that $(w_+, \widetilde{s}_0) \in ({\cal C} \cap
L^{\infty}) ([1, \infty ), H^{k+1} \oplus X^{\ell + 1})$ or equivalently defined by (5.4) with
$(w_+, s_0(1)) \in H^{k+1} \oplus X^{\ell + 1}$. Let $(w, s)$ be the solution of (5.1) (5.2)
defined in Proposition 5.7, part (1), and let $(w'_+, s'_0)$ be the solution of (5.3) defined
from $(w, s)$ in Propositions 5.3 and 5.5. Then $w'_+ = w_+$ and $s'_0 = s_0$.} \par

{\it (2) Let $(w, s)$ be a solution of (5.1) (5.2) such that $(w, \widetilde{s}) \in ({\cal C} \cap
L^{\infty}) ([T, \infty ), H^{k+2} \oplus X^{\ell + 2})$ for some $T$, $1 \leq T < \infty$. Let
$(w_+, s_0)$ be the solution of (5.3) defined in Propositions 5.3 and 5.5, so that $(w_+,
\widetilde{s}_0) \in ({\cal C} \cap L^{\infty})([1, \infty ), H^{k+1} \oplus X^{\ell + 1})$ and
let $(w',s')$ be the solution of (5.1) (5.2) defined from $(w_+, s_0)$ in Proposition 5.7, part
(1), so that $(w', \widetilde{s}') \in ({\cal C} \cap L^{\infty})([T', \infty ), H^k \oplus
X^{\ell})$ for some $T'$, $1 \leq T' < \infty$. Then $w' = w$ and $s' = s$ for $t \geq T \vee
T'$.} \\

\noi {\bf Proof.} {\bf Part (1).} From Proposition 5.7, esp. (5.69) we obtain 
$$|w - w_+|_k \leq M \ t^{- \gamma} \qquad , \quad |s - s_0|_{\ell} \leq M \ t^{1 - 2 \gamma}$$
\noi for some $M$ depending on $(w_+, s_0)$. From Propositions 5.3 and 5.5, esp. (5.33) and
(5.43), we obtain similarly
$$|w - w'_+]_{k-1} \leq M \ t^{- \gamma} \quad , \quad |s - s'_0|_{\ell - 1} \leq M \ t^{1 - 2
\gamma} \quad .$$
\noi Taking the limit $t \to \infty$ shows that $w'_+ = w_+$, so that by (5.4) $s'_0 - s_0$ is
constant in time, and therefore zero. \\

\noi {\bf Part (2).} From Propositions 5.3 and 5.5, we obtain in the same way 
$$|w - w_+|_{k+1} \leq M \ t^{-\gamma} \quad , \quad |s - s_0|_{\ell + 1} \leq M \ t^{1 - 2
\gamma} \quad .$$
\noi From Proposition 5.7, we then obtain
$$|w'- w_+|_k \leq M \ t^{- \gamma} \quad , \quad |s' - s_0|_{\ell} \leq M \ t^{1 - 2 \gamma}$$
\noi so that
$$|w - w'|_k \leq M \ t^{-\gamma} \quad , \quad |s - s'|_{\ell} \leq M \ t^{1 - 2 \gamma}$$
\noi for $t \geq T \vee T'$. The result then follows from Proposition 5.2. \par \nobreak
\hfill $\sq$ \\

We conclude this section with a brief discussion of the other systems of equations that can be
used instead of (5.1) (5.2) (5.3) and on their drawbacks as compared with the latter. That
discussion was briefly sketched at the end of Section 2 and can now be resumed at a more
technical level. \par
(i) Instead of the system (5.1) (5.2) corresponding to (2.25) (2.26) and to the choice $w =
w_3$, we could have used a system corresponding to (2.19) (2.20) and to the choice $w = w_2$,
with an additional term $(2t^2)^{-1} \Delta w$ in the equation for $w$. That term would make
no difference in the energy estimate (3.40) of Lemma 3.7 and therefore in Proposition 5.1.
However it would produce an additional term $Ct^{-2}|w_+|_{k+2}$ in the energy estimate (5.58)
for $|w - w_+|_k$, and therefore a loss of two derivatives instead of one on $w$ in the
crucial Proposition 5.6. \par

(ii) Instead of the free auxiliary equation (5.3) corresponding to (3.7), we could have used
the more accurate HJ equation corresponding to (3.6). Since however we need the assumption
$\gamma > 1/2$ in a crucial way in order to obtain estimates uniform in $t_0$ for $t \leq
t_0$, both in Proposition 5.4 (see esp. (5.39)) and in Proposition 5.6 (see esp. (5.54)) using
(3.6) instead of (3.7) would not produce any improvement of the results. On the other hand it
would make the treatment of $s_0$ more complicated, since the equation (3.6) is hardly
simpler than the system (3.4) (3.5). In particular, instead of the explicit solution (5.4), we
would need a proposition similar to Proposition 5.1 in order to solve (3.6) for $s_0$, with a
result valid for large time only. Similarly, in order to estimate $s_0(t)$ uniformly in $t_0$
for $t \leq t_0$ and to take the limit $t_0 \to \infty$ so as to prove the existence of
asymptotic states, we would have to replace the relatively simple Propositions 5.4 and 5.5 by
more complicated ones of the same degree of complication as needed to prove the existence of
wave operators, namely Propositions 5.6 and 5.7. Finally it would no longer be possible, in
the definition of the wave operators, to characterize the asymptotic solution $(w_+, s_0)$ by
an initial condition at a fixed time $t = 1$ independent of $w_+$.

\section{The auxiliary system at infinite time. Asymptotics~II.}
\hspace*{\parindent}
In this section, we perform a construction similar to that of Section 5, and we essentially
construct local wave operators at infinity for the auxiliary system (5.1) (5.2), now however
compared with the auxiliary free equation (3.8) which is better suited than (3.7)$\equiv$(5.3) for
the study of gauge invariance, as is explained in Section 2. We recall that $g_0$ and $g$ are
defined by (3.1) (3.2). Most of the results of this section will be obtained from those of Section
5 and will involve a comparison of solutions of (5.3) and (3.8). In order to make that comparison
more transparent, we shall use exclusively $g_0$ and refrain from using $g$ in this and the next
section. For brevity we shall also use the short hand notation $g_0(w) = g_0(w, w)$ for the
diagonal restriction of $g_0$. With that notation
$$g(w, w) = g_0 \left ( U^*(1/t)w \right ) \quad .$$

In order to distinguish solutions of (3.8) from those of (5.3) considered in the previous section,
we shall use the notation $s_{02}$ for the former, the additional subscript 2 referring to the
fact that the nonlinearity in (3.8) is that of the free auxiliary equation (2.22) naturally
associated with the system (2.19) (2.20) for $(w_2, \varphi_2)$. With the previous notation, the
equation (3.8) is rewritten as
$$\partial_t s_{02} = t^{-\gamma} \ \nabla \ g_0 (w_+) \eqno(6.1)$$
\noi and is trivially solved by
$$s_{02}(t) = s_{02}(1) + (1 - \gamma)^{-1} \ (t^{1 - \gamma} - 1) \ \nabla g_0 (w_+) \eqno(6.2)$$
\noi to be compared with the general solution (5.4) of (5.3). It follows from (6.2) and Corollary
3.1 that for admissible $(k , \ell )$ and for $(w_+, s_{02}(1)) \in H^k \oplus X^{\ell}$, $
\widetilde{s}_{02} \equiv t^{\gamma - 1} s_{02} \in ({\cal C} \cap L^{\infty})([1, \infty ),
X^{\ell})$ and
$$\parallel \widetilde{s}_{02} ; L^{\infty}([1, \infty ), X^{\ell})\parallel \ \leq \ |
s_{02}(1)|_{\ell} + C (1 - \gamma )^{-1} \ |w_+|_k^2 \quad . $$

We first compare solutions of (5.3) and (6.1) which coincide in a suitable sense at some time
$t_0$. \\

\noi {\bf Lemma 6.1.} {\it Let $1/2 < \gamma < 1$ and let $(k, \ell )$ be an admissible pair. Let
$w_0 \in H^k$. Let $1 \leq t_0 \leq \infty$. Let $s_0$ be a solution of (5.3) with $w_+ = w_0$ and
let $s_{02}$ be a solution of (6.1) with $w_+ = U^*(1/t_0) w_0$ such that $s_0 (t_0) =
s_{02}(t_0)$, so that} 
$$s_0(t) - s_{02}(t) = \int_{t_0}^t dt' \ t'^{-\gamma} \left ( \nabla g_0 (U^*(1/t')w_0) - \nabla
g_0 (U^*(1/t_0)w_0) \right ) \quad . \eqno(6.3)$$
\noi {\it (For $t_0 = \infty$, coincidence at $t_0$ is defined by (6.3) and justified by the
estimates to follow). Then} \par
{\it (1) The following estimates hold~:}
$$|s_0(t) - s_{02}(t)|_{\ell - 1} \leq C(1 - \gamma )^{-1}\  |w_0|_k^2 \ t_0^{-1/2} \ t^{1 - \gamma}
\quad (t \geq t_0) \eqno(6.4)$$
$$|s_0(t) - s_{02}(t)|_{\ell - 1} \leq C (2 \gamma - 1)^{-1}\  |w_0|_k^2 \ t^{1/2 - \gamma} \quad (t
\leq t_0) \quad . \eqno(6.5)$$

{\it (2) Assume in addition that $w_0 \in H^{k+1}$. Then the following estimates hold~:}
$$|s_0(t) - s_{02}(t)|_{\ell + 1} \leq C(1 - \gamma )^{-1}\  |w_0|_{k+1}^2 \ t_0^{-1/2} \ t^{1 -
\gamma} \quad (t \geq t_0) \eqno(6.6)$$
$$|s_0(t) - s_{02}(t)|_{\ell + 1} \leq C (2 \gamma - 1)^{-1} \ |w_0|_{k+1}^2 \ t^{1/2 - \gamma} \quad
(t \leq t_0) \quad . \eqno(6.7)$$

\noi {\bf Proof.} {\bf Part (1).} From (6.3) and from Corollary 3.1, we obtain 
$$|s_0(t) - s_{02}(t)|_{\ell - 1} \leq C \int_{t_0}^t dt' \ t'^{-\gamma} \Big |(U^*(1/t') -
U^*(1/t_0)) w_0 \Big |_{k-1} \ |w_0|_k \quad . \eqno(6.8)$$
\noi Now 
$$\Big | \left ( U^*(1/t') - U^*(1/t_0) \right ) w_0 \Big |_{k-1} \leq |1/t' - 1/t_0 |^{1/2} \
|w_0|_k \quad . \eqno(6.9)$$
\noi Substituting (6.9) into (6.8) and integrating over time yields (6.4) (6.5). \\
\noi {\bf Part (2).} We estimate similarly by (6.3) and Corollary 3.1
$$|s_0(t) - s_{02}(t)|_{\ell + 1} \leq C \int_{t_0}^t dt' \ t'^{- \gamma} \big | \left ( U^*(1/t') -
U^*(1/t_0) \right ) w_0 \Big |_k \ |w_0|_{k+1} \eqno(6.10)$$
\noi from which (6.6) (6.7) follow as previously. \par \nobreak
\hfill $\sq$ \\

We now follow step by step the constructions performed in Section 5 with $s_0$, leading to the
existence of asymptotic states (Propositions 5.4 and 5.5) and to the existence of local wave
operators at infinity (Propositions 5.6 and 5.7). We first consider a fixed solution $(w, s)$ of
the system (5.1) (5.2) as constructed in Proposition 5.1 and we look for a solution of (6.1) which is
asymptotic to $s(t)$ at infinity. For that purpose we take some $t_0$ large enough and we define
the solution of (6.1) which coincides with $s(t)$ at $t_0$ by
$$s_{02,t_0}(t) = s(t_0) + \int_{t_0}^t dt' \ t'^{-\gamma} \ \nabla g_0 (U^*(1/t_0) w_0)
\eqno(6.11)$$
\noi or equivalently
$$s_{02,t_0}(t) = s(t) - \int_{t_0}^t dt' \left \{  t'^{-2} (s \cdot \nabla ) s + t'^{-\gamma}
\left (  \nabla g_0 (U^*(1/t') w) - \nabla g_0(U^*(1/t_0) w_0) \right ) \right \} \quad
\eqno(6.12)$$
\noi with $w_0 = w(t_0)$ (compare with (5.34) (5.35)). As in Section 5, $s_{02,t_0}(t)$ satisfies
the analogue of the estimate (5.37), which is not uniform in $t_0$, and the first task is to obtain
an estimate uniform in $t_0$. This is done in the following proposition, which is the analogue of
Proposition 5.4.\\

\noi {\bf Proposition 6.1.} {\it Let $1/2 < \gamma < 1$ and let $(k, \ell )$ be an admissible pair.
Let $(w, s)$ be a solution of the system (5.1) (5.2) such that $(w, \widetilde{s}) \in ({\cal C}
\cap L^{\infty})([T, \infty), H^k \oplus X^{\ell})$ for some $T > 0$ and define $a$ and $b$ by
(5.30). Let $t_0 \geq T$ and $w_0 = w(t_0)$. Then $s_{02,t_0}$ defined by (6.11) satisfies the
estimates}
$$|s_{02,t_0}(t) - s(t)|_{\ell - 1} \leq C \left ( ( b^2 + (1 - \gamma )^{-1} a^2b ) t_0^{-
\gamma} + (1 - \gamma )^{-1} \ a^2 \ t_0^{-1/2}  \right ) t^{1 - \gamma} \eqno(6.13)$$ 
\noi {\it for $t \geq t_0$ and $t_0^{\gamma} \geq Cb$,}
$$|s_{02,t_0}(t) - s(t)|_{\ell - 1} \leq C (2 \gamma - 1)^{-1} \left ( ( b^2 + a^2b )
t^{1 - 2 \gamma} + a^2\ t^{1/2 - \gamma}  \right ) \eqno(6.14)$$
\noi {\it for $T \leq t \leq t_0$ and $t^{\gamma} \geq Cb$.} \par
{\it Assume in addition that $(1 - \gamma) t_0^{\gamma} \geq C a^2$, $(1 - \gamma ) t_0^{1/2}
\geq C \ a^2 \ b^{-1}$, $(2 \gamma - 1) t^{\gamma} \geq C (a^2 + b)$ and $(2 \gamma - 1) t^{1/2}
\geq C \ a^2 \ b^{-1}$. Then the following estimate holds~:}
$$|s_{02,t_0}(t)|_{\ell - 1} \leq C \ b \ t^{1 - \gamma} \quad . \eqno(6.15)$$

\noi {\bf Proof.} The result follows from Proposition 5.4 and Lemma 6.1, part (1), applied to
$s_{02,t_0}$ and $s_{0,t_0}$ defined by (6.11) and (5.34). \par \nobreak
\hfill $\sq$ \\

We now want to prove that $s_{02,t_0}$ has a limit when $t_0 \to \infty$. Following the method of
Section 5 would lead us to define the limiting function $s_{02}$ by taking the formal limit $t_0
\to \infty$ in (6.12), namely
$$s_{02}(t) = s(t) + \int_t^{\infty} dt' \left \{ t'^{-2} (s \cdot \nabla ) s + t'^{- \gamma}
\left ( \nabla g_0 (U^*(1/t') w) - \nabla g_0 (w_+) \right ) \right \} \eqno(6.16)$$
\noi where $w_+$ is the limit of $w(t)$ as $t \to \infty$ obtained in Proposition 5.3. It is
however simpler to take advantage of the results of Section 5, esp. Proposition 5.5 and to
define $s_{02}(t)$ in terms of $s_0(t)$ obtained in the latter and defined by (5.36), namely to
define $s_{02}(t)$ by
$$s_{02}(t) = s_0(t) + \int_t^{\infty} dt' \ t'^{-\gamma} \left ( \nabla g_0 (U^*(1/t') w_+) -
\nabla g_0(w_+) \right ) \quad . \eqno(6.17)$$
\noi The following proposition is the analogue for $s_{02}$ of Proposition 5.5. \\

\noi {\bf Proposition 6.2.} {\it Let $1/2 < \gamma < 1$ and let $(k, \ell )$ be an admissible pair.
Let $(w, s)$ be a solution of the system (5.1) (5.2) such that $(w, \widetilde{s}) \in ({\cal C}
\cap L^{\infty})([T, \infty), H^k \oplus X^{\ell})$ for some $T \geq 1$ and define $a$ and $b$ by
(5.30). Let $w_+$ be the limit of $w(t)$ when $t \to \infty$ obtained in Proposition 5.3. Then}
\par
{\it (1) The integral in (6.16) is absolutely convergent in $X^{\ell - 1}$ and defines a
solution $s_{02}$ of the equation (6.1) such that $\widetilde{s}_{02} \in ({\cal C} \cap
L^{\infty})([1, \infty ), X^{\ell - 1})$ and such that the following estimate holds for $t \geq
T$, $t^{\gamma} \geq Cb$~:} 
$$|s_{02}(t) - s(t) |_{\ell - 1} \leq C (2 \gamma - 1)^{-1} \left ( ( b^2 + a^2b) t^{1 - 2
\gamma} + a^2 t^{1/2 - \gamma} \right ) \quad . \eqno(6.18)$$
\noi {\it If in addition $(2 \gamma - 1) t^{\gamma} \geq C (b + a^2)$ and $(2 \gamma - 1)
t^{1/2} \geq C a^2 b^{-1}$, the following estimate holds~:}
$$|s_{02}(t)|_{\ell - 1} \leq C \ b \ t^{1 - \gamma} \quad . \eqno(6.19)$$

{\it (2) The function $s_{02,t_0}$ defined by (6.11) converges to $s_{02}$ in norm in $X^{\ell -
1}$ when $t_0 \to \infty$ for $t \geq T$, $t^{\gamma} \geq Cb$, uniformly in compact intervals,
and the following estimate holds for $t \wedge t_0 \geq T, (t \wedge t_0)^{\gamma} \geq Cb$~:}
$$|s_{02,t_0}(t) - s_{02}(t)|_{\ell - 1} \leq C \left ( (1 - \gamma)^{-1} + (2 \gamma - 1)^{-1}
\right ) \left ( (b^2 + a^2 b) t_0^{- \gamma} + a^2 \ t_0^{-1/2} \right )  (t \vee t_0)^{1 - \gamma}
 \quad . \eqno(6.20)$$

\noi {\bf Proof. Part(1).} Let $s_0(t)$ be defined by (5.36) supplemented by Proposition 5.5,
part (1) and define $s_{02}(t)$ by (6.17). The result now follows from Proposition 5.5, part (1)
and Lemma 6.1, part (1). In particular (6.18) follows from (5.43) and (6.5). \\

\noi {\bf Part (2).} For $t \geq t_0$, we estimate
$$|s_{02,t_0}(t) - s_{02}(t) |_{\ell - 1} \leq |s_{02,t_0}(t) - s(t)|_{\ell - 1} + |s_{02}(t) -
s(t) |_{\ell - 1}$$
\noi and we estimate the first norm in the RHS by (6.13) and the second norm by (6.18). This
yields (6.20) for $t \geq t_0$.  \par

For $t \leq t_0$, we obtain from (6.11) and (6.1)
$$s_{02,t_0}(t) - s_{02}(t) = s(t_0) - s_{02}(t_0) + \int_t^{t_0} dt' \ t'^{- \gamma} \left (
\nabla g_0(w_+) - \nabla g_0 (U^*(1/t_0) w(t_0) \right ) \quad . \eqno(6.21)$$
\noi We estimate the first difference in the RHS by (6.18) with $t = t_0$, and the integral by the
same method as in the proof of Lemma 6.1, part (1), and by the use of (5.33) and (6.9) with $t' =
\infty$, so that 
$$\Big | \int_t^{t_0} dt ' \cdots \Big |_{\ell - 1} \leq C (1 - \gamma )^{-1} \left ( a^2 \ b \
t_0^{1 - 2 \gamma} + a^2 \ t_0^{1/2 - \gamma} \right )$$
\noi which completes the proof of (6.20) for $t \leq t_0$. \par

The convergence stated in Part (2) follows from the estimate (6.20). \par \nobreak
\hfill $\sq$ \\

\noi {\bf Remark 6.1.} As announced in Section 2, the approximation of $s$ by $s_{02}$ for
solutions $(w, s)$ of the system (5.1) (5.2) is not as good as the approximation by $s_0$ obtained
in Section 5. This can be seen on (6.18) where the term $a^2 \ t^{1/2 - \gamma}$ is dominant for
large $t$ as compared with the term $(b^2 + a^2b) t^{1 - 2 \gamma}$ obtained from (5.43). \\

We now turn to the converse construction, namely to the construction of a solution $(w, s)$ of
the system (5.1) (5.2) which is asymptotic to a fixed solution $(w_+, s_{02})$ of (6.1), defined
by (6.2). The first step consists in constructing a solution $(w, s)$ which coincides with $(w_+,
s_0)$ at some time $t_0$. The next result is the analogue of Proposition 5.6. \\

\noi {\bf Proposition 6.3.} {\it Let $1/2 < \gamma < 1$ and let $(k, \ell )$ be an admissible
pair. Let $(w_+, s_{02}(1)) \in H^{k+1} \oplus X^{\ell + 1}$, define $s_{02}$ by (6.2) and let} 
$$a = |w_+|_{k+1} \qquad , \quad b = \parallel \widetilde{s}_{02} ; L^{\infty} ([1 , \infty ),
X^{\ell + 1} ) \parallel \quad . \eqno(6.22)$$
\noi {\it Then, there exist $T_0$ and $T$, $1 \leq T_0$, $T < \infty$, depending only on $(\gamma , a
, b)$, such that for all $t_0 \geq T_0 \vee T$, the system (5.1) (5.2) with initial data $w(t_0)
= U(1/t_0) w_+$, $s(t_0) = s_{02}(t_0)$, has a unique solution $(w_{t_0}, s_{t_0})$ such that
$(w_{t_0}, \widetilde{s}_{t_0}) \in ({\cal C} \cap L^{\infty} ) ([T, \infty ), H^k \oplus
X^{\ell})$. One can take}
$$T_0 = C \left \{ (b + (1 - \gamma )^{-1} a^2 )^{1/\gamma} \vee (1 - \gamma )^{-2} \ a^4 \ b^{-2}
\right \} \eqno(6.23)$$
$$T = C (2\gamma - 1)^{-2} \left ( (b + a^2)^{1/\gamma} \vee a^4 \ b^{-2} \right ) \eqno(6.24)$$
\noi and the solution satisfies the estimates
$$\left \{ \begin{array}{l} |w_{t_0} (t) - w_+ |_k \leq C \left ( a\ b \
t_0^{-\gamma} + a \ t_0^{-1/2} \right ) \\ \\ |s_{t_0}(t) - s_{02}(t)|_{\ell} \leq
C \left ( (b^2 + (1 - \gamma )^{-1} a^2 \ b ) t_0^{-\gamma} + (1 - \gamma
)^{-1} a^2 \ t_0^{-1/2} \right ) t^{1 - \gamma} \end{array}\right . \eqno(6.25)$$
\noi {\it for $t \geq t_0$},
$$\left \{ \begin{array}{l} |w_{t_0} (t) - w_+ |_k \leq C \left ( a\ b \
t^{-\gamma} + a \ t_0^{-1/2} \right ) \\ \\ |s_{t_0}(t) - s_{02}(t)|_{\ell} \leq
C (2 \gamma - 1)^{-1} \left ( (b^2 + a^2 \ b ) t^{1-2\gamma} + a^2 \ t^{1/2- \gamma} \right ) 
 \end{array}\right . \eqno(6.26)$$ \noi {\it for $T \leq t \leq t_0$, and}
$$|w_{t_0}(t)|_k \leq C a \qquad , \quad  |s_{t_0}(t)|_{\ell} \leq C \ b \ t^{1 - \gamma} \eqno(6.27)$$
\noi {\it for all $t \geq T$.} \\

\noi {\bf Proof.} Let $w_0 = U(1/t_0)w_+$ and define $s_0 (t)$ by (6.3) so that $s_0(t)$ solves (5.3)
and satisfies $s_0(t_0) = s_{02}(t_0)$. By Lemma 6.1, part (2), $s_0 \in {\cal C}([1, \infty ),
X^{\ell + 1})$ and $s_0$ satisfies
$$\parallel \widetilde{s}_0 ; L^{\infty}([T, \infty ), X^{\ell + 1} ) \parallel \ \leq C \ b
\eqno(6.28)$$
\noi provided $(1 - \gamma )t_0^{1/2} \geq Ca^2b^{-1}$ and $(2 \gamma - 1)T^{1/2} \geq Ca^2b^{-1}$
which follow from (6.23) (6.24) for $t_0 \geq T_0$. We now apply Proposition 5.6 with $s_0$ just
defined. Let $(w_{t_0}, s_{t_0})$ be the solution of the system (5.1) (5.2) thereby obtained under
the conditions (5.51) (5.52) which also follow from (6.23) (6.24). That solution satisfies the
required initial condition $w_{t_0}(t_0) = w_0$ and $s_{t_0} (t_0) = s_0(t_0) = s_{02}(t_0)$.
Furthermore it satisfies the estimates (5.53) (5.54) (5.55) with however $w_+$ replaced by $w_0$. The
estimates (6.25) (6.26) (6.27) follow from the previous ones, from Lemma 6.1, part (2) and from the
estimate
$$|w_0 - w_+|_k = |(U(1/t_0) - 1)w_+|_k \leq t_0^{-1/2} \ |w_+|_{k+1} \quad . \eqno(6.29)$$
\par \nobreak
\hfill $\sq$ \\

We now take the limit $t_0 \to \infty$ of the solution $(w_{t_0}, s_{t_0})$ constructed in
Proposition 6.3 for fixed $(w_+, s_{02})$. The next result is the analogue of Proposition 5.7. \\

\noi {\bf Proposition 6.4.} {\it Let $1/2 < \gamma < 1$ and let $(k, \ell )$ be an admissible pair.}
\par

{\it (1) Let $(w_+ , s_{02}(1)) \in H^{k+1} \oplus X^{\ell + 1}$ and define $s_{02}$, $a$, $b$ by
(6.2) (6.22). Then there exists $T$, $1 \leq T < \infty$, depending only on $(\gamma , a, b)$ and a
unique solution $(w, s)$ of the system (5.1) (5.2) such that $(w, \widetilde{s}) \in ({\cal C} \cap
L^{\infty})([T, \infty ), H^k \oplus X^{\ell})$ and such that the following estimates hold for all
$t \geq T$~:}
$$\left \{ \begin{array}{l} |w(t) - w_+|_k \leq C \ ab \ t^{- \gamma} \quad , \\ \\
|s(t) - s_{02}(t)|_{\ell} \leq C (2 \gamma - 1)^{-1} \left ( (b^2 + a^2b ) t^{1 - 2
\gamma} + a^2 \ t^{1/2 - \gamma} \right ) \quad , \end{array}\right . \eqno(6.30)$$

$$|w(t)|_k \leq C \ a \quad , \quad |s(t)|_{\ell} \leq C \ b \ t^{1 - \gamma} \quad . \eqno(6.31)$$
\noi {\it One can take}
$$T = C (2 \gamma - 1)^{-2} \left ( (b + a^2)^{1/\gamma} \vee a^4 \ b^{-2} \right ) \quad .
\eqno(6.32)$$

{\it (2) Let $(w_{t_0}, s_{t_0})$ be the solution of the system (5.1) (5.2) constructed in
Proposition 6.3 for $t_0 \geq T_0 \vee T$ and in particular such that $(w_{t_0} ,
\widetilde{s}_{t_0}) \in ({\cal C} \cap L^{\infty}) ([T, \infty ), H^k \oplus X^{\ell})$. Then
$(w_{t_0}, s_{t_0})$ converges to $(w, s)$ in norm in $L^{\infty}(J, H^{k-1} \oplus X^{\ell - 1})$
and in the weak-$*$ sense in $L^{\infty}(J, H^k \oplus X^{\ell})$ for any compact interval $J
\subset [T, \infty )$, and in the weak-$*$ sense in $H^k \oplus X^{\ell}$ pointwise in $t$.} \par

{\it (3) The map $(w_+, s_{02}(1)) \mapsto (w, s)$ defined in Part (1) is continuous on the bounded
sets of $H^{k+1} \oplus X^{\ell + 1}$ from the norm topology of $(w_+, s_{02}(1))$ in $H^{k-1}
\oplus X^{\ell - 1}$ to the norm topology of $(w, s)$ in $L^{\infty}(J, H^{k-1} \oplus X^{\ell -
1})$ and to the weak-$*$ topology in $L^{\infty}(J, H^k \oplus X^{\ell})$ for any compact interval
$J \subset [T, \infty )$, and to the weak-$*$ topology in $H^k \oplus X^{\ell}$ pointwise in $t$.} \\

\noi {\bf Proof. Part (1).} Take now (6.17) as the definition of $s_0$ in terms of $s_{02}$, so that
$s_0$ satisfies (6.28) with $T$ given by (6.32) by Lemma 6.1, part (2). Let $(w, s)$ be the solution
of the system (5.1) (5.2) constructed in Proposition 5.7. Then $(w, s)$ satisfies the properties
stated in Part (1). In particular (6.30) follows from (5.69) and from (6.7). Uniqueness of $(w, s)$
follows from Proposition 5.2. \\

\noi {\bf Part (2).} We estimate $y = |w_{t_0}(t) - w(t)|_{k-1}$ and $z = |s_{t_0} (t) - s(t)|_{\ell
- 1}$ in the same way as in the proof of Proposition 5.7, part (2) by using Lemma 3.9 and Lemma 5.1,
with initial conditions $y(t_0)$ and $z(t_0)$ estimated by (6.26) and (6.30) at time $t_0$. The rest
of the proof is identical with that of Proposition 5.7, part (2). \\

\noi {\bf Part (3).} The proof is almost identical with that of Proposition 5.7, part (3) with $s_0$
replaced everywhere by $s_{02}$. The only difference is the appearance of an additional term $(2
\gamma - 1)^{-1}a^2t^{1/2 - \gamma}$ in the RHS of the second inequality in (5.79), coming from
(6.30). This leads to the choice $t_0^{-1/2} = \varepsilon$ instead of $t_0^{-\gamma} = \varepsilon
b^{-1}$ so that in (5.81) (5.82) the factor $\varepsilon^{(2 \gamma - 1)/\gamma}$ is replaced by
$\varepsilon^{2\gamma - 1}$. \par \nobreak
\hfill $\sq$ \\

Proposition 5.8 applies mutatis mutandis to the map $(w_+, s_{02}(1)) \to (w, s)$ because the map
$(w_+, s_0(1)) \to (w_+, s_{02} (1))$ is bijective since it is defined by the explicit formula
(6.17) and exactly preserves the relevant regularity by Lemma 6.1. 

\section{Wave operators and asymptotics for $u$.}
\hspace*{\parindent}
In this section we complete the construction of the wave operators for the equation (1.1) and we
derive asymptotic properties of solutions in their range. The construction relies in an essential
way on those of Sections 5 and 6, esp. Proposition 6.4, and will require a discussion of the gauge
invariance of those constructions. The first task is to supplement them with the determination of
the phase $\varphi$, which appears so far only through its gradient $s$ (see (2.26) (2.31)).
Actually the treatment in Sections 4-6 only involved the variable $s$, and dit not even assume
that $s$ was a gradient. \par

We recall that $g_0$, $g$ are defined by (3.1) (3.2) and we continue to use the short hand
notation $g_0(w, w) = g_0(w)$ so that
$$g(w, w) = g_0 (U^*(1/t) w) \quad .$$
\noi We are interested in solving the system (2.25) (2.26) which we rewrite as
$$\partial_t w = (2t^2)^{-1} U(1/t) (2 \nabla \varphi \cdot \nabla + (\Delta \varphi )) U^*(1/t) w
\eqno(2.25) \equiv (7.1)$$
$$\partial_t \varphi = (2t^2)^{-1} |\nabla \varphi |^2 + t^{-\gamma} \ g_0(U^*(1/t) w) \quad .
\eqno(2.26) \sim (7.2)$$
\noi In analogy with (5.5), we shall use the notation $\widetilde{\varphi} = t^{\gamma - 1}\varphi$.
The relevant spaces for the phase $\varphi$ are the spaces $Y^{\ell}$ defined by
$$Y^{\ell} = \left \{ \varphi : \varphi \in L^{\infty} \quad \hbox{and} \quad \nabla \varphi \in
X^{\ell} \right \} \quad . \eqno(7.3)$$

We first consider the Cauchy problem with finite initial time $t_0$ for the system (7.1) (7.2) with
initial data $(w_0, \varphi_0) \in H^k \oplus Y^{\ell}$. That problem is solved by first solving
the Cauchy problem for the system (5.1) (5.2) with initial data $(w_0, \nabla \varphi_0)$ at
time $t_0$ for $(w, s)$ by Proposition 4.1 or 5.1 and then recovering $\varphi$ from $(w, s)$ by
$\varphi (t_0) = \varphi_0$ and
$$\varphi (t) = \varphi (t_0) + \int_{t_0}^t dt' \left \{ (2t'^2)^{-1} |s|^2 + t'^{-\gamma} \
g_0 (U^*(1/t')w) \right \} \quad . \eqno(7.4)$$
\noi As mentioned in Section 2, it follows from (5.2) that the vorticity $\omega = \nabla
\times s$ satisfies the equation (2.35). Furthermore $\omega (t_0) = 0$. Under the available
regularity properties of $s$, it follows from Lemma 3.2 with $p = \infty$, $u = \omega$, $v =
s$ and from Gronwall's Lemma that $\omega (t) = 0$ for all $t$. On the other hand, from (5.2)
and (7.2), it follows that
$$(s - \nabla \varphi )(t) = \int_{t_0}^t dt' \ t'^{-2} (s \times \omega ) (t') \eqno(7.5)$$
\noi where $(s \times \omega)_i = \sum\limits_j s_j \omega_{ji}$, and therefore $s(t) = \nabla
\varphi (t)$ for all $t$ for which $(w, s)$ is defined since $\omega = 0$. \par

If $(w, \varphi )$ is a solution of the system (7.1) (7.2) as obtained from Proposition 5.1
and from the previous argument, it follows from that proposition and from Corollary 3.1 that
$(w, \widetilde{\varphi}) \in ({\cal C} \cap L^{\infty})([T, \infty ), H^k \oplus Y^{\ell})$.
\par

We next consider the Cauchy problem with infinite initial time covered by Propositions 5.5 and
5.7. There we have established a correspondence between solutions $(w, s)$ of (5.1) (5.2) and
solutions $(w_+, s_0)$ of (5.3), and we now extend it to a correspondence between solutions
$(w, \varphi )$ of (7.1) (7.2) and solutions $(w_+, \varphi_0)$ of the auxiliary free equation
$$\partial_t \varphi_0 = t^{-\gamma} \ g_0(U^*(1/t) w_+) \quad . \eqno(2.28) \sim (7.6)$$
\noi The general solution of (7.6) can be written as
$$\varphi_0(t) = \varphi_0 (1) + \int_1^t dt' \ t'^{-\gamma} \ g_0(U^*(1/t') w_+) \eqno(7.7)$$
\noi and if $s_0(t)$ and $\varphi_0(t)$ are defined by (5.4) (7.7) with $s_0(t) = \nabla
\varphi_0(t)$ for one $t$ (for instance $t = 1$), the same relation holds for all $t$. We
supplement the correspondence established in Propositions 5.5 and 5.7 with the relation
$$\varphi (t) = \varphi_0 (t) - \int_t^{\infty} dt' \left \{ (2t'^2)^{-1} |s|^2 + t'^{-\gamma}
\left ( g_0 (U^*(1/t') w) - g_0 (U^*(1/t') w_+) \right ) \right \} \eqno(7.8)$$
\noi which will be used to define $\varphi (t)$ and/or $\varphi_0 (t)$ in terms of each other.
>From (5.36) and (7.8) it follows that
$$(s - \nabla \varphi )(t) - (s_0 - \nabla \varphi_0)(t) = - \int_t^{\infty} dt' \ t'^{-2} (s
\times \omega ) (t') \quad . \eqno(7.9)$$

Corresponding to the situation of Proposition 5.5, let $(w, \varphi ) \in ({\cal C} \cap
L_{loc}^{\infty}) ([T, \infty ), H^k \oplus Y^{\ell})$ be a solution of (7.1) (7.2) such that
$(w, s = \nabla \varphi )$ satisfies the assumptions of Proposition 5.5. Define $(w_+, s_0)$
by Proposition 5.3, by (5.36) and by Proposition 5.5, part (1), and define $\varphi_0$ by
(7.8). Then it follows from (7.9) that $s_0(t) = \nabla \varphi_0 (t)$ for all $t$. \par

Conversely, corresponding to the situation of Proposition 5.7, let $(w_+, \varphi_0(1)) \in
H^{k+1} \oplus Y^{\ell + 1}$, define $\varphi_0(t)$ by (7.7) and $s_0(t)$ by $s_0(t) = \nabla
\varphi_0(t)$, define $(w, s)$ by Proposition 5.7, part (1) and define $\varphi (t)$ by
(7.8). Let $(w_{t_0}, s_{t_0})$ be defined by Proposition 5.6. From the finite initial time
results it follows that $\omega_{t_0} (t) \equiv \nabla \times s_{t_0}(t) = 0$ for all $t$
and $t_0$, while by Proposition 5.7, part (2) $\omega_{t_0}$ converges to $\omega$ in
$X^{\ell - 2}$ uniformly in compact intervals, so that $\omega = 0$. It then follows again
from (7.9) that $s(t) = \nabla \varphi (t)$ for all $t$. This proves that the correspondence
established in Propositions 5.5 and 5.7 extends to a correspondence between solutions $(w,
\varphi )$ of (7.1) (7.2) and solutions $(w_+ , \varphi_0)$ of (7.6) preserving the relations
$s = \nabla \varphi$, $s_0 = \nabla \varphi_0$. \par

The same discussion applies mutatis mutandis to the situation of Propositions 6.2 and 6.4 and
allows for an extension of the correspondence between solutions $(w, s)$ of (5.1) (5.2) and
solutions $(w_+, s_{02})$ of (6.1) established there to a correspondence between solutions
$(w, \varphi )$ of (7.1) (7.2) and solutions $(w_+, \varphi_{02}$) of the equation 
$$\partial_t \ \varphi_{02} = t^{- \gamma} \ g_0(w_+) \eqno(2.29) \sim (7.10)$$
\noi which preserves the relations $s = \nabla \varphi$, $s_{02} = \nabla \varphi_{02}$. The
relation (7.8) between $\varphi (t)$ and $\varphi_0(t)$ is replaced by the relation
$$\varphi (t) = \varphi_{02} (t) - \int_t^{\infty} dt' \left \{ (2t'^2)^{-1} |s|^2 + t'^{-
\gamma} \left (  g_0(U^*(1/t') w) - g_0(w_+) \right ) \right \} \eqno(7.11)$$
\noi between $\varphi (t)$ and $\varphi_{02}(t)$, so that $\varphi_0 (t)$ and $\varphi_{02}
(t)$ are related by 
$$\varphi_0(t) = \varphi_{02} (t) - \int_t^{\infty} dt' \ t'^{- \gamma}
\left ( g_0 (U^*(1/t') w_+) - g_0(w_+) \right ) \eqno(7.12)$$
\noi in agreement with (6.16) and (6.17). \par

The equation (7.10) is trivially solved by
$$\varphi_{02} (t) = \varphi_{02} (1) + (1 - \gamma)^{-1} \ (t^{1 - \gamma} - 1) \ g_0(w_+) \quad
. \eqno(7.13)$$

We can now embark on the explicit construction of wave operators, and we first construct a
wave operator $W$ for the auxiliary system (5.1) (5.2). This could be based either on
Proposition 5.7 or on Proposition 6.4. We choose the latter since as mentioned before it is
better suited than the former for the discussion of gauge invariance to be performed next. \\

{\bf Definition 7.1.} We define $W$ as the map
$$W : (w_+, \varphi_{02} (1)) \to (w, \varphi ) \eqno(7.14)$$
\noi from $H^{k+1} \oplus Y^{\ell + 1}$ to the set of $(w, \varphi )$ such that $(w,
\widetilde{\varphi } ) \in ({\cal C} \cap L^{\infty}) ([T, \infty ), H^k \oplus Y^{\ell})$
for some $T$, $1 \leq T < \infty$, as follows. Define $\varphi_{02} (t)$ by (7.13) and
$s_{02} (t)$ by $s_{02} (t) = \nabla \varphi_{02} (t)$. Define $(w, s)$ by Proposition 6.4,
part (1) and finally define $\varphi$ by (7.11), so that $s = \nabla \varphi$ by the
previous discussion. Then $W$ is well defined by (7.14) as a map between the spaces
indicated. \par

Before defining the wave operators for $u$, we now study the gauge invariance of $W$, which
plays an important role in justifying that definition, as was explained in Section 2. For
that purpose we need some information on the Cauchy problem for the equation (1.1) at
finite times. In addition to the operators $M = M(t)$ and $D = D(t)$ defined by (2.5)
(2.6), we introduce the operator
$$J = J(t) = x + it \nabla \quad , \eqno(7.15)$$
\noi the generator of Galilei transformations. The operators $M$, $D$, $J$ satisfy the
commutation relation
$$i \ M \ D \ \nabla = J \ M \ D \quad . \eqno(7.16)$$
\noi For any interval $I \subset [1, \infty )$ and any nonnegative integer, we define the
space
$$\begin{array}{ll} {\cal X}^k(I) &= \left \{ u : D^* M^* u \in
{\cal C} (I, H^k) \right \}  \\
& \\
&= \left \{ u : <J(t)>^k u \in {\cal C} (I, L^2) \right \} 
\end{array} \eqno(7.17)$$ 
\noi where $<\lambda > = (1 + \lambda^2)^{1/2}$ for any real number or self-adjoint
operator $\lambda$ and where the second equality follows from (7.16). Then  \\

\noi {\bf Proposition 7.1.} {\it Let $k$ be a positive integer and let $0 < \mu < n \wedge 2k$.
Then the Cauchy problem for the equation (1.1) with initial data $u(t_0) = u_0$ such that
$<J(t_0)>^k$ $u_0 \in L^2$ at some initial time $t_0 \geq 1$ is locally well posed in
${\cal X}^k(\cdot )$, namely} \par
{\it (1) There exists $T > 0$ such that (1.1) has a unique solution with initial data $u(t_0) = u_0$
in ${\cal X}^k ([1 \vee (t_0 - T), t_0 + T])$.} \par
{\it (2) For any interval $I$, $t_0 \in I \subset [1 , \infty )$, (1.1) with initial data $u(t_0) =
u_0$ has at most one solution in ${\cal X}^k(I)$.} \par
{\it (3) The solution of Part (1) depends continuously on $u_0$ in the norms considered there.} \\

\noi {\bf Proof.}  The proof is obtained by minor variations of the corresponding results in
\cite{8r}. The differences come from the factor $t^{\mu - \gamma}$ in (1.2), which is irrelevant for
the present problem, and from the replacement of ${\cal C}(I, H^k)$ by ${\cal X}^k(I)$. That
replacement is made possible by the properties of the operator $J(t)$ and especially the commutation
relation (7.16), which implies that $J(t)$ behaves as a derivative on gauge invariant functions. See
for instance \cite{2r}. \par \nobreak
\hfill $\sq$ \\

In the study of gauge invariance for $W$ we shall actually need only the uniqueness statement, Part
(2) of Proposition 7.1.  \par

We recall that in the transition from the system (5.1) (5.2) to the equation (1.1), $u$ should be
defined by (2.14) with $(w, \varphi ) = (w_3, \varphi_3)$ and accordingly we define the map
$$\Phi : (w, \varphi ) \to u = M \ D \exp [-i \varphi ] \ U^*(1/t) w \quad . \eqno(7.18)$$
\noi That map satisfies the following property. \\

\noi {\bf Lemma 7.1.} {\it The map $\Phi$ defined by (7.18) is bounded and continuous from ${\cal C}
(I, H^k \oplus Y^{\ell})$ to ${\cal X}^k(I)$ for any admissible pair $(k, \ell )$ and any interval $I
\subset [1, \infty )$.} \\

\noi {\bf Proof.} An immediate consequence of Lemma 3.5. \par \nobreak
\hfill $\sq$ \\

We can now make the following definition. \\

\noi {\bf Definition 7.2.} Let $(k, \ell )$ be an admissible pair and let $(w, \varphi )$ and
$(w', \varphi ')$ be two solutions of the system (7.1) (7.2) in ${\cal C} (I, H^k \oplus
Y^{\ell})$ for some interval $I \subset [1, \infty)$. We say that $(w, \varphi )$ and $(w',
\varphi ')$ are gauge equivalent if they give rise to the same $u$, namely $\Phi ((w, \varphi ))
= \Phi ((w', \varphi '))$, or equivalently if
$$\exp [ - i \varphi (t) ] \ U^*(1/t) \ w(t) = \exp [- i \varphi '(t)] \ U^*(1/t) \ w'(t)
\eqno(7.19)$$ \noi for all $t \in I$. \par

A sufficient condition for gauge equivalence is given by the following Lemma. \\   

\noi {\bf Lemma 7.2.} {\it Let $(k, \ell )$ be an admissible pair and let $(w, \varphi )$
and $(w', \varphi ')$ be two solutions of the system (7.1) (7.2) in ${\cal C} (I, H^k
\oplus Y^{\ell})$. In order that $(w, \varphi )$ and $(w', \varphi ')$ be gauge equivalent,
it is sufficient that (7.19) holds for one $t \in I$.} \\

\noi {\bf Proof.} An immediate consequence of Lemma 7.1, of Proposition 7.1, part (2), and
of the fact that $(k, \ell)$ admissible implies $k > 1 + \mu /2$. \par \nobreak
\hfill $\sq$ \\

The gauge covariance properties of $W$ will be expressed by the following two propositions.
\\

\noi {\bf Proposition 7.2.} {\it Let $0 < \gamma < 1$ and let $(k, \ell )$ be an admissible
pair. Let $(w, \varphi )$ and $(w', \varphi ')$ be two solutions of the system (7.1) (7.2)
such that $(w, \widetilde{\varphi}), (w', \widetilde{\varphi}') \in ({\cal C} \cap
L^{\infty})([T, \infty ), H^k \oplus Y^{\ell})$ for some $T \geq 1$, and assume that $(w,
\varphi )$ and $(w', \varphi ')$ are gauge equivalent. Then} \par

{\it (1) There exists $\sigma \in Y^{\ell - 1}$ such that $\varphi ' (t) - \varphi (t)$
converges to $\sigma$ strongly in $Y^{\ell - 2}$ and in the weak-$*$ sense in $Y^{\ell -
1}$. The following estimates holds~:}
$$\parallel \varphi ' (t) - \varphi (t) - \sigma ; Y^{\ell - 2} \parallel \ \leq A \ t^{-
\gamma} \eqno(7.20)$$
\noi {\it for some constant $A$ depending on $T$ and on the norms of
$\widetilde{\varphi}$, $\widetilde{\varphi}'$ in $L^{\infty} (\cdot , Y^{\ell})$, with
the exception of the case $n$ even, $\ell = n/2 + 1$ where the $L^{\infty}$ norm of
$\nabla \varphi$ satisfies only}
$$\parallel \nabla \varphi '(t) - \nabla \varphi (t) - \nabla \sigma \parallel_{\infty}
\ \leq A \ t^{- \gamma /2} \quad . \eqno(7.21)$$

{\it (2) Assume in addition that $\gamma > 1/2$. Then $\varphi '(t) - \varphi (t)$
converges to $\sigma$ in norm in $Y^{\ell - 1}$ and the following estimate holds~:}
$$\parallel \varphi '(t) - \varphi (t) - \sigma ; Y^{\ell - 1} \parallel \ \leq A \
t^{1 - 2 \gamma} \quad . \eqno(7.22)$$
\noi {\it Furthermore $w'_+ = w_+ \exp [i \sigma ]$ where $w_+$, $w'_+$ are the limits of
$w(t)$, $w'(t)$ as $t \to \infty$ obtained in Proposition 5.3.} \par

{\it (3) Assume in addition that $\gamma > 1/2$ and that $(w, \varphi )$, $(w',
\varphi ') \in {\cal R} (W)$. Then $\varphi '_{02}(t) = \varphi_{02}(t) + \sigma$ for
all $t \geq 1$. In particular $\sigma \in Y^{\ell +1}$. } \\

\noi {\bf Proof. Part (1).} Define $\varphi_- = \varphi ' - \varphi$, $s = \nabla
\varphi$, $s' = \nabla \varphi '$, $s_{\pm} = s' \pm s$, and 
$$b = \parallel \widetilde{s};L^{\infty}([T, \infty ), X^{\ell}) \parallel \ \vee \
\parallel \widetilde{s}' ; L^{\infty} ([T, \infty ), X^{\ell}) \parallel \quad .
\eqno(7.23)$$
\noi From (5.2) and gauge equivalence it follows that
$$\partial_t s_- = t^{-2} \left ( (s_- \cdot \nabla ) s_+ + (s_+ \cdot \nabla ) s_-
\right ) \eqno(7.24)$$
\noi and therefore by Lemma 3.9
$$\partial_t |s_-|_{\ell - 1} \leq C \ t^{-2} |s_-|_{\ell - 1} \ |s_+|_{\ell} \leq C
\ b \ t^{-1-\gamma} \ |s_-|_{\ell - 1}$$
\noi so that by Gronwall's Lemma, for all $t \geq T$,
$$|s_-(t)|_{\ell - 1} \leq |s_-(T)|_{\ell - 1} \exp ( C \ b\ \gamma^{-1} \
T^{-\gamma} ) \equiv A_0 \quad ,$$
\noi namely
$$\parallel s_-; L^{\infty} ([T, \infty ), X^{\ell - 1})\parallel \ \leq A_0 \quad .
\eqno(7.25)$$

>From (7.2) and gauge equivalence, it follows that
$$\partial_t \ \varphi_- = (2t^2)^{-1} \ (s_- \cdot s_+) $$
\noi and therefore by (7.23) (7.25) for any $t \geq t_0 \geq T$
$$\parallel \varphi_-(t) - \varphi_-(t_0) \parallel_{\infty} \ \leq \int_{t_0}^t dt'
\ b \ A_0 \ t'^{-1 - \gamma} \leq b \ A_0 \ \gamma^{-1} \ t_0^{- \gamma} \quad .
\eqno(7.26)$$
\noi This implies that $\varphi_-(t)$ converges in norm in $L^{\infty}$ to some
$\sigma \in L^{\infty}$ and that
$$\parallel \varphi_- (t) - \sigma \parallel_{\infty} \ \leq b \ A_0 \ \gamma^{-1} \
t^{- \gamma} \eqno(7.27)$$
\noi which is the part of (7.20) involving $\varphi_-$ itself (and not $s_-$ only). From the
uniform estimate (7.25) and standard compactness arguments, it follows that $\nabla \sigma \in
X^{\ell - 1}$, that $|\nabla \sigma |_{\ell - 1} \leq A_0$ and that $s_-$ converges to $\nabla
\sigma$ in the weak-$*$ sense in $X^{\ell - 1}$, which together with (7.27) implies weak-$*$
convergence of $\varphi_-$ to $\sigma$ in $Y^{\ell - 1}$.  \par

We finally prove the strong convergence of $s_-$ to $\nabla \sigma$ in $X^{\ell - 2}$. From (7.24)
which we rewrite as
$$\partial_t (s_- - \nabla \sigma ) = t^{-2} \left \{ \left ( (s_- - \nabla \sigma ) \cdot \nabla
\right ) s_+ + (s_+ \cdot \nabla ) (s_- - \nabla \sigma ) + (\nabla \sigma \cdot \nabla ) s_+ + (s_+
\cdot \nabla ) \nabla \sigma \right \}\eqno(7.28)$$

\noi and by the same estimates as in the proof of Lemma 3.11, we obtain
$$\begin{array}{ll} \partial_t |s_- - \nabla \sigma |_{\ell - 2} &\leq C \ t^{-2} \left \{ |s_- -
\nabla \sigma |_{\ell - 2} \ |s_+|_{\ell} \ + |\nabla \sigma |_{\ell - 1} \ |s_+|_{\ell} \right \}
\\ &\\ &\leq C \ b \ t^{-1 - \gamma} \left ( |s_- - \nabla \sigma |_{\ell - 2} + A_0 \right ) 
\end{array}\eqno(7.29)$$
\noi with the only exception of the case $n$ even, $\ell = n/2 + 1$, where the $L^{\infty}$ norm of
$(s_- - \nabla \sigma )$ which occurs in the $X^{\ell - 2}$ norm, is not estimated as in (7.29)
because $\parallel \partial^2 \sigma \parallel_{\infty}$ is not controlled by $|\nabla \sigma
|_{\ell - 1}$. From (7.29) we obtain by Gronwall's Lemma
$$|s_-(t) - \nabla \sigma |_{\ell - 2} \leq A_0 \left \{ \exp (C \ b \ \gamma^{-1} \ t^{- \gamma} )
- 1 \right \} \eqno(7.30)$$
\noi which together with (7.27) completes the proof of (7.20). \par

For $n$ even, $\ell = n/2 + 1$, we estimate simply
$$\parallel s_-(t) - \nabla \sigma \parallel_{\infty} \ \leq C \parallel \varphi_-(t) - \sigma
\parallel_{\infty}^{1/2} \ \parallel s_-(t) - \nabla \sigma ; \dot{H}^{n/2+1} \parallel^{1/2}$$ \noi
and (7.21) follows from (7.27) and from the $\dot{H}^{n/2}$ part of (7.30). \\

\noi {\bf Part (2).} Let $t \geq t_0 \geq T$, $s_- (t) = s_-$ and $s_- (t_0) = s_0$. We rewrite
(7.24) as 
$$\partial_t(s_- - s_0) = t^{-2} \left \{ ((s_- - s_0)\cdot \nabla )s_+ + (s_+ \cdot \nabla ) (s_-
- s_0) + (s_0 \cdot \nabla ) s_+ + (s_+ \cdot \nabla ) s_0 \right \} \quad . $$ \noi By Lemma 3.9, we
estimate
$$\begin{array}{ll} \partial_t |s_- - s_0|_{\ell - 1} &\leq C \ t^{-2} \left \{ |s_- - s_0|_{\ell - 1}
\ |s_+|_{\ell} + |s_-|_{\ell - 1} \ |s_+|_{\ell} + |s_+|_{\ell - 1} \ |s_0|_{\ell} \right \} \\ & \\
&\leq C \ t^{-2} \left ( |s_- - s_0|_{\ell - 1} + |s_0|_{\ell} \right ) |s_+|_{\ell} \\ & \\
&\leq C \ b \ t^{-1-\gamma} \left ( |s_- - s_0|_{\ell - 1} + b \ t_0^{1 - \gamma} \right )
\end{array}$$ \noi and therefore by Gronwall's Lemma
$$|s_-(t) - s_-(t_0)|_{\ell - 1} \leq b \ t_0^{1 - \gamma} \left ( \exp (C \ b \ \gamma^{-1} \
t_0^{-\gamma} ) - 1 \right ) \quad .$$
\noi This yields a separate proof of the convergence of $s_-(t)$ in $X^{\ell - 1}$, together with
the estimate
$$|s_-(t) - \nabla \sigma |_{\ell - 1} \leq b \ t^{1 - \gamma} \ \left ( \exp (C\ b\ \gamma^{-1} \
t^{-\gamma} )- 1 \right ) \quad ,$$
\noi which together with (7.27) completes the proof of (7.22). \par

We now prove that $w'_+ = w_+ \exp (i \sigma )$. For that purpose we estimate
$$|w'_+ - w_+ e^{i \sigma}|_{k-1} \leq |w'_+ - U^*(1/t) w'(t)|_{k-1} + \Big | U^*(1/t) w'(t) - \exp
[i \varphi_- (t) ] \ U^* (1/t) w(t) \Big |_{k-1}$$ $$+ \Big | \left ( \exp [i (\varphi_- (t) - \sigma
) ] - 1 \right ) \exp (i \sigma ) \ U^*(1/t) w(t) \Big |_{k-1} + \Big | \exp (i \sigma ) \left
(U^*(1/t) w(t) - w_+ \right ) \Big |_{k-1} \quad . \eqno(7.31)$$ 
\noi We estimate the first norm in the RHS by Proposition 5.3, esp. (5.33) and by (6.9) as 
$$\Big | w'_+ - U^*(1/t) w'(t) \Big |_{k-1} \leq A \left ( t^{-\gamma} + t^{-1/2} \right ) \quad .
\eqno(7.32)$$ 
\noi The second norm in the RHS of (7.31) is zero by gauge equivalence. The third norm
is estimated by Lemma 3.5 as
$$|\cdot | \leq C \left ( \parallel \varphi_- - \sigma \parallel_{\infty} + |\nabla (\varphi_- -
\sigma )|_{\ell - 1} \left ( 1 + |\nabla (\varphi_- - \sigma )|_{\ell - 1} \right )^{k-2} \right )
 \left ( 1 + |\nabla \sigma|_{\ell -1} \right )^{k-1} \ |w|_{k-1} \eqno(7.33)$$  
\noi and the last norm is estimated by Lemma 3.5 again followed by the analogue of (7.32) for $w$ as
$$|\cdot | \leq \left ( 1 + |\nabla \sigma|_{\ell - 1} \right )^{k-1} \ A\left ( t^{-\gamma} +
t^{-1/2} \right ) \quad . \eqno(7.34)$$

Collecting (7.32) (7.33) (7.34) and using (7.25) and (7.22) shows that the RHS of (7.31) tends to
zero when $t \to \infty$ and therefore that the LHS is zero since it is time independent. \\ 

\noi {\bf Part (3).} We recall that in the situation of Proposition 6.4 and of the definition of
$W$, $\varphi (t)$ and $\varphi_{02}(t)$ are related by (7.11), and by the estimates (3.35) (6.30)
(6.31)
$$\parallel \varphi (t) - \varphi_{02} (t)\parallel_{\infty} \ \leq A \ t^{1/2 - \gamma}
\eqno(7.35)$$ \noi for some $A$ depending only on $a$, $b$ defined by (6.22) and for $t$
sufficiently large. Similarly
$$\parallel \varphi '(t) - \varphi '_{02} (t)\parallel_{\infty} \ \leq A \ t^{1/2 - \gamma} \quad . \eqno(7.36)$$
\noi It follows then from (7.35) (7.36) and (7.20) (or (7.22)) that 
$$\parallel \varphi '_{02}(t) - \varphi_{02} (t) - \sigma \parallel_{\infty} \ \to 0 \quad
\hbox{when} \ t \to \infty \quad .$$
\noi On the other hand from (7.10) (or (7.13)) and from the condition $w'_+ = w_+ \exp (i \sigma )$
it follows that $\varphi '_{02}(t) - \varphi_{02} (t)$ is constant in time. Therefore $\varphi
'_{02} (t) - \varphi_{02}(t) = \sigma$ for all $t$. \par \nobreak
\hfill $\sq$ \\

Proposition 7.2 prompts us to make the following definition of gauge equivalence for asymptotic
states. \\

\noi {\bf Definition 7.3.} Two pairs $(w_+, \varphi_{02}(1))$ and $(w'_+ , \varphi '_{02}(1))$ are
gauge equivalent if there exists a real function $\sigma = \sigma (x)$ such that $w'_+ = w_+ \exp (i
\sigma )$ and $\varphi '_{02} (1) - \varphi_{02} (1) = \sigma$. \\

Two gauge equivalent pairs generate two solutions $(w_+, \varphi_{02})$ and $(w'_+, \varphi '_{02})$
of (7.10) such that $\varphi '_{02} (t) - \varphi_{02} (t) = \sigma$ for all $t \geq 1$. Those two
solutions will also be said to be gauge equivalent. \par

In Definition 7.3, we have not specified the regularity of $(w_+, \varphi_{02}(1))$ and $(w'_+,
\varphi '_{02}(1))$. This can be done easily, depending on the needs, and possibly with the help of
Lemma 3.5. \par

With the previous definition, Proposition 7.2 should be understood to mean that two gauge equivalent
solutions of the system (7.1) (7.2) in ${\cal R}(W)$ are images of two gauge equivalent solutions of
(7.10). The next proposition states that conversely two gauge equivalent solutions of (7.10) have
gauge equivalent images under $W$. \\

\noi {\bf Proposition 7.3.} {\it Let $1/2 < \gamma < 1$ and let $(k, \ell )$ be an admissible pair.
Let $(w_+, \varphi_{02}(1))$, $(w'_+ , \varphi '_{02} (1)) \in H^{k+1} \oplus Y^{\ell + 1}$ be
gauge equivalent, and let $(w, \varphi )$, $(w', \varphi ')$ be their images under $W$. Then $(w,
\varphi )$ and $(w' , \varphi ')$ are gauge equivalent, and $\lim\limits_{t \to \infty} \varphi
'(t) - \varphi (t) = \sigma \equiv \varphi '_{02} (1) - \varphi_{02}(1)$.} \\

\noi {\bf Proof.} Let $t_0$ be sufficiently large and let $(w_{t_0}, \varphi_{t_0})$ and
$(w'_{t_0}, \varphi '_{t_0})$ be the solutions of the system (7.1) (7.2) constructed by
Proposition 6.3 supplemented with (7.4) with the appropriate initial conditions at $t_0$, namely
$$\varphi_{t_0} (t_0) = \varphi_{02}(t_0) \quad , \quad \varphi '_{t_0} (t_0) = \varphi '_{02}(t_0)
\quad . \eqno(7.37)$$
\noi From (7.37) and from the initial conditions
$$w_{t_0} (t_0) = U(1/t_0) w_+ \quad , \quad w'_{t_0}(t_0) = U(1/t_0) w'_+$$
\noi imposed in Proposition 6.3 it follows that 
$$\exp [- i \varphi_{t_0} (t_0) ] \ U^*(1/t_0) \ w_{t_0}(t_0) = \exp [- i \varphi '_{t_0} (t_0)]
\ U^*(1/t_0) \ w'_{t_0}(t_0)$$
\noi and therefore by Lemma 7.2, $(w_{t_0}, \varphi_{t_0})$ and $(w'_{t_0} , \varphi '_{t_0})$ are
gauge equivalent, namely satisfy
$$\exp [- i \varphi_{t_0} (t)] \ U^*(1/t) \ w_{t_0} (t) = \exp [- i \varphi '_{t_0} (t)] \ U^*(1/t)
\ w'_{t_0} (t) \eqno(7.38)$$
\noi for all $t$ for which they are defined.\par

We now take the limit $t_0 \to \infty$ for fixed $t$ in (7.38). By Proposition 6.4, part (2)
supplemented with similar estimates on $\varphi_{t_0}$ and $\varphi '_{t_0}$, for fixed $t$,
$(w_{t_0}, \varphi_{t_0})$ and $(w'_{t_0}, \varphi '_{t_0})$ converge respectively to $(w, \varphi
)$ and $(w', \varphi ')$ in norm in $H^{k-1} \oplus Y^{\ell - 1}$. By an easy application of Lemma
3.5, it follows therefrom that one can take the limit $t_0 \to \infty$ in (7.38), thereby
obtaining (7.19), so that $(w, \varphi )$ and $(w', \varphi ')$ are gauge equivalent. \par

The last statement of Proposition 7.3 is a repetition of Proposition 7.2, part (3). \par \nobreak 
\hfill $\sq$ \\

We can now define the wave operator for $u$. We recall from the heuristic discussion in Section 2
that we want to exploit the operator $W$ defined in Definition 7.1, reconstruct $u$ through the
map $\Phi$ defined by (7.18), and eliminate the arbitrariness in $\varphi_{02}$ by fixing some
initial condition for it, namely $\varphi_{02}(1) = 0$, thereby purporting to ensure the
injectivity of the wave operator without restricting its range. That program is implemented by the
following definition and proposition. \\

\noi {\bf Definition 7.4.} We define the wave operator $\Omega$ as the map 
$$\Omega : u_+ \to u = (\Phi \circ W) (Fu_+, 0) \eqno(7.39)$$
\noi from $FH^{k+1}$ to ${\cal X}^k([T, \infty))$ for some $T$, $1 \leq T < \infty$, where $k$ is the
first element of an admissible pair, and $W$, $\Phi$ are defined by Definition 7.1 and by (7.18).
\par

The fact that $\Omega$ acts between the spaces indicated follows from Proposition 6.4 and from
Lemma 7.1. The value of $T$ depends on $u_+$ and can be taken according to (6.32) with $b = C(1 -
\gamma )^{-1} a^2$, namely
$$T = C \left ( (1 - \gamma )^{-2} + (2 \gamma - 1)^{-2} \right ) \left ( |Fu_+|_{k+1} \vee 1
\right )^{2/\gamma } \quad . \eqno(7.40)$$

\noi {\bf Proposition 7.4.} \par

{\it (1) The map $\Omega$ is injective.}  \par
{\it (2) ${\cal R}(\Omega ) = {\cal R} (\Phi \circ W)$.} \\

\noi {\bf Proof. Part (1).} Let $u = \Omega (u_+) = \Omega (u'_+)$ and let $(w, \varphi ) =
W(Fu_+, 0)$, $(w', \varphi ') = W(Fu'_+, 0)$, so that $u= \Phi (w, \varphi ) = \Phi (w', \varphi
')$, namely $(w, \varphi )$ and $(w', \varphi ')$ are gauge equivalent. By Proposition 7.2, part
(3), also $(Fu_+, 0)$ and $(Fu'_+, 0)$ are gauge equivalent in the sense of Definition 7.3, and
therefore $u_+ = u'_+$. \\

\noi {\bf Part (2).} Let $u = (\Phi \circ W ) (w_+ , \varphi_{02} (1)) \in {\cal R} (\Phi \circ W)$.
Then by Proposition 7.3, also $u = (\Phi \circ W) (w'_+, 0)$ where $w'_+ = w_+ \exp [- i \varphi_{02}
(1)]$, and $w'_+ \in H^{k+1}$ by Lemma 3.5. Therefore $u = \Omega (F^* w'_+) \in {\cal R}
(\Omega )$. \par \nobreak
\hfill $\sq$ \\

Note in particular that Proposition 7.4 part (2) means that we have not restricted the range of
the wave operators from $W$ to $\Omega$ by arbitrarily imposing $\varphi_{02} (1) = 0$. \par

We now collect all the available information on the solutions of the original equation (1.1) so
far constructed, namely existence through the previous definition of $\Omega$, some partial
form of uniqueness coming from Proposition 5.2, and asymptotic decay estimates coming from
Propositions 5.7 and 6.4. In order to state the result we need the phases $\varphi_{02}(t)$ and
$\varphi_0(t)$ defined now by (7.13) (7.12) with $w_+ = Fu_+$ and $\varphi_{02}(1) = 0$, namely
$$\varphi_{02} (t) = (1 - \gamma )^{-1} \ (t^{1 - \gamma} - 1) \ g_0(Fu_+) \eqno(7.41)$$
$$\varphi_0 (t) = \varphi_{02}(t) - \int_t^{\infty} dt' \ t'^{-\gamma} \left ( g_0 (FMu_+) -
g_0 (Fu_+) \right ) \quad . \eqno(7.42)$$
\noi (We recall that $M = M(t)$ is defined by (2.5) and satisfies (2.11)). The main result of
this paper can now be stated as follows. \\

\noi {\bf Proposition 7.5.} {\it Let $n \geq 3$, $0 < \mu \leq n - 2$, $1/2 < \gamma < 1$ and
let $(k, \ell )$ be an admissible pair. Let $u_+ \in F H^{k+1}$ and define $\varphi_{02}(t)$ and
$\varphi_0(t)$ by (7.41) and (7.42). Then} \par {\it (1) There exists a unique solution $u \in {\cal
X}^k([T, \infty ))$ of the equation (1.1) which can be represented as}
$$u = MD \exp (- i \varphi ) U^*(1/t) w \eqno(7.43)$$
\noi {\it where $(w, \varphi )$ is a solution of the system (7.1) (7.2) such that $(w,
\widetilde{\varphi}) \in ({\cal C} \cap L^{\infty})([T, \infty ), H^k \oplus Y^{\ell})$ and such
that}  $$|w(t) - F \ u_+|_{k-1} \ t^{1 - \gamma} \to 0 \eqno(7.44)$$ $$\parallel \varphi (t) -
\varphi_{02} (t) ; Y^{\ell - 1} \parallel \to 0 \eqno(7.45)$$ \noi {\it when $t \to \infty$. The
time $T$ depends on $\gamma$ and $u_+$ and can be taken in the form (7.40).} \par
{\it (2) The solution $u$ is obtained as $u = \Omega (u_+)$ where the map $\Omega$ is defined in
Definition 7.4. The map $\Omega$ is injective.} \par
{\it (3) The map $\Omega$ is continuous on the bounded sets of $FH^{k+1}$ from the norm topology in
$F H^{k-1}$ for $u_+$ to the norm topology in ${\cal X}^{k-1}(J)$ and to the weak-$*$ topology in
${\cal X}^k(J)$ for $u$ for any compact interval $J \subset [T, \infty )$, and to the weak
topology in $FH^k$ pointwise in $t$.} \par {\it (4) The solution $u$ satisfies the following
estimates for $t \geq T$~:} $$\parallel <J(t)>^k \left ( \exp [i \varphi_{02} (t, x/t) ] u(t) - M(t)
\ D(t) \ Fu_+ \right ) \parallel_2 \ \leq A(2 \gamma - 1)^{-1} \ t^{1/2 - \gamma} \quad ,
\eqno(7.46)$$  $$\parallel <J(t)>^k
\left ( \exp [i \varphi_{0} (t, x/t) ] u(t) - U(t) \ u_+ \right ) \parallel_2 \ \leq A(2 \gamma -
1)^{-1} \ t^{1 -2 \gamma} \quad , \eqno(7.47)$$              
\noi {\it for some constant $A$ depending on $\gamma$ and $u_+$ and bounded in $\gamma$ for fixed
$u_+$ and $\gamma$ away from 1.} \\

\noi {\bf Proof. Parts (1) and (2).} All the results except uniqueness follow from Proposition 6.4
supplemented with the reconstruction of $\varphi$ and from the subsequent definition of $\Omega$.
In particular (7.43) is essentially (7.18) and the injectivity of $\Omega$ is Proposition 7.4, part
(1). \par

Uniqueness is an immediate consequence of Proposition 5.2, given the asymptotic behaviour of $(w,
\varphi )$ that follows from Proposition 6.4. \\

\noi {\bf Part (3)} follows from Proposition 6.4, part (3) and from Lemma 7.1. \\

\noi {\bf Part 4.} From Proposition 5.7, part (1), esp. (5.69) and from Proposition 6.4, part (1),
esp. (6.30), supplemented by similar estimates on $\varphi - \varphi_0$ and on $\varphi -
\varphi_{02}$ easily obtained from (7.8) (7.11), it follows that
$$|w(t) - Fu_+|_k \leq A_0 \ t^{- \gamma} \eqno(7.48)$$
$$\parallel \varphi (t) - \varphi_0(t); Y^{\ell} \parallel \ \leq A_0 (2 \gamma - 1)^{-1} \ t^{1 - 2
\gamma} \eqno(7.49)$$
$$\parallel \varphi (t) - \varphi_{02} (t) ; Y^{\ell} \parallel \ \leq A_0(2 \gamma - 1)^{-1} \
t^{1/2 - \gamma} \eqno(7.50)$$
\noi for some constant $A_0$ of the type stated for $A$. From the definition (7.15) of $J$, from
the commutation relation (7.16), from (7.43) and from Lemma 3.5, we obtain
$$\parallel <J(t)>^k \left ( \exp [i \varphi_{02}(t, x/t) ] u(t) - M\ D\ F \ u_+ \right )
\parallel_2 \ = \Big | \exp [i (\varphi_{02} - \varphi )] U^*(1/t) w - Fu_+ \Big |_k$$
$$\leq \left \{ \parallel \varphi_{02} - \varphi \parallel_{\infty} \ + |\nabla (\varphi_{02} -
\varphi )|_{\ell - 1} \left ( 1 + |\nabla (\varphi_{02} - \varphi )|_{\ell - 1} \right )^{k-1}
\right \} |w|_k$$
$$+ |w - Fu_+|_k + |(U^*(1/t) - 1) F u_+ |_k \eqno(7.51)$$
\noi which yields immediately (7.46) by the use of (7.48) (7.50) and (6.9). Similarly
$$\parallel <J(t)>^k \left ( \exp [i \varphi_{0}(t, x/t) ] u(t) - U(t) \ u_+ \right )
\parallel_2 \ = \Big | \exp [i (\varphi_{0} - \varphi )] U^*(1/t) w - U^*(1/t) \ Fu_+ \Big |_k$$
$$\leq \left \{ \parallel \varphi_{0} - \varphi \parallel_{\infty} \ + |\nabla (\varphi_{0} -
\varphi )|_{\ell - 1} \left ( 1 + |\nabla (\varphi_{0} - \varphi )|_{\ell - 1} \right )^{k-1}
\right \} |w|_k + |w - Fu_+|_k  \eqno(7.52)$$
\noi from which (7.47) follows by the use of (7.48) (7.49). \par \nobreak
\hfill $\sq$ \\

\noi {\bf Remark 7.1.} The uniqueness statement in Proposition 7.5 is rather restrictive because it
requires the representation of $u$ by (7.43). It would be more satisfactory to have uniquess under
assumptions bearing directly on $u$, for instance under (7.46). However (7.46) seems insufficient to
derive the asymptotic conditions on $w$ and $\varphi$ separately which are required in Proposition
5.2. \\

\noi {\bf Remark 7.2.} The estimate (7.46) states that $u$ behaves asymptotically as expected, namely
$$u(t) \sim \exp \left [ - i \varphi_{02}(t, x/t) + i x^2 (2t)^{-1} \right ] (it)^{-n/2} (Fu_+)(x/t)
\quad . \eqno(7.53)$$
\noi In order to have the strongest possible statement however, one has to shift the phase factor
to $u$ before taking derivatives as contained in $<J(t)>^k$. On the other hand, this is not
necessary if one wants only asymptotic estimates in $L^r$. In fact one has the following corollary.
\\ 
 
\noi {\bf Corollary 7.1.} {\it Let $u$ be the solution of (1.1) obtained in Proposition 7.5. Let
$r$ satisfy $0 \leq \delta (r) \leq k \wedge n/2$, $\delta (r) < n/2$ if $k = n/2$. Then $u$
satisfies the following estimates}
$$\parallel u(t) - \exp [- i \varphi_{02} (t, x/t)] M \ D \ Fu_+\parallel_r \ \leq A(2 \gamma -
1)^{-1} \ t^{- \delta (r) + 1/2 - \gamma} \quad , \eqno(7.54)$$
 $$\parallel u(t) - \exp [- i \varphi_{0} (t, x/t)] U(t) \ u_+\parallel_r \ \leq A(2 \gamma -
1)^{-1} \ t^{- \delta (r) + 1 - 2\gamma} \quad . \eqno(7.55)$$

\noi {\bf Proof.} An immediate consequence of (7.46) (7.47) and of the inequality 
$$\begin{array}{ll} 
\parallel f \parallel_r \ = t^{- \delta (r)} \parallel D^*M^*f\parallel_r & \leq C \ t^{-
\delta (r)} \parallel < \nabla >^k D^* M^* f\parallel_2 \\ & \\
&= C\ t^{- \delta (r)} \parallel <J(t)>^k f\parallel_2 \end{array}$$
\noi which follows from Lemma 3.1 and from the commutation relation (7.16). \par \nobreak
\hfill $\sq$ \\

\noi {\bf Remark 7.3.} As already mentioned before, the phase $\varphi_0$ used in Proposition 5.7
produces a better asymptotic approximation of $u$ than the simpler phase $\varphi_{02}$ used in
Proposition 6.4. This shows up in the estimate behaving as $t^{1 - 2 \gamma}$ in the RHS of (7.47)
as compared with $t^{1/2 - \gamma}$ in the RHS of (7.46). In the same spirit, we have compared $u$
with the solution $U(t) u_+$ of the free Schr\"odinger equation in (7.47), and with the standard
asymptotic form $MDFu_+$ thereof in (7.46), obtained in dropping the second $M$ in $U = MDFM$. The
difference between the two is
$$|F(M(t) - 1) u_+|_k \leq t^{-1/2} |Fu_+|_{k+1}$$
\noi and shows up if one uses the latter, in the form of the last norm in (7.51). Since $t^{-1/2} <
t^{1/2 - \gamma}$ for $\gamma < 1$, that term is smaller than the other terms in the RHS of
(7.51), thereby preserving the estimate (7.46). This justifies the use of the simple explicit form
$MDFu_+$ in that case. On the other hand $t^{-1/2} > t^{1-2 \gamma}$ for $\gamma > 3/4$, and we
have therefore preferred to keep the more precise $U(t) u_+$ in (7.47) in order to preserve the
$t^{1 - 2 \gamma}$ decay in the RHS for all $\gamma < 1$. \par

So far we have constructed local solutions of the equation (1.1) in a neighborhood of in\-fi\-ni\-ty
associated with given asymptotic states $u_+$ and defined the local wave operator at infinity
$\Omega$. In order to complete the construction of the standard wave operators, it remains only to
extend the previous solutions from a neighborhood of infinity by using the results on the global
Cauchy problem at finite times. This can be done with the help of the following result which is
essentially contained in \cite{25r}. \\

\noi {\bf Proposition 7.6.} {\it Let $k$ be a positive integer and let $0 < \mu < 2$. Then the
Cauchy problem for the equation (1.1) with initial data $u(t_0) = u_0$ such that $<J(t_0)>^k u_0
\in L^2$ at some initial time $t_0 \geq 1$ is globally well posed in ${\cal X}^k([1, \infty ))$,
namely the local solutions of Proposition 7.1 can be extended to $[1, \infty )$.} \\

\noi {\bf Proof.} A minor variation of Proposition 2.1 part (1) in \cite{25r}. \par \nobreak
\hfill $\sq$ \\

We can now define the standard wave operator $\Omega_1$ for the equation (1.1). \\

\noi {\bf Definition 7.5.} Under the assumptions of Proposition 7.5 supplemented with $\mu < 2$,
we define the wave operator $\Omega_1$ as the map $\Omega_1 : u_+ \to u(1)$ where $u$ is the
solution of the equation (1.1) obtained by continuing $\Omega (u_+)$ down to $t = 1$ with the
help of Proposition 7.6. \\

>From Propositions 7.5 and 7.6 it follows that $\Omega_1$ is an injective map from $FH^{k+1}$ to
the space
$$K^k = \left \{ u : \exp (-ix^2/2) u \in H^k \right \}$$
\noi and that $\Omega_1$ satisfies continuity properties easily obtained from Proposition 7.5,
part (3). Since all the interesting information is already contained in that proposition, we
refrain from a more formal statement.

\newpage
\section*{Appendix A}
{\bf Sketch of the proof of Lemma 3.2} \par

We discuss only the case where $p < \infty$, the case $p = \infty$ being obtainable from the
previous one by a limiting procedure. The proof proceeds by first writing a regularized equation for
a suitable approximant of the function $|u|^p$, then by proving an analogue of the estimate (3.13)
for that function and finally by removing the regularisation. The method used is known and applied
in \cite{5r} to the equation (3.12) with $\eta = 0$. For this reason, even though the assumptions on
$u$ and $v$ used there are slightly different from those of Proposition 3.2, we refer to \cite{5r}
for estimating the second and third term in the RHS of (3.12). Here below we continue with the
analysis of the first ($\eta$ dependent) term. \par
Let $\varphi \in {\cal C}_0^{\infty} ({I\hskip-1truemm R}^n)$ denote a regularizing sequence of
functions of the space variable, namely a sequence converging to the $\delta$ function. Let $\phi$
denote the operator of convolution with $\varphi$, i.e. $\phi f = \varphi * f$, let $\varepsilon >
0$, let $u_1 \equiv \phi * u$, and let $w \equiv (|u_1|^2 + \varepsilon^2)^{1/2}$. The function $w$
is uniformly bounded away from zero, is ${\cal C}^{\infty}$ in the space variable and satisfies
$$\lim_{\varphi \to \delta} \ \lim_{\varepsilon \to 0} \int w(t, x)^p dx = \ \parallel u(t)
\parallel_p^p \quad .$$
\noi We compute $\partial_t w^p$ by using (3.12), thereby obtaining 
$$\begin{array}{ll}
\partial_t w^p &= 2w^{p-2} \ {\rm Re} (\bar{u}_1 \  \partial_t \ u_1)\\ & \\
&= p \left \{ \eta \ J_1 +J_2 +J_3 \right \}
\end{array}$$
\noi where
$$\begin{array}{l}
J_1 = w^{p-2} \ {\rm Re} \ \bar{u}_1 \ \Delta u_1 \\ \\
J_2 = w^{p-2} \ {\rm Re} \ \bar{u}_1 \ \phi \ \nabla (uv) \\ \\
J_3 = w^{p-2} \ {\rm Re} \ \bar{u}_1 \phi \ h \quad .
\end{array}$$
\noi The terms $J_2$ and $J_3$ are treated as in \cite{5r}. For the completion of the proof it is
sufficient to show that the space integral of $J_1$ is negative, so that $J_1$ does not
contribute to the estimates (3.13) (3.14). Repeated application of the Leibnitz rule yields
$$\begin{array}{ll}
J_1 &= - w^{p-2}  | \nabla u_1 |^2 + 1/2 \ w^{p-2} \left ( \nabla \cdot \nabla |u_1|^2 \right
) \\ & \\ &= - w^{p-2} |\nabla u_1|^2 - (1/4) (p - 2) w^{p-4} \left ( \nabla |u_1|^2 \right )^2 +
(1/2) \nabla \cdot \left \{ w^{p-2} \nabla |u_1|^2 \right \} \end{array}$$
\noi so that
$$J_1 \leq (1/2) \nabla \cdot \left \{ w^{p-2} \nabla |u_1|^2 \right \}$$
\noi which essentially implies that
$$\int J_1 \ dx \leq 0 \quad .$$
\noi Actually, in order to ensure integrability of $J_1$, the limit $\varepsilon \to 0$ should be
taken before performing the integration over the space variables up to infinity. \par \nobreak
\hfill $\sq$

\end{document}